\documentclass[12pt]{amsart}

\usepackage{amssymb}
\usepackage{graphicx}
\usepackage{ifthen}
\usepackage{mathrsfs}
\usepackage{pict2e}
\usepackage{xargs}
\usepackage{xspace}

\newcommand{\definedsymbol}[1]{$#1$}
\newcommand{\definedterm}[1]{\emph{#1}}

\makeatletter
\@addtoreset{section}{part}
\makeatother

\newcommand{\addblanklinetocontents}{\addtocontents{toc}{\vspace*{12pt}}}

\newcommand{\absolutelycontinuous}[1][]{\ll_{#1}}

\newcommand{\action}{\curvearrowright}

\newcommand{\ASL}[2]{#2^2 \semidirectproduct \mathrm{SL}_{#1}(#2)}
\newcommand{\Bairespace}[1][]{
  \ifthenelse{\equal{#1}{}}{\functions{\N}{\N}}{\functions{#1}{\N}}
}
\newcommand{\Bairetree}[1][]{
  \ifthenelse{\equal{#1}{}}{\functions{<\N}{\N}}{\functions{#1}{\N}}
}
\newcommand{\ball}[2]{\calB(#1, #2)}
\newcommand{\bermetric}[2]{e_{#2}}
\newcommand{\Borelfunctions}[2]{L(#1, #2)}
\newcommand{\calA}{\mathscr{A}}
\newcommand{\calB}{\mathscr{B}}
\newcommand{\calD}{\mathscr{D}}
\newcommand{\calE}{\mathscr{E}}
\newcommand{\calF}{\mathscr{F}}
\newcommand{\calL}{\mathscr{L}}
\newcommand{\calN}{\mathcal{N}}

\newcommand{\Cantorspace}[1][\N]{
  \functions{#1}{2}
}
\newcommand{\Cantortree}[1][]{
  \ifthenelse{\equal{#1}{}}{\functions{<\N}{2}}{\functions{#1}{2}}
}
\newcommand{\cardinality}[1]{|#1|}

\newcommand{\closedinterval}[2]{[#1, #2]}
\newcommand{\closedopeninterval}[2]{[#1, #2)}
\newcommand{\codedfunction}[1]{\phi_{#1}}
\newcommand{\composition}{\circ}
\newcommandx{\concatenation}[2][1 =, 2 =]{
  \ifthenelse{\equal{#1}{}}{{}^\smallfrown}{
    \ifthenelse{\equal{#2}{}}{\bigoplus #1}{\bigoplus_{#1} #2}
  }
}
\newcommand{\continuousfunctions}[2]{C(#1, #2)}
\newcommand{\continuum}{2^{\aleph_0}}
\newcommand{\convergesto}{\rightarrow}
\newcommand{\cost}[2]{C_{#1}(#2)}
\newcommand{\countingmeasure}{\mu_c}

\newcommand{\determinant}[1]{\mathrm{det}(#1)}
\newcommand{\diagonal}[1]{\Delta(#1)}
\newcommand{\differenceset}[2]{D(#1, #2)}
\newcommandx{\disjointunion}[2][1 =, 2 =]{
  \ifthenelse{\equal{#1}{}}{\sqcup}{
    \ifthenelse{\equal{#2}{}}{\bigsqcup #1}{{\bigsqcup_{#1} #2}}
  }
}
\newcommandx{\disjunction}[2][1 =, 2 =]{
  \ifthenelse{\equal{#1}{}}{\vee}{
    \ifthenelse{\equal{#2}{}}{\bigvee #1}{\bigvee_{#1} #2}
   }
}
\newcommand{\domain}[1]{\mathrm{dom}(#1)}

\newcommand{\equivalenceclass}[2]{[#1]_{#2}}
\newcommand{\ergodic}[2]{\mathcal{E}_{#2}}

\newcommand{\ergodicinvariant}[2]{\mathcal{EI}_{#2}}
\newcommand{\ergodicquasiinvariant}[2]{\mathcal{EQ}_{#2}}
\newcommand{\extendedby}{\sqsubseteq}
\newcommand{\extensions}[1]{\calN_{#1}}
\newcommand{\Ezero}{\mathbb{E}_0}

\newcommand{\forconullmany}[1]{\forall^*_{#1}}
\newcommand{\from}{\colon}
\newcommandx{\functions}[3][2 = undefined]{
  \ifthenelse{\equal{#2}{undefined}}{#3^{#1}}{
    \ifthenelse{\equal{#2}{measured}}{L(#1, #3)}{L(#1, #2, #3)}
  }
}
\newcommand{\goesto}{\rightarrow}
\newcommand{\graph}[1]{\mathrm{graph}(#1)}
\newcommand{\heightcorrection}[1]{\raisebox{0pt}[0pt][0pt]{#1}}
\newcommand{\homomorphisms}[3]{
  \mathrm{Hom}(#1, #2, #3)
}
\newcommand{\horizontalsection}[2]{#1^{#2}}
\newcommand{\hyperfinite}[2]{\mathcal{H}_{#2}}
\newcommand{\id}{\mathrm{id}}
\newcommand{\image}[2]{#1(#2)}
\newcommandx{\intersection}[2][1 =, 2 =]{
  \ifthenelse{\equal{#1}{}}{\cap}{
    \ifthenelse{\equal{#2}{}}{\bigcap #1}{{\bigcap_{#1} #2}}
  }
}
\newcommand{\invariant}[2]{\mathcal{I}_{#2}}
\newcommand{\inverse}[1]{#1^{-1}}

\newcommand{\mathand}{\text{ and }}
\newcommand{\mathcomma}{\text{, }}
\newcommand{\mathcommaand}{\text{, and }}
\newcommand{\mathor}{\text{ or }}
\renewcommand{\matrix}[4]{
  \left(
    \begin{smallmatrix}
      #1 & #2 \\
      #3 & #4
    \end{smallmatrix}
  \right)
}
\newcommand{\measureequivalence}[1][]{\sim_{#1}}

\newcommand{\N}{\mathbb{N}}
\newcommand{\openclosedinterval}[2]{(#1, #2]}
\newcommand{\openinterval}[2]{(#1, #2)}
\newcommand{\orbitequivalencerelation}[2]{E_{#1}^{#2}}
\newcommand{\pair}[2]{(#1, #2)}
\newcommand{\preimage}[2]{#1^{-1}(#2)}
\newcommand{\probabilitymeasures}[1]{P(#1)}
\newcommandx{\projection}[2][1 =, 2 =]{
  \ifthenelse{\equal{#1}{}}{\mathrm{proj}}{
    \ifthenelse{\equal{#2}{}}{\projection_{#1}}{
      \projection[#1](#2)
    }
  }
}
\newcommand{\pushforward}[2]{{#1}_* #2}
\newcommand{\Q}[1][]{
  \ifthenelse{\equal{#1}{}}{\mathbb{Q}}{\mathbb{Q}^{#1}}
}
\newcommand{\quasiinvariant}[2]{\mathcal{Q}_{#2}}
\newcommand{\R}[1][]{
  \ifthenelse{\equal{#1}{}}{\mathbb{R}}{\mathbb{R}^{#1}}
}
\renewcommand{\restriction}[2]{#1 \upharpoonright #2}
\newcommand{\reverseabsolutelycontinuous}[1][]{\gg_{#1}}
\newcommand{\Rplus}{\openinterval{0}{\infty}}
\newcommand{\saturation}[2]{[#1]_{#2}}
\newcommand{\semidirectproduct}{\rtimes}
\newcommandx{\sequence}[2][2 = undefined]{
  \ifthenelse{\equal{#2}{undefined}}{(#1)}{
    (#1)_{#2}
  }
}
\newcommandx{\set}[2][2 = undefined]{
  \ifthenelse{\equal{#2}{undefined}}{\{ #1 \}}{
    \{ #1 \suchthat #2 \}
  }
}
\newcommandx{\sets}[4][3 = undefined, 4 = undefined]{
  \ifthenelse{\equal{#3}{undefined}}{[#2]^{#1}}{
    \ifthenelse{\equal{#4}{undefined}}{[#2]^{#1}_{#3}}{[#2]^{#1}_{#3 / #4}}
  }
}
\newcommand{\setcomplement}[1]{\twiddle #1}
\newcommand{\sigmaclass}[1]{\sigma(#1)}
\newcommandx{\Sigmaclass}[2][1=,2=]{
  \ifthenelse{\equal{#2}{}}{\mathbf{\Sigma}_{#1}}{\mathbf{\Sigma}^{#1}_{#2}}
}

\newcommand{\singleton}[1]{\set{#1}}
\newcommand{\singletonsequence}[1]{(#1)}
\newcommand{\SL}[2]{\mathrm{SL}_{#1}(#2)}
\renewcommand{\square}[1]{I(#1)}

\newcommand{\suchthat}{\mid}
\newcommand{\symmetricdifference}{\mathrel{\triangle}}
\newcommand{\T}[1][]{
  \ifthenelse{\equal{#1}{}}{\mathbb{T}}{\mathbb{T}^{#1}}
}
\newcommand{\tailequivalencerelation}[1]{E_t(#1)}
\newcommand{\textexponent}[2]{$#1^{\mathrm{#2}}$}
\newcommand{\trace}[1]{\mathrm{tr}(#1)}

\newcommand{\twiddle}{\raisebox{1pt}{\scalebox{.75}{$\mathord{\sim}$}}}
\newcommand{\uniformmetric}[1]{d_{#1}}
\newcommandx{\union}[2][1 =, 2 =]{
  \ifthenelse{\equal{#1}{}}{\cup}{
    \ifthenelse{\equal{#2}{}}{\bigcup #1}{{\bigcup_{#1} #2}}
  }
}

\renewcommand{\vector}[3][]{
  \ifthenelse{\equal{#1}{}}{
    \left(
      \begin{smallmatrix}
        #2 \\
        #3
      \end{smallmatrix}
    \right)
  }{
    \begin{smallmatrix}
      #2 \\
      #3
    \end{smallmatrix}
  }
}
\newcommand{\verticalsection}[2]{#1_{#2}}
\newcommand{\Z}[1][]{
  \ifthenelse{\equal{#1}{}}{\mathbb{Z}}{\mathbb{Z}^{#1}}
}

\newcommand{\Adams}{Adams\xspace}
\newcommand{\Baire}{Baire\xspace}
\newcommand{\Banach}{Ban\-ach\xspace}
\newcommand{\Borel}{Bor\-el\xspace}
\newcommand{\Burgess}{Bur\-gess\xspace}

\newcommand{\Ditzen}{Dit\-zen\xspace}
\newcommand{\Dougherty}{Dough\-er\-ty\xspace}
\newcommand{\Dye}{Dye\xspace}
\newcommand{\Effros}{Eff\-ros\xspace}
\newcommand{\Feldman}{Feld\-man\xspace}
\newcommand{\Fubini}{Fu\-bi\-ni\xspace}
\newcommand{\Gaboriau}{Gab\-or\-i\-au\xspace}
\newcommand{\Galvin}{Gal\-vin\xspace}
\newcommand{\Glimm}{Glimm\xspace}
\newcommand{\Harrington}{Har\-ring\-ton\xspace}
\newcommand{\Hjorth}{Hjorth\xspace}
\newcommand{\Hopf}{Hopf\xspace}
\newcommand{\Jackson}{Jack\-son\xspace}

\newcommand{\Kechris}{Kech\-ris\xspace}
\newcommand{\Krieger}{Krie\-ger\xspace}
\newcommand{\Kuratowski}{Kur\-at\-ow\-ski\xspace}
\newcommand{\Lebesgue}{Leb\-esgue\xspace}
\newcommand{\Louveau}{Lou\-veau\xspace}
\newcommand{\Lusin}{Lu\-sin\xspace}
\newcommand{\Mauldin}{Maul\-din\xspace}

\newcommand{\Moore}{Moore\xspace}
\newcommand{\Mycielski}{My\-ciel\-ski\xspace}
\newcommand{\Nikodym}{Nik\-od\'{y}m\xspace}
\newcommand{\Novikov}{No\-vik\-ov\xspace}
\newcommand{\Polish}{Po\-lish\xspace}
\newcommand{\Radon}{Ra\-don\xspace}
\newcommand{\Rokhlin}{Rokh\-lin\xspace}
\newcommand{\Segal}{Se\-gal\xspace}
\newcommand{\Silver}{Sil\-ver\xspace}
\newcommand{\Slaman}{Sla\-man\xspace}
\newcommand{\Steel}{Steel\xspace}
\newcommand{\Ulam}{U\-lam\xspace}

\newcommand{\vonNeumann}{von Neu\-mann\xspace}
\newcommand{\Weiss}{Weiss\xspace}
\newcommand{\Woodin}{Wood\-in\xspace}

\newenvironment{lemmaproof}{
  
  \begin{proof}
}{\end{proof}}

\newenvironment{propositionproof}{
  
  \begin{proof}
}{\end{proof}}

\newenvironment{sublemmaproof}{
  
  \begin{proof}
}{\end{proof}}

\newenvironment{theoremproof}{
  
  \begin{proof}
}{\end{proof}}

\newtheorem{lemma}{Lemma}[section]
\newtheorem{proposition}[lemma]{Proposition}
\newtheorem{sublemma}[lemma]{Sublemma}
\newtheorem{theorem}[lemma]{Theorem}

\newtheorem{introconjecture}{Conjecture}
\newtheorem{introtheorem}[introconjecture]{Theorem}

\theoremstyle{definition}

\newtheorem{remark}[lemma]{Remark}
\newtheorem*{acknowledgments}{Acknowledgments}

\DeclareRobustCommand{\SkipTocEntry}[4]{}

\newcommand{\intropart}[1]{
  \addtocontents{toc}{\SkipTocEntry}
  \part*{#1}
}

\newcommand{\introsection}[1]{
  \addtocontents{toc}{\SkipTocEntry}
  \section*{#1}
}

\begin{document}


\begin{abstract}
  We show that every basis for the countable \Borel equivalence relations strictly
  above $\Ezero$ under measure reducibility is uncountable, thereby ruling out
  natural generalizations of the \Glimm-\Effros dichotomy. We also push many
  known results concerning the abstract structure of the measure reducibility
  hierarchy to its base, using arguments substantially simpler than those
  previously employed.
\end{abstract}

\author[C.T. Conley]{Clinton T. Conley}

\address{
  Clinton T. Conley \\
  Department of Mathematics \\
  Carnegie Mellon University \\
  Pittsburgh, PA 15213 \\
  USA
}

\email{clintonc@andrew.cmu.edu}

\urladdr{
  http://www.math.cmu.edu/math/faculty/Conley
}

\thanks{The first author was supported in part by NSF Grant DMS-1500906.}

\author[B.D. Miller]{Benjamin D. Miller}

 \address{
  Benjamin D. Miller \\
  Kurt G\"{o}del Research Center for Mathematical Logic \\
  University of Vienna \\
  W\"{a}hringer Stra{\ss}e 25 \\
  1090 Vienna \\
  Austria
 }

\email{benjamin.miller@univie.ac.at}

\urladdr{
  http://www.logic.univie.ac.at/benjamin.miller
}

\thanks{The second author was supported in part by DFG SFB Grant 878
  and FWF Grant P28153.}



\title[Measure reducibility]{Measure reducibility of countable Borel equivalence relations}

\maketitle

\intropart{Introduction}

Over the last few decades, the notion of \Borel reducibility of equivalence
relations has been used to identify obstacles of definability inherent in
classification problems throughout mathematics. While there are far too
many such applications to provide an exhaustive list here, a few notable
examples include the classifications of torsion-free abelian groups \cite
{Hjorth:TFA, AdamsKechris, Thomas:JAMS, Thomas:ICM}, ergodic
measure-preserving transformations \cite{Hjorth:Measure, ForemanWeiss,
ForemanRudolphWeiss:Paris, ForemanRudolphWeiss:Annals}, separable
\Banach spaces \cite{FerencziLouveauRosendal, Rosendal:Notices}, and
separable $C^*$-algebras \cite{FarahTomsTornquist:IMRN,
FarahTomsTornquist:RAM, Sabok}. In order to better understand such
results, one must obtain insight into the abstract structure of the
\Borel reducibility hierarchy. Unfortunately, this has turned out to be a very
difficult task.

The first of the two main lines of research into the abstract structure of the
\Borel reducibility hierarchy concerns its base. The first such result appeared
in \cite{Silver}, where it was shown that equality on $\R$ is the immediate
successor of equality on $\N$ within the co-analytic equivalence relations.
Building upon this and operator-algebraic work in \cite{Glimm, Effros},
it was shown in \cite{HarringtonKechrisLouveau}
that the relation $\Ezero$ on $\Cantorspace$, given by
$x \mathrel{\Ezero} y \iff \exists n \in \N \forall m \ge n \ x(m) = y(m)$,
is the immediate successor of equality on $\R$ within the \Borel equivalence relations.
Work in this direction stalled shortly thereafter, with \cite[Theorem 2]{KechrisLouveau}
ruling out further such results within the \Borel equivalence relations. However, the
question of whether there are further such results within the countable \Borel equivalence
relations remains open.

The first of the two main goals of this paper is to show that every basis for the countable
\Borel equivalence relations strictly above $\Ezero$ under measure reducibility is
uncountable.

The second of the two main lines of research into the abstract structure of the \Borel
reducibility hierarchy concerns exotic properties appearing beyond its base. The first
such result, due originally to \Woodin and later refined in \cite{LouveauVelickovic},
was the existence of uncountable families of pairwise incomparable \Borel equivalence
relations. However, the underlying arguments depended heavily upon \Baire category
techniques, and \cite[Theorem 6.2]{HjorthKechris} ensures that such an approach
cannot yield incomparability of countable \Borel equivalence relations.

This difficulty was eventually overcome in \cite{AdamsKechris}, yielding the
existence of uncountable families of pairwise incomparable countable \Borel
equivalence relations, in addition to myriad further results concerning the complexity
of the \Borel reducibility hierarchy. The arguments behind these theorems marked a
sharp departure from earlier approaches, relying upon sophisticated supperrigidity
machinery for actions of linear algebraic groups.

Soon thereafter, similar techniques were used in \cite{Adams, Thomas} to obtain many
striking new properties of the \Borel reducibility hierarchy, such as the existence of
countable \Borel equivalence relations $E$ to which the disjoint union of two copies
of $E$ is not \Borel reducible. While many of the underlying arguments were later simplified in 
\cite{HjorthKechris:Rigidity}, even these refinements depended upon complex
rigidity phenomena. And while the still simpler arguments of \cite{Hjorth:Antichain}
gave rise to pairwise incomparable treeable countable \Borel equivalence relations,
they still gave little sense of how far one must travel beyond the base of the
\Borel reducibility hierarchy before encountering such extraordinary behavior.

The second of the two main goals of this paper is to show that such phenomena appear
just beyond $\Ezero$ under measure reducibility.

We obtain our results by introducing a measureless notion of rigidity, which we establish
directly for the usual action of $\SL{2}{\Z}$ on $\T[2]$. In the presence of a measure, this yields
strong separability properties of the induced orbit equivalence relation. Many of our results follow
rather easily from the latter, while others require an additional graph-theoretic stratification
theorem, also established via elementary methods.

\introsection{Basic notions}

A set $X$ is \definedterm{countable} if there is an injection $\phi \from X \to \N$.
A sequence $\sequence{X_r}[r \in \R]$ of sets is \definedterm{increasing} if $X_r
\subseteq X_s$ for all real numbers $r \le s$.

Suppose that $X$ and $Y$ are standard \Borel spaces. We say that a sequence
$\sequence{x_y}[y \in Y]$ of points of $X$ is \definedterm{\Borel} if $\set{\pair{x_y}
{y}}[y \in Y]$ is a \Borel subset of $X \times Y$, and more generally, a sequence
$\sequence{X_y}[y \in Y]$ of subsets of $X$ is \definedterm{\Borel} if $\set{\pair{x}{y}}
[x \in X_y \mathand y \in Y]$ is a \Borel subset of $X \times Y$.

Suppose that $E$ is a \Borel equivalence relation on $X$. We say that $E$ is
\definedterm{aperiodic} if all of its classes are infinite, $E$ is \definedterm{countable} if
all of its classes are countable, and $E$ is \definedterm{finite} if all of its classes are finite.
A \definedterm{subequivalence relation} of $E$ is a subset of $E$ that is an equivalence
relation on $X$. The \definedterm{$E$-saturation} of a set $W \subseteq X$, or
\definedsymbol{\saturation{W}{E}}, is the smallest $E$-invariant set containing $W$. The
\definedterm{orbit equivalence relation} induced by an action of a group $\Gamma$ on
$X$ is the equivalence relation on $X$ given by $x \mathrel{\orbitequivalencerelation
{\Gamma}{X}} y \iff \exists \gamma \in \Gamma \ \gamma \cdot x = y$.

Suppose that $F$ is a \Borel equivalence relation on $Y$. A \definedterm
{homomorphism} from $E$ to $F$ is a function $\phi \from X \to Y$ sending
$E$-equivalent points to $F$-equivalent points, a \definedterm{reduction}
of $E$ to $F$ is a homomorphism sending $E$-inequivalent points to
$F$-inequivalent points, and an \definedterm{embedding} of $E$ into $F$ is an
injective reduction.

A \definedterm{graph} on $X$ is an irreflexive symmetric set $G \subseteq
X \times X$. A \definedterm{path} through $G$ is a sequence $\sequence
{x_i}[i \le n]$ with the property that $\forall i < n \ x_i \mathrel{G} x_{i+1}$,
in which case $n$ is the \definedterm{length} of the path. A graph is
\definedterm{acyclic} if there is at most one injective path between any
two points. We say that $G$ is a \definedterm{graphing} of $E$ if $E$ is the
smallest equivalence relation on $X$ containing $G$. When $G$ is
acyclic, we also say that $G$ is a \definedterm{treeing} of $E$. We say that
$E$ is \definedterm{treeable} if there is a \Borel treeing of $E$.

Suppose that $\mu$ is a \Borel measure on $X$. We say that $\mu$ is
\definedterm{$E$-ergodic} if every $E$-invariant \Borel set is $\mu$-null or
$\mu$-conull, $\mu$ is \definedterm{$E$-invariant} if $\mu(B) = \mu(\image
{T}{B})$ for all \Borel sets $B \subseteq X$ and all \Borel injections $T \from
B \to X$ whose graphs are contained in $E$, and $\mu$ is \definedterm
{$E$-quasi-invariant} if the $E$-saturation of every $\mu$-null set is $\mu$-null.

We say that $E$ is \definedterm{$\mu$-nowhere reducible} to $F$ if there is
no $\mu$-positive \Borel set $B \subseteq X$ for which $\restriction{E}{B}$
is \Borel reducible to $F$, $E$ is \definedterm{$\mu$-reducible} to $F$ if
there is a $\mu$-conull \Borel set $C \subseteq X$ for which $\restriction{E}
{C}$ is \Borel reducible to $F$, $E$ is \definedterm{invariant-measure reducible}
to $F$ if $\restriction{E}{B}$ is $\mu$-reducible to $F$ for every \Borel set $B
\subseteq X$ and every $(\restriction{E}{B})$-invariant \Borel probability
measure $\mu$ on $B$, and $E$ is \definedterm{measure reducible}
to $F$ if $E$ is $\mu$-reducible to $F$ for every \Borel probability measure
$\mu$ on $X$. The corresponding notions of \definedterm{invariant-measure
embeddability} and \definedterm{measure embeddability} are defined analogously.
It is straightforward to check that invariant-measure embeddability, measure
embeddability, and measure reducibility are transitive (and only marginally more
difficult to check that invariant-measure reducibility is transitive). 

We say that $E$ is \definedterm{hyperfinite} if it is a union of an increasing
sequence $\sequence{E_n}[n \in \N]$ of finite \Borel subequivalence relations,
$E$ is \definedterm{$\mu$-nowhere hyperfinite} if there is no $\mu$-positive \Borel
set $B \subseteq X$ for which $\restriction{E}{B}$ is hyperfinite, $E$ is \definedterm
{$\mu$-hyperfinite} if there is a $\mu$-conull \Borel set $C \subseteq X$ for
which $\restriction{E}{C}$ is hyperfinite, $E$ is \definedterm{invariant-measure
hyperfinite} if $\restriction{E}{B}$ is $\mu$-hyper\-finite for every \Borel set
$B \subseteq X$ and every $(\restriction{E}{B})$-invariant \Borel probability
measure $\mu$ on $B$, and $E$ is \definedterm{measure hyperfinite} if $E$ is
$\mu$-hyperfinite for every \Borel probability measure $\mu$ on $X$. As a
countable \Borel equivalence relation is hyperfinite if and only if it is \Borel
reducible to $\Ezero$ (see Theorem \ref
{preliminaries:hyperfiniteness:theorem:reduction}), it immediately follows that a countable
\Borel equivalence relation is invariant-measure hyperfinite if and only if it is
invariant-measure reducible to $\Ezero$, and measure hyperfinite if and only
if it is measure reducible to $\Ezero$.

\introsection{Bases}

A \definedterm{quasi-order} on $Q$ is a reflexive transitive binary relation $\le$
on $Q$. A \definedterm{basis} for $Q$ under $\le$ is a set $B \subseteq Q$
such that $\forall q \in Q \exists b \in B \ b \le q$.

Here we seek to elucidate the extent to which measure theory can shed light on the
structure of the \Borel reducibility hierarchy just beyond $\Ezero$. But given our
limited knowledge of the structure of the hierarchy, the appropriate meaning of ``just
beyond'' is not entirely clear. We will focus on properties that hold of some relation in
every basis for the non-measure-hyperfinite countable \Borel equivalence relations
under measure reducibility. One should first strive to understand the structure of such
bases, the original motivation for this paper.

\begin{introtheorem} \label{introduction:theorem:antibasis}
  Every basis for the non-measure-hyperfinite countable \Borel equivalence relations
  under measure reducibility is uncountable.
\end{introtheorem}

\introsection{Separability}

Although we will later give a somewhat different definition, for the sake of the introduction we will 
say that $F$ is \definedterm{projectively separable} if whenever $X$ is a standard \Borel space, 
$E$ is a countable \Borel equivalence relation on $X$, and $\mu$ is a \Borel probability measure on
$X$ for which $E$ is $\mu$-nowhere hyperfinite, there is a \Borel set $R \subseteq X \times
Y$, whose vertical sections are countable, such that $\mu(\set{x \in B}[\neg x \mathrel{R}
\phi(x)]) = 0$ for every \Borel set $B \subseteq X$ and every countable-to-one \Borel
homomorphism $\phi \from B \to Y$ from $\restriction{E}{B}$ to $F$. It is easy to see that
measure-hyperfinite \Borel equivalence relations are projectively separable.

Recall that $\SL{2}{\Z}$ is the group of all two-by-two matrices with integer entries and
determinant one. The natural action of $\SL{2}{\Z}$ on $\R[2]$ factors over $\Z[2]$ to
an action of $\SL{2}{\Z}$ on the quotient space $\T[2]$. It is well-known that the orbit
equivalence relation induced by this action is not measure hyperfinite, although it is
treeable (see Propositions \ref{preliminaries:SL2:proposition:nonhyperfinite} and 
\ref{preliminaries:SL2:proposition:treeable}). Our primary new tool here is the following.

\begin{introtheorem}
  The orbit equivalence relation induced by the action of $\SL{2}{\Z}$ on $\T[2]$ is
  projectively separable.
\end{introtheorem}

We obtain Theorem \ref{introduction:theorem:antibasis} by showing that if $E$ is a
non-measure-hyper\-finite projectively-separable treeable countable \Borel equivalence
relation, then every basis for the non-measure-hyperfinite \Borel subequivalence
relations of $E$ under measure reducibility is uncountable.

Ultimately, one would like to have the analogous result for bases for the
non-measure-hyperfinite countable \Borel equivalence relations measure
reducible to $E$. We show that $E$ is a counterexample if and only if it
is a countable disjoint union of successors of $\Ezero$ under
measure reducibility. While the existence of such successors remains open,
we show that if there are any at all, then there are uncountably many.

As projective separability and treeability are closed downward under \Borel reducibility, every basis for
the non-measure-hyper\-finite countable \Borel equivalence relations under measure
reducibility contains a relation whose restriction to some \Borel set is not measure hyperfinite,
but is projectively separable and treeable. In particular, if we wish to prove that every such
basis contains a relation whose restriction to some \Borel set has a given property, then it
is sufficient to show that the property holds of every non-measure-hyperfinite projectively-separable
treeable countable \Borel equivalence relation.

\introsection{Antichains}

The existence of a \Borel sequence $\sequence{E_r}[r \in \R]$ of pairwise non-measure-reducible
treeable countable \Borel equivalence relations was first established in \cite[Theorem 1.1]
{Hjorth:Antichain}. In light of the above observations, the following yields a simple new proof
of this result, while simultaneously pushing it to the base of the reducibility hierarchy.

\begin{introtheorem} \label{introduction:theorem:antichain}
  Suppose that $X$ is a standard \Borel space and $E$ is a non-measure-hyperfinite
  projectively-separable treeable countable \Borel equivalence relation on $X$. Then
  there is a \Borel sequence $\sequence{E_r}[r \in \R]$ of pairwise non-measure-reducible
  subequivalence relations of $E$.
\end{introtheorem}

As with our anti-basis theorem, one would like to have the analogous result
in which each $E_r$ is measure reducible to $E$, rather than contained in $E$.
We show that $E$ is a counterexample if and only if there is a finite
family $\calF$ of successors of $\Ezero$ under measure reducibility for which $E$ is a
countable disjoint union of countable \Borel equivalence relations measure
bi-reducible with those in $\calF$.

In particular, it follows that the existence of a sequence $\sequence{E_n}[n \in \N]$ of
pairwise non-measure-reducible countable \Borel equivalence relations measure
reducible to $E$ is equivalent to the existence of a \Borel sequence $\sequence{E_r}[r \in \R]$
of pairwise non-measure-reducible countable equivalence relations measure
reducible to $E$. Moreover, the nonexistence of such sequences implies the stronger fact that
every sequence $\sequence{E_n}[n \in \N]$ of countable \Borel equivalence relations
measure reducible to $E$ contains an infinite subsequence that is increasing under measure
reducibility.

\introsection{Complexity}
In \cite{AdamsKechris}, the existence of perfect families of pairwise incomparable countable
\Borel equivalence relations with distinguished ergodic \Borel probability measures was used to
establish a host of complexity results. We obtain simple new proofs of these results, while
simultaneously pushing them to the base of the reducibility hierarchy, by establishing the
following strengthening of Theorem \ref{introduction:theorem:antichain}.

\begin{introtheorem}
  Suppose that $X$ is a standard \Borel space and $E$ is a non-measure-hyperfinite
  projectively-separable treeable countable \Borel equivalence relation on $X$. Then
  there are \Borel sequences $\sequence{E_r}[r \in \R]$ of subequivalence relations of
  $E$ and $\sequence{\mu_r}[r \in \R]$ of \Borel probability measures on $X$ such that:
  \begin{enumerate}
    \item Each $\mu_r$ is $\verticalsection{E}{r}$-ergodic and $\verticalsection{E}{r}$-quasi-invariant.
    \item For all distinct $r, s \in \R$, the relation $\verticalsection{E}{r}$ is $\mu_r$-nowhere
      reducible to the relation $\verticalsection{E}{s}$.
  \end{enumerate}
\end{introtheorem}

While this result is somewhat technical, the complexity results of \cite{AdamsKechris} are
all obtained as abstract consequences of its conclusion.

Again, one would like to have the analog for which each $E_r$ is measure reducible to $E$,
rather than contained in $E$. We show that $E$ is a counterexample if and only if it is a
countable disjoint union of successors of $\Ezero$ under measure reducibility.

\introsection{Products}
The existence of non-measure-hyperfinite treeable countable \Borel equivalence relations
which do not measure reduce every treeable countable \Borel equivalence relation was
originally established in \cite[Theorem 1.6]{Hjorth}. Identify \definedsymbol{E \times F} with the
equivalence relation on $X \times Y$ given by $\pair{x_1}{y_1} \mathrel{(E \times F)}
\pair{x_2}{y_2} \iff (x_1 \mathrel{E} x_2 \mathand y_1 \mathrel{F} y_2)$, and let
\definedsymbol{\diagonal{X}} denote the diagonal on $X \times X$. The following
yields a simple new proof of the aforementioned result.

\begin{introtheorem}
  Suppose that $X$ is a standard \Borel space and $E$ is a non-measure-hyperfinite
  projectively-separable countable \Borel equivalence relation on $X$. Then $E \times
  \diagonal{\R}$ is not measure reducible to a \Borel subequivalence relation of $E$.
\end{introtheorem}

In \cite[Theorem 3.3a]{Thomas}, the rigidity results behind \cite{AdamsKechris} were used to
establish the existence of countable \Borel equivalence relations $E$ with the property that for
no $n \in \N$ is $E \times \diagonal{n+1}$ measure reducible to $E \times \diagonal{n}$. While
there are non-measure-hyperfinite projectively-separable countable \Borel equivalence relations
that do not have this property, in light of our observations on bases, the following yields a
simple new proof of this result, while simultaneously pushing it to the base of the reducibility
hierarchy.

\begin{introtheorem}
  Suppose that $X$ is a standard \Borel space and $E$ is a non-measure-hyperfinite
  projectively-separable countable \Borel equivalence relation on $X$. Then there is a \Borel
  set $B \subseteq X$ such that for no $n \in \N$ is $(\restriction{E}{B}) \times
  \diagonal{n+1}$ measure reducible to $(\restriction{E}{B}) \times \diagonal{n}$.
\end{introtheorem}

\introsection{Containment versus reducibility}
In \cite{Adams}, the rigidity results behind \cite{AdamsKechris} were used to establish the existence
of countable \Borel equivalence relations $E \subseteq F$ on the same space such that $E$ is not measure
reducible to $F$. This was strengthened by the proof of \cite[Theorem 1.1]{Hjorth:Antichain}, which
actually provided an increasing \Borel sequence $\sequence{E_r}[r \in \R]$ of pairwise
non-measure-reducible treeable countable equivalence relations on the same space. In light of our
observations on bases, the following yields a simple new proof of this fact, while
simultaneously pushing it to the base of the reducibility hierarchy.

\begin{introtheorem}
  Suppose that $X$ is a standard \Borel space and $E$ is a non-measure-hyperfinite
  projectively-separable treeable countable \Borel equivalence relation on $X$. Then
  there is an increasing \Borel sequence $\sequence{E_r}[r \in \R]$ of pairwise
  non-measure-reducible subequivalence relations of $E$.
\end{introtheorem}

\introsection{Embeddability versus reducibility}
In \cite[Theorem 3.3b]{Thomas}, the rigidity results behind \cite{AdamsKechris} were used to
establish the existence of aperiodic countable \Borel equivalence relations $E$ and $F$ for
which $E$ is \Borel reducible to $F$, but $E$ is not measure embeddable into $F$. In fact,
such examples were produced with $E = F \times I(2)$, where $I(X) = X \times X$.

If $E$ is invariant-measure hyperfinite and $F$ is aperiodic and countable, then $E$ is measure
reducible to $F$ if and only if $E$ is measure embeddable into $F$ (see Proposition
\ref{embeddability:proposition:sufficientcondition:measurable}). In particular, if
$E$ is aperiodic and invariant-measure hyperfinite, then $E \times \square{\N}$ is measure
embeddable into $E$. In light of our observations on bases, the following yields a simple
new proof of the aforementioned result, while simultaneously pushing it to the base of the
reducibility hierarchy.

\begin{introtheorem}
  Suppose that $X$ is a standard \Borel space and $E$ is an aperiodic
  non-invariant-measure-hyperfinite projectively-separable treeable countable \Borel
  equivalence relation on $X$. Then there is an aperiodic \Borel subequivalence
  relation $F$ of $E$ with the property that for no $n \in \N$ is $F \times I(n+1)$
  measure embeddable into $F \times I(n)$.
\end{introtheorem}

\introsection{Refinements}

We have taken great care to state our results in forms which make both the theorems and the
underlying arguments as clear as possible. Nevertheless, by utilizing several additional ideas,
one can obtain many generalizations and strengthenings.

In particular, by establishing analogs of our results for orbit equivalence relations induced by free
\Borel actions of countable discrete non-abelian free groups, one can rule out strong dynamical
forms of the \vonNeumann conjecture, while simultaneously providing an elementary proof of the
existence of continuum-many pairwise incomparable such relations, as found, for example, in \cite
{GaboriauPopa}. Moreover, as the notion of comparison we consider is far weaker than those
typically appearing in ergodic theory, our results are correspondingly stronger.

One can also obtain similar results for substantial weakenings of measure reducibility, as well as
for broader classes of equivalence relations. We plan to explore such developments in future
papers.

\tableofcontents

\part{Preliminaries}

We assume familiarity with the basic results and terminology of descriptive set theory, as
found in \cite{Kechris}. We provide here all additional standard definitions
and previously known results utilized throughout the paper. Although we mainly give references to the
relevant arguments, we provide proofs when they are particularly short or difficult to find in
the literature. For the sake of simplicity, we assume the axiom of choice throughout. However, with
only one slight exception (see \S\ref{bases}), our results go through with only the axiom of
dependent choice.

\section{Borel equivalence relations} \label{preliminaries:Borelequivalencerelations}

A \definedterm{partial transversal} of an equivalence relation is a subset of its domain
intersecting every equivalence class in at most one point.

\begin{theorem}[\Silver] \label{preliminaries:Borelequivalencerelations:theorem:perfect}
  Suppose that $X$ is a \Polish space and $E$ is a co-analytic equivalence relation on $X$. Then
  exactly one of the following holds:
  \begin{enumerate}
    \item The relation $E$ has only countably-many classes.
    \item There is a continuous injection of $\Cantorspace$ into a partial transversal of $E$.
  \end{enumerate}
\end{theorem}

\begin{theoremproof}
  See \cite{Silver}.
\end{theoremproof}

We say that $E$ is \definedterm{smooth} if it is \Borel reducible to
equality on a standard \Borel space. The following fact ensures that
$\Ezero$ is the minimal non-smooth \Borel equivalence relation
under \Borel reducibility.

\begin{theorem}[\Harrington-\Kechris-\Louveau]
  \label{preliminaries:Borelequivalencerelations:proposition:E0:embedding}
  Suppose that $X$ is a \Polish space and $E$ is a \Borel equivalence relation on $X$. Then
  exactly one of the following holds:
  \begin{enumerate}
    \item The relation $E$ is smooth.
    \item There is a continuous embedding of $\Ezero$ into $E$.
  \end{enumerate}
\end{theorem}

\begin{theoremproof}
  See \cite[Theorem 1.1]{HarringtonKechrisLouveau}.
\end{theoremproof}

\section{Countable Borel equivalence relations} \label
  {preliminaries:countableBorelequivalencerelations}

We begin by considering smoothness in the presence of countability.

\begin{proposition} \label{preliminaries:countableBorelequivalencerelations:proposition:smooth:finite}
  Suppose that $X$ is a standard \Borel space and $E$ is a finite \Borel equivalence relation on
  $X$. Then $E$ is smooth.
\end{proposition}

\begin{propositionproof}
  By the isomorphism theorem for standard \Borel spaces (see, for example, \cite[Theorem 15.6]
  {Kechris}), there is a \Borel linear ordering $\le$ of $X$. But then the \Lusin-\Novikov uniformization
  theorem (see, for example, \cite[Theorem 18.10]{Kechris}) ensures that the function $\phi \from X
  \to X$, given by $\phi(x) = \min_\le \equivalenceclass{x}{E}$, is a \Borel reduction of $E$ to equality.
\end{propositionproof}

\begin{remark} \label{preliminaries:countableBorelequivalencerelations:remark:smooth:countable}
  We say that a subset of $X$ is \definedterm{$E$-complete} if it intersects every $E$-class.
  A \definedterm{selector} for $E$ is a reduction of $E$ to equality on $X$ whose graph is contained
  in $E$, and a \definedterm{transversal} of $E$ is an $E$-complete partial transversal of $E$.
  Although the above argument actually yields the apparently stronger fact that every finite
  \Borel equivalence relation has a \Borel selector, the \Lusin-\Novikov uniformization theorem implies
  that if $E$ is countable, then smoothness, the existence of a \Borel selector, the existence of a \Borel
  transversal, and the existence of a partition $\sequence{B_n}[n \in \N]$ of $X$ into \Borel partial
  transversals are all equivalent. Moreover, in the special case that $E$ is aperiodic, they are also
  equivalent to the existence of a partition $\sequence{B_n}[n \in \N]$ of $X$ into \Borel
  transversals.
\end{remark}

\begin{proposition}
  \label{preliminaries:countableBorelequivalencerelations:proposition:smooth:closure}
  Suppose that $X$ is a standard \Borel space and $E$ is a smooth countable \Borel equivalence
  relation on $X$. Then every \Borel subequivalence relation of $E$ is smooth.
\end{proposition}

\begin{propositionproof}
  By Remark \ref{preliminaries:countableBorelequivalencerelations:remark:smooth:countable}, there
  is a partition of $X$ into countably-many \Borel partial transversals of $E$. As every partial
  transversal of $E$ is a partial transversal of all of its subequivalence relations, one more
  application of Remark \ref
  {preliminaries:countableBorelequivalencerelations:remark:smooth:countable} yields the desired
  result.
\end{propositionproof}

A function $I \from X \to X$ is an \definedterm{involution} if $I^2 = \id$.

\begin{theorem}[\Feldman-\Moore]
  \label{preliminaries:countableBorelequivalencerelations:theorem:involutions}
  Suppose that $X$ is a standard \Borel space and $R \subseteq X \times X$ is a reflexive symmetric
  \Borel set whose vertical sections are all countable. Then there are \Borel involutions $I_n \from X \to X$
  with the property that $R = \union[n \in \N][\graph{I_n}]$.
\end{theorem}

\begin{theoremproof}
  This follows from the proof of \cite[Theorem 1]{FeldmanMoore}.
\end{theoremproof}

The following can be viewed as generalizations of \Rokhlin's Lemma.

\begin{proposition}[\Slaman-\Steel]
  \label{preliminaries:countableBorelequivalencerelations:proposition:markers}
  Suppose that $X$ is a standard \Borel space and $E$ is an aperiodic countable \Borel
  equivalence relation on $X$. Then there is a decreasing sequence $\sequence{B_n}[n
  \in \N]$ of $E$-complete \Borel subsets of $X$ with empty intersection.
\end{proposition}

\begin{propositionproof}
  By the isomorphism theorem for standard \Borel spaces, we can assume that $X =
  \Cantorspace$. For each $n \in \N$ and $x \in \Cantorspace$, let $s_n(x)$ be the lexicographically
  least $s \in \functions{n}{2}$ for which $\extensions{s} \intersection \equivalenceclass
  {x}{E}$ is infinite. The \Lusin-\Novikov uniformization theorem ensures that each of the
  functions $s_n$ is \Borel, thus so too is each of the sets $A_n = \set{x \in \Cantorspace}
  [s_n(x) \extendedby x]$. It follows that the sets $B_n = A_n \setminus \intersection[n \in \N][A_n]$
  are as desired.
\end{propositionproof}

\begin{proposition}
  \label{preliminaries:countableBorelequivalencerelations:proposition:markers:disjoint}
  Suppose that $X$ is a standard \Borel space and $E$ is an aperiodic countable \Borel
  equivalence relation on $X$. Then there is a sequence $\sequence{B_n}[n \in \N]$ of
  pairwise disjoint $E$-complete \Borel subsets of $X$.
\end{proposition}

\begin{propositionproof}
  By Proposition \ref{preliminaries:countableBorelequivalencerelations:proposition:markers},
  there is a decreasing sequence $\sequence{A_n}[n \in \N]$ of $E$-complete \Borel
  subsets of $X$ with empty intersection. Recursively define functions $k_n \from X \to \N$
  by first setting $k_0(x) = 0$, and then defining $k_{n+1}(x) = \min \set{k \in \N}[(A_{k_n(x)}
  \setminus A_k) \intersection \equivalenceclass{x}{E} \neq \emptyset]$. The
  \Lusin-\Novikov uniformization theorem ensures that these functions are \Borel, so the
  sets $B_n = \set{x \in X}[x \in A_{k_n(x)} \setminus A_{k_{n+1}(x)}]$ are as desired.
\end{propositionproof}

\section{Hyperfiniteness} \label{preliminaries:hyperfiniteness}

We begin with the most basic properties of hyperfiniteness.

\begin{proposition}[\Dougherty-\Jackson-\Kechris] \label{preliminaries:hyperfiniteness:proposition:closure:unions}
  Suppose that $X$ is a standard \Borel space and $E$ is a \Borel equivalence relation
  on $X$. Then the family of \Borel sets on which $E$ is hyperfinite is closed under countable
  unions.
\end{proposition}

\begin{propositionproof}
  See, for example, \cite[Proposition 5.2]{DoughertyJacksonKechris}.
\end{propositionproof}

\begin{proposition}[\Dougherty-\Jackson-\Kechris]
  \label{preliminaries:hyperfiniteness:proposition:closure:homomorphisms}
  The family of hyperfinite \Borel equivalence relations is closed downward under
  countable-to-one \Borel homomorphism.
\end{proposition}

\begin{propositionproof}
  This follows, for example, from \cite[Proposition 5.2]{DoughertyJacksonKechris}.
\end{propositionproof}

\begin{proposition}[\Jackson-\Kechris-\Louveau] \label
  {preliminaries:hyperfiniteness:proposition:aperiodic}
  Suppose that $X$ is a standard \Borel space and $E$ is an aperiodic countable \Borel
  equivalence relation on $X$. Then there is an aperiodic hyperfinite \Borel subequivalence
  relation $F$ of $E$.
\end{proposition}

\begin{propositionproof}
  See, for example, \cite[Lemma 3.25]{JacksonKechrisLouveau}.
\end{propositionproof}

We say that a countable discrete group $\Gamma$ is \definedterm{hyperfinite} if whenever $X$ is
a standard \Borel space and $\Gamma \action X$ is a \Borel action, the induced orbit equivalence
relation $\orbitequivalencerelation{\Gamma}{X}$ is hyperfinite.

\begin{proposition}[\Slaman-\Steel, \Weiss]
  \label{preliminaries:hyperfiniteness:proposition:Z}
  The group $\Z$ is hyperfinite.
\end{proposition}

\begin{propositionproof}
  See, for example, \cite[Lemma 1]{SlamanSteel}.
\end{propositionproof}

We say that $E$ is \definedterm{hypersmooth} if it is the union of an increasing
sequence $\sequence{E_n}[n \in \N]$ of smooth \Borel subequivalence relations.

\begin{theorem}[\Dougherty-\Jackson-\Kechris]
  \label{preliminaries:hyperfiniteness:theorem:hypersmooth}
  Suppose that $X$ is a standard \Borel space and $E$ is a hypersmooth countable \Borel
  equivalence relation on $X$. Then $E$ is hyperfinite.
\end{theorem}

\begin{theoremproof}
  See, for example, the beginning of \cite[\S8]{DoughertyJacksonKechris}.
\end{theoremproof}

The \definedterm{tail equivalence relation} induced by a function $T \from X \to X$
is the equivalence relation on $X$ given by
\begin{equation*}
  x \mathrel{\tailequivalencerelation{T}} y \iff \exists m, n \in \N \ T^m(x) = T^n(y).
\end{equation*}

\begin{theorem}[\Dougherty-\Jackson-\Kechris]
  \label{preliminaries:hyperfiniteness:theorem:tailrelation}
  Suppose that $X$ is a standard \Borel space and $T \from X \to X$ is \Borel. Then
  $\tailequivalencerelation{T}$ is hypersmooth.
\end{theorem}

\begin{theoremproof}
  See, for example, \cite[Theorem 8.1]{DoughertyJacksonKechris}.
\end{theoremproof}

We now mention several further facts concerning reducibility.

\begin{theorem}[\Dougherty-\Jackson-\Kechris]
  \label{preliminaries:hyperfiniteness:theorem:embedding}
  All hyperfinite \Borel equivalence relations on standard \Borel spaces are \Borel
  embeddable into all non-smooth \Borel equivalence relations on standard \Borel
  spaces.
\end{theorem}

\begin{theoremproof}
  This follows from Theorem \ref{preliminaries:Borelequivalencerelations:proposition:E0:embedding}
  and \cite[Theorem 1]{DoughertyJacksonKechris}.
\end{theoremproof}

\begin{theorem}[\Dougherty-\Jackson-\Kechris]
  \label{preliminaries:hyperfiniteness:theorem:reduction}
  Suppose that $X$ is a standard \Borel space and $E$ is a countable \Borel equivalence
  relation on $X$. Then the following are equivalent:
  \begin{enumerate}
    \item The relation $E$ is hyperfinite.
    \item The relation $E$ is \Borel reducible to $\Ezero$.
  \end{enumerate}
\end{theorem}

\begin{theoremproof}
  To see $(1) \implies (2)$, note that $\Ezero$ is non-smooth, and appeal to 
  Theorem \ref{preliminaries:hyperfiniteness:theorem:embedding}.
  To see $(2) \implies (1)$, note that $\Ezero$ is hyperfinite, so Proposition \ref
  {preliminaries:hyperfiniteness:proposition:closure:homomorphisms} ensures that so too is every
  countable \Borel equivalence relation \Borel reducible to $\Ezero$.
\end{theoremproof}

\begin{theorem}[\Dougherty-\Jackson-\Kechris] \label{preliminaries:hyperfiniteness:theorem:linearity}
  All hyperfinite \Borel eq\-uivalence relations on standard \Borel spaces are comparable under
  \Borel reducibility.
\end{theorem}

\begin{theoremproof}
  As all standard \Borel spaces are comparable under \Borel embeddability, and Remark \ref
  {preliminaries:countableBorelequivalencerelations:remark:smooth:countable} implies that
  smooth countable \Borel equivalence relations have \Borel transversals, it follows from the
  \Lusin-\Novikov uniformization theorem that all smooth countable \Borel equivalence
  relations are comparable under \Borel reducibility. But the desired result then follows from
  Theorem \ref{preliminaries:hyperfiniteness:theorem:embedding}.
\end{theoremproof}

\begin{proposition} \label{preliminaries:hyperfiniteness:proposition:smoothavoidance}
  Suppose that $X$ is a standard \Borel space and $E$ is a non-smooth countable \Borel
  equivalence relation on $X$. Then there is a \Borel reduction $\pi \from X \to X$ of $E$ to
  $E$ such that $E$ is non-smooth off of $\saturation{\image{\pi}{X}}{E}$.
\end{proposition}

\begin{propositionproof}
  By Theorem \ref{preliminaries:Borelequivalencerelations:proposition:E0:embedding}, it is sufficient
  to establish the proposition for $\Ezero$. Towards this end, observe that the function $\pi \from
  \Cantorspace \to \Cantorspace$, given by
  \begin{equation*}
    \pi(x)(n) =
      \begin{cases}
        x(m) & \text{if $n = 2m$, and} \\
        0 & \text{if $n$ is odd,}
      \end{cases}
  \end{equation*}
  is as desired.
\end{propositionproof}

\section{Treeability} \label{preliminaries:treeability}

Here we note the analog of Proposition \ref{preliminaries:hyperfiniteness:proposition:closure:homomorphisms}
for treeability.

\begin{proposition}[\Jackson-\Kechris-\Louveau] \label{preliminaries:treeability:proposition:closure}
  The family of treeable countable \Borel equivalence relations is closed downward under
  count\-able-to-one \Borel homomorphism.
\end{proposition}

\begin{propositionproof}
  See \cite[Proposition 3.3]{JacksonKechrisLouveau}.
\end{propositionproof}

\section{Measures} \label{preliminaries:measures}

Following \cite[\S17]{Kechris}, we use \definedsymbol{\probabilitymeasures{X}} to
denote the standard \Borel space of all \Borel probability measures on $X$, and
when $X$ is a \Polish space, we use the same notation to denote the \Polish
space of all \Borel probability measures on $X$. Two \Borel
measures $\mu$ and $\nu$ are \definedterm{orthogonal} if there is a \Borel
set which is $\mu$-null and $\nu$-conull.

\begin{theorem}[\Burgess-\Mauldin] \label{preliminaries:measures:theorem:perfect}
  Suppose that $X$ is a standard \Borel space and $A \subseteq \probabilitymeasures{X}$ is an
  uncountable analytic set of pairwise orthogonal measures. Then there are \Borel sequences
  $\sequence{B_c}[c \in \Cantorspace]$ of pairwise disjoint subsets of $X$ and $\sequence
  {\mu_c}[c \in \Cantorspace]$ of \Borel probability measures on $X$ in $A$ such that $\mu_c(B_c)
  = 1$ for all $c \in \Cantorspace$.
\end{theorem}

\begin{theoremproof}
  By the isomorphism theorem for standard \Borel spaces, we can assume that $X$ is
  a zero-dimensional \Polish space. Fix a countable clopen basis $\calA$ for $X$.
  By Theorem \ref{preliminaries:Borelequivalencerelations:theorem:perfect}, there is a
  continuous injection $\pi \from \Cantorspace \to A$. Fix real numbers $\epsilon_n > 0$
  such that $\sum_{n \in \N} \epsilon_n < \infty$, and appeal to the regularity of \Borel
  probability measures on \Polish spaces (see, for example, \cite[Theorem 17.10]{Kechris})
  to recursively obtain $k_n \in \N$, $\phi_n \from \Cantorspace[n] \to \Cantorspace
  [k_n]$, and $A_n \from \Cantorspace[n] \to \calA$ with the following properties:
  \begin{enumerate}
    \item $\forall n \in \N \forall s \in \Cantorspace[n] \ \phi_{n+1}(s \concatenation
      \singletonsequence {0}) \neq \phi_{n+1}(s \concatenation \singletonsequence{1})$.
    \item $\forall i < 2 \forall n \in \N \forall s \in \Cantorspace[n] \ \phi_n(s) \extendedby \phi_{n+1}
      (s \concatenation \singletonsequence{i})$.
    \item $\forall n \in \N \forall s, t \in \Cantorspace[n] \ ( s = t \iff A_n(s) \intersection A_n(t) \neq
      \emptyset )$.
    \item $\forall n \in \N \forall s \in \Cantorspace[n] \forall \mu \in \image{\pi}{\extensions{\phi_n(s)}}
      \ \mu(A_n(s)) \ge 1 - \epsilon_n$.
  \end{enumerate}
  Define $\phi \from \Cantorspace \to \Cantorspace$ by $\phi(c) = \union[n \in \N][\phi_n
  (\restriction{c}{n})]$, and for each $c \in \Cantorspace$, define $B_c = \union[n \in \N]
  [{\intersection[m \ge n][A_n(\restriction{c}{n})]}]$ and $\mu_c = (\pi \composition \phi)(c)$.
\end{theoremproof}

We now describe a means of coding \Borel functions, modulo sets which
are null with respect to \Borel probability measures, which is uniform in both the function
and the measure in question. Let \definedsymbol{\continuousfunctions{X}{Y}} denote the
space of continuous functions from $X$ to $Y$ (see, for example, \cite[\S4.E]{Kechris}).
In order to keep our coding as transparent as possible, we will assume that $X$, $Y$,
$\continuousfunctions{X}{Y}$, and $\continuousfunctions{Y}{X}$ are \Polish, and that
every continuous partial function from $X$ to $Y$ has a continuous total extension. This
holds, for example, when $X = Y = \Cantorspace$.

\begin{proposition} \label{preliminaries:measures:proposition:evaluation:continuous}
  Suppose that $X$ is a compact \Polish space and $Y$ is a \Polish space. Then the function $\phi
  \from \continuousfunctions{X}{Y} \times X \to Y$ given by $\phi(f, x) = f(x)$ is continuous.
\end{proposition}

\begin{propositionproof}
  It is sufficient to show that if $U \subseteq Y$ is open and $\phi(f, x) \in U$, then there are open
  neighborhoods $V$ and $W$ of $f$ and $x$ such that $\image{\phi}{V \times W} \subseteq U$.
  Towards this end, fix a \Polish metric $d$ on $Y$ compatible with its underlying topology. As
  $U$ is open, there exists $\epsilon > 0$ such that $\ball{f(x)}{\epsilon} \subseteq U$. As $f$ is
  continuous, there is an open neighborhood $W$ of $x$ such that $\image{f}{W} \subseteq
  \ball{f(x)}{\epsilon/2}$. Fix an open neighborhood $V$ of $f$ such that $\forall g \in V \forall x \in X
  \ d(f(x), g(x)) \le \epsilon/2$. It only remains to note that if $\pair{g}{y} \in V \times W$, then $d(f(x),
  g(y)) \le d(f(x), f(y)) + \epsilon / 2 < \epsilon$, thus $g(y) \in U$.
\end{propositionproof}

We refer to elements $c$ of $\functions{\N}{\continuousfunctions{X}{Y}}$ as \definedterm{codes}
for measurable functions. Proposition \ref{preliminaries:measures:proposition:evaluation:continuous}
ensures that the sets
\begin{equation*}
  D_n = \set{\pair{c}{x} \in \functions{\N}{\continuousfunctions{X}{Y}} \times X}[\forall m \ge n \ c(m)
    (x) = c(n)(x)]
\end{equation*}
and $D = \union[n \in \N][D_n]$ are \Borel. We associate with each $c \in \functions{\N}
{\continuousfunctions{X}{Y}}$ the map $\codedfunction{c} \from \verticalsection{D}{c} \to Y$,
where $\codedfunction{c}(x)$ is the eventual value of $\sequence{c(n)(x)}[n \in \N]$.

\begin{proposition} \label{preliminaries:measures:proposition:evaluation:codes}
  Suppose that $X$ and $Y$ are standard \Borel spaces. Then the function $\phi \from
  D \to Y$ given by $\phi(c, x) = \codedfunction{c}(x)$ is \Borel.
\end{proposition}

\begin{propositionproof}
  As $\phi(c, x) = y \iff \exists n \in \N \forall m \ge n \ c(m)(x) = y$, the graph of $\phi$ is \Borel,
  so $\phi$ is \Borel (see, for example, \cite[Theorem 14.12]{Kechris}).
\end{propositionproof}

The \definedterm{push-forward} of a \Borel measure $\mu$ on $X$ through a \Borel function $\phi
\from X \to Y$ is given by $(\pushforward{\phi}{\mu})(B) = \mu(\preimage{\phi}{B})$.

\begin{proposition} \label{preliminaries:measures:proposition:pushforward}
  Suppose that $X$ and $Y$ are standard \Borel spaces. Then the function $\phi \from \set{\pair{c}
  {\mu} \in \functions{\N}{\continuousfunctions{X}{Y}} \times \probabilitymeasures{X}}
  [\mu(\verticalsection{D}{c}) = 1] \to \probabilitymeasures{Y}$ given by $\phi(c, \mu) = 
  \pushforward{(\codedfunction{c})}{\mu}$ is \Borel.
\end{proposition}

\begin{propositionproof}
  It is sufficient to show that if $B \subseteq Y$ is \Borel and $F \subseteq \R$ is of the form
  $\openclosedinterval{a}{b}$, where $a < b$ are in $\R$, then the intersection of the sets
  \begin{equation*}
    R = \set{\pair{c}{\mu} \in \functions{\N}{\continuousfunctions{X}{Y}} \times \probabilitymeasures
      {X}}[\mu(\verticalsection{D}{c}) = 1]
  \end{equation*}
  \centerline{and}
  \begin{equation*}
    S = \set{\pair{c}{\mu} \in \functions{\N}{\continuousfunctions{X}{Y}} \times \probabilitymeasures
      {X}}[\pushforward{(\codedfunction{c})}{\mu}(B) \in F]
  \end{equation*}
  is \Borel. But $R$ is clearly \Borel (see, for example, \cite[Theorem 17.25]{Kechris}), and
  to see that $S$ is \Borel, observe that $\pair{c}{\mu} \in S$ if and only if $\exists n \in \N \forall
  m \ge n \ \mu(\preimage{c(m)}{B} \intersection \verticalsection{(D_n)}{c}) \in F$.
\end{propositionproof}

\section{Measured equivalence relations} \label{preliminaries:measuredequivalencerelations}

Here we consider countable \Borel equivalence relations in the presence of measures.

\begin{proposition} \label{preliminaries:measuredequivalencerelationsproposition:null}
  Suppose that $X$ is a standard \Borel space, $E$ is a non-smooth \Borel equivalence relation
  on $X$, and $\mu$ is a \Borel probability measure on $X$. Then there is a $\mu$-null \Borel
  set on which $E$ is non-smooth.
\end{proposition}

\begin{propositionproof}
  By Theorem \ref{preliminaries:Borelequivalencerelations:proposition:E0:embedding}, there is a
  continuous embedding $\pi \from \Cantorspace \to X$ of $\Ezero$ into $E$. For each
  $c \in \Cantorspace$, the function $\pi_c \from \Cantorspace \to \Cantorspace$, given by
  \begin{equation*}
    \pi_c(d)(n) =
      \begin{cases}
        c(m) & \text{if $n = 2m$, and} \\
        d(m) & \text{if $n = 2m + 1$},
      \end{cases}
  \end{equation*}
  is a continuous embedding of $\Ezero$ into $\Ezero$. As the sets of the form $\image{\pi_c}
  {\Cantorspace}$ are pairwise disjoint, it follows that for all but countably many $c \in
  \Cantorspace$, the function $\pi \composition \pi_c$ is as desired.
\end{propositionproof}

Suppose that $\rho \from E \to \Rplus$ is a \definedterm{cocycle}, in the sense that
$\rho(x, z) = \rho(x, y)\rho(y, z)$ whenever $x \mathrel{E} y \mathrel{E} z$. For
each set $Y \subseteq \equivalenceclass{x}{E}$, define $\rho(Y, x) = \sum_{y \in
Y} \rho(y, x)$. We say that $Y$ is \definedterm{$\rho$-finite} or \definedterm
{$\rho$-infinite} according to whether $\rho(Y, x)$ is finite or infinite. Our assumption
that $\rho$ is a cocycle ensures that the $\rho$-finiteness of $Y$ does not depend on the choice
of $x \in \saturation{Y}{E}$. We say that $\rho$ is \definedterm{finite} if every equivalence
class of $E$ is $\rho$-finite, and $\rho$ is \definedterm{aperiodic} if every equivalence class of
$E$ is $\rho$-infinite. Given $Y, Z \subseteq \equivalenceclass{x}{E}$, define $\rho(Y, Z) = \rho(Y,
x) / \rho(Z, x)$. Again, our assumption that $\rho$ is a cocycle ensures that $\rho(Y, Z)$ does not
depend on the choice of $x \in \saturation{Y}{E}$.

\begin{proposition} \label{preliminaries:measuredequivalencerelationsproposition:periodic}
  Suppose that $X$ is a standard \Borel space, $E$ is a countable \Borel equivalence relation on
  $X$, and there is a finite \Borel cocycle $\rho \from E \to \Rplus$. Then $E$ is smooth.
\end{proposition}

\begin{propositionproof}
  See, for example, \cite[Proposition 2.1]{Miller:Measures}.
\end{propositionproof}

\begin{theorem}[\Ditzen] \label{preliminaries:measuredequivalencerelations:theorem:Borel}
  Suppose that $X$ is a standard \Borel space and $E$ is a countable \Borel equivalence relation
  on $X$. Then the set of $E$-ergodic $E$-quasi-invariant \Borel probability measures on $X$ is \Borel.
\end{theorem}

\begin{theoremproof}
  See \cite[Theorem 2 of Chapter 2]{Ditzen}.
\end{theoremproof}

We say that $\mu$ is \definedterm{$\rho$-invariant} if $\mu(\image{T}{B}) = \int_B \rho(T(x),
x) \ d\mu(x)$, for all \Borel sets $B \subseteq X$ and all \Borel injections $T \from B \to X$
whose graphs are contained in $E$.

\begin{proposition} \label{preliminaries:measuredequivalencerelationsproposition:cocycles}
  Suppose that $X$ is a standard \Borel space, $E$ is a countable \Borel equivalence relation on
  $X$, and $\mu$ is an $E$-quasi-invariant \Borel probability measure on $X$. Then there is a
  \Borel cocycle $\rho \from E \to \Rplus$ with respect to which $\mu$ is invariant.
\end{proposition}

\begin{propositionproof}
  See, for example, \cite[\S8]{KechrisMiller}.
\end{propositionproof}

We say that $E$ is \definedterm{$\mu$-nowhere smooth} if there is no
$\mu$-positive \Borel set $B \subseteq X$ for which $\restriction{E}{B}$ is smooth.

\begin{proposition} \label{preliminaries:measuredequivalencerelationsproposition:aperiodic}
  Suppose that $X$ is a standard \Borel space, $E$ is a countable \Borel equivalence relation on
  $X$, $\rho \from E \to \Rplus$ is an aperiodic \Borel cocycle, and $\mu$ is a $\rho$-invariant
  \Borel probability measure on $X$. Then $E$ is $\mu$-nowhere smooth.
\end{proposition}

\begin{propositionproof}
  See, for example, \cite[Proposition 2.1]{Miller:Measures}.
\end{propositionproof}

The following fact usually allows us to assume quasi-invariance.

\begin{proposition} \label{preliminaries:measuredequivalencerelationsproposition:quasiinvariant}
  Suppose that $X$ is a standard \Borel space, $E$ is a countable \Borel equivalence relation on
  $X$, and $\mu$ is a \Borel probability measure on $X$. Then there is an $E$-quasi-invariant
  \Borel probability measure $\nu$ on $X$ such that $\mu \absolutelycontinuous \nu$ and the two
  measures take the same values on all $E$-invariant \Borel sets.
\end{proposition}

\begin{propositionproof}
  By Theorem \ref{preliminaries:countableBorelequivalencerelations:theorem:involutions}, there is a
  countable group $\Gamma = \set{\gamma_n}[n \in \N]$ of \Borel automorphisms of $X$ whose
  induced orbit equivalence relation is $E$. Define $\nu = \sum_{n \in \N} \pushforward
  {(\gamma_n)}{\mu} / 2^{n+1}$.

  To see that $\nu$ is a \Borel probability measure, simply note that it is a convex combination of
  \Borel probability measures. Moreover, if $B \subseteq X$ is an $E$-invariant \Borel set,
  then $\mu(B) = (\pushforward{\gamma}{\mu})(B)$ for all $\gamma \in \Gamma$, thus $\nu(B) =
  \sum_{n \in \N} \mu(B) / 2^{n+1} = \mu(B)$. And if $N \subseteq X$ is a
  $\nu$-null \Borel set, then $\mu(N) \le \sum_{n \in \N} \pushforward{(\gamma_n)}{\mu}(N) = 0$,
  thus $\mu \absolutelycontinuous \nu$.
\end{propositionproof}

Note that any two $E$-ergodic $E$-quasi-invariant \Borel probability measures are either
orthogonal or equivalent; the following gives a sufficient condition to strengthen equivalence
to equality.

\begin{proposition} \label{preliminaries:measuredequivalencerelations:proposition:equality}
  Suppose that $X$ is a standard \Borel space, $E$ is a countable \Borel equivalence relation on
  $X$, $\rho \from E \to \Rplus$ is a \Borel cocycle, and $\mu \absolutelycontinuous \nu$ are
  $E$-ergodic $\rho$-invariant \Borel probability measures on $X$. Then $\mu = \nu$.
\end{proposition}

\begin{propositionproof}
  The \Radon-\Nikodym theorem (see, for example, \cite[\S17.A]{Kechris}) yields a
  \Borel function $\phi \from X \to \closedopeninterval{0}{\infty}$ such that $\mu(B) = \int \phi(x)
  \ d\nu(x)$ for all \Borel sets $B \subseteq X$. As $\mu(X) = \nu(X) = 1$, to see that $\mu =
  \nu$, it is sufficient to show that $\phi$ is constant on a $\mu$-conull \Borel set. Suppose,
  towards a contradiction, that there are $\mu$-positive \Borel sets $A, B \subseteq X$ with
  the property that $\forall x \in A \forall y \in B \ \phi(x) < \phi(y)$. As $E$ is countable, Theorem
  \ref{preliminaries:countableBorelequivalencerelations:theorem:involutions} yields a countable
  group $\Gamma$ of \Borel automorphisms of $X$ whose induced orbit equivalence relation
  is $E$. As $\mu$ is $E$-ergodic, there exists $\gamma \in \Gamma$ such that the set $A' =
  A \intersection \preimage{\gamma}{B}$ is $\mu$-positive, so
  \begin{equation*}
    \mu(\image{\gamma}{A'})
      = \int_{A'} \rho(\gamma \cdot x, x) \ d\mu(x) 
      = \int_{A'} \phi(x) \rho(\gamma \cdot x, x) \ d\nu(x)
  \end{equation*}
  \centerline{and}
  \begin{equation*}
      \mu(\image{\gamma}{A'})
        = \int_{\image{\gamma}{A'}} \phi(x) \ d\nu(x)
        = \int_{A'} \phi(\gamma \cdot x) \rho(\gamma \cdot x, x) \ d\nu(x),
  \end{equation*}
  the desired contradiction.
\end{propositionproof}

A \definedterm{\Borel disintegration} of a \Borel probability measure $\mu$ on $X$ through a \Borel
function $\phi \from X \to Y$ is a \Borel sequence $\sequence{\mu_y}[y \in Y]$ of \Borel probability
measures on $X$ with the property that $\mu = \int \mu_y \ d(\pushforward{\phi}{\mu})(y)$ and
$\mu_y(\preimage{\phi}{y}) = 1$ for all $y \in Y$. The existence of such sequences is noted, for
example, in \cite[Exercise 17.35]{Kechris}.

A \definedterm{\Borel ergodic decomposition} of a \Borel cocycle $\rho \from E \to \Rplus$ is a
\Borel sequence $\sequence{\mu_x}[x \in X]$ of \Borel probability measures on $X$ such that
$\mu_x = \mu_y$ for all $\pair{x}{y} \in E$, $\mu(\set{x \in X}[\mu = \mu_x]) = 1$ for all $E$-ergodic
$\rho$-invariant \Borel probability measures $\mu$, and $\mu = \int \mu_x \ d\mu(x)$ for all 
$\rho$-invariant \Borel probability measures $\mu$.
\begin{theorem}[\Ditzen]
  \label{preliminaries:measuredequivalencerelations:theorem:ergodicdecomposition}
  Suppose that $X$ is a standard \Borel space, $E$ is a \Borel equivalence relation on $X$, and
  $\rho \from E \to \Rplus$ is a \Borel cocyle. Then there is a \Borel ergodic decomposition of
  $\rho$.
\end{theorem}

\begin{theoremproof}
  See \cite[Theorem 6 of Chapter 2]{Ditzen}.
\end{theoremproof}

A \definedterm{compression} of $E$ is a \Borel injection $T \from X \to X$, whose graph is
contained in $E$, such that the complement of $\image{T}{X}$ is $E$-complete. We say that
$E$ is \definedterm{compressible} if there is a \Borel compression of $E$.

\begin{proposition} \label{preliminaries:measuredequivalencerelations:proposition:compression}
  Suppose that $X$ is a standard \Borel space, $E$ is a countable \Borel equivalence relation on
  $X$, and $B \subseteq X$ is a \Borel $E$-complete set for which $\restriction{E}{B}$ is
  compressible. Then there is a \Borel injection $T \from X \to B$ whose graph is contained in $E$.
\end{proposition}

\begin{propositionproof}
  Fix a \Borel compression $\phi \from B \to B$ of $\restriction{E}{B}$.
  The \Lusin-\Novikov uniformization theorem yields a \Borel function
  $\psi \from X \to B \setminus \image{\phi}{B}$ whose graph is contained
  in $E$, as well as a \Borel function $\xi \from X \to \N$ such that $\psi
  \times \xi$ is injective. Set $\pi(x) = \phi^{\xi(x)} \composition \psi(x)$.
\end{propositionproof}

We say that $E$ is \definedterm{$\mu$-nowhere compressible} if there is no
$\mu$-positive \Borel set $B \subseteq X$ for which $\restriction{E}{B}$ is
compressible.

\begin{theorem}[\Hopf] \label{preliminaries:measuredequivalencerelations:theorem:existence}
  Suppose that $X$ is a standard \Borel space, $E$ is a countable \Borel equivalence relation on
  $X$, and $\mu$ is an $E$-quasi-invariant $\sigma$-finite \Borel measure on $X$. If $E$ is
  $\mu$-nowhere compressible, then there is an $E$-invariant \Borel probability measure
  $\nu \measureequivalence \mu$.
\end{theorem}

\begin{theoremproof}
  See, for example, \cite[\S10]{Nadkarni}.
\end{theoremproof}

When $\mu$ is an $E$-invariant \Borel probability measure, the \definedterm{$\mu$-cost}
of a graphing $G$ of $E$ is given by
\begin{equation*}
  \cost{\mu}{G} = \frac{1}{2} \int \cardinality{\verticalsection{G}{x}} \ d\mu(x).
\end{equation*}
The \definedterm{$\mu$-cost} of $E$ is the infimum of the
costs of its \Borel graphings.

\begin{proposition}[\Gaboriau]
  \label{preliminaries:measuredequivalencerelations:proposition:costformula}
  Suppose that $X$ is a standard \Borel space, $E$ is a countable \Borel equivalence relation on
  $X$, $\mu$ is an $E$-invariant \Borel probability measure on $X$, $B \subseteq X$ is an
  $E$-complete \Borel set, and $\mu_B$ is the \Borel probability measure on $B$ given by
  $\mu_B(D) = \mu(D) / \mu(B)$. Then $\cost{\mu}{E} - 1 = \mu(B)(\cost{\mu_B}{\restriction{E}{B}} -
  1)$. In particular, it follows that $\cost{\mu}{E} \le \cost{\mu_B}{\restriction{E}{B}}$.
\end{proposition}

\begin{propositionproof}
  See, for example, \cite[Theorem 21.1]{KechrisMiller}.
\end{propositionproof}

\begin{proposition}[\Gaboriau]
  \label{preliminaries:measuredequivalencerelations:proposition:costbound}
  Suppose that $X$ is a standard \Borel space, $E$ is an aperiodic treeable countable \Borel
  equivalence relation on $X$, and $\mu$ is an $E$-invariant \Borel probability measure on $X$
  for which $E$ is not $\mu$-hyperfinite. Then $\cost{\mu}{E} > 1$.
\end{proposition}

\begin{propositionproof}
  See, for example, \cite[Corollary 27.12]{KechrisMiller}.
\end{propositionproof}

An $E$-ergodic measure $\mu$ is \definedterm{$\pair{E}{F}$-ergodic} if there is no
$\mu$-null-to-one \Borel homomorphism from $E$ to $F$.

\begin{proposition} \label{preliminaries:measuredequivalencerelations:proposition:E0ergodic}
  Suppose that $X$ is a standard \Borel space, $E$ is a countable \Borel equivalence relation on
  $X$, $\mu$ is an $\pair{E}{\Ezero}$-ergodic \Borel probability measure 
  on $X$, and $\sequence{E_n}[n \in \N]$ is an increasing sequence of countable \Borel equivalence
  relations on $X$ whose union is $E$. Then for all $\epsilon > 0$, there is a \Borel set $B
  \subseteq X$ of $\mu$-measure at least $1 - \epsilon$ on which $\mu$ is $E_n$-ergodic for all
  sufficiently large $n \in \N$.
\end{proposition}

\begin{propositionproof}
  See, for example, \cite[Proposition 2.2]{Miller:Hjorth}.
\end{propositionproof}

For the following, recall the definition of codes for measurable functions given just before
Proposition \ref{preliminaries:measures:proposition:evaluation:codes}.

\begin{proposition} \label{preliminaries:measuredequivalencerelations:proposition:analytic}
  Suppose that $X$ and $Y$ are compact \Polish spaces and $E$ and $F$ are countable \Borel
  equivalence relations on $X$ and $Y$. Then the set of pairs $\pair{c}{\mu} \in \functions{\N}
  {\continuousfunctions{X}{Y}} \times \probabilitymeasures{X}$ for which $\codedfunction{c}$ is a
  reduction of $E$ to $F$ on an $E$-invariant $\mu$-conull \Borel set is analytic.
\end{proposition}

\begin{propositionproof}
  By Theorem \ref{preliminaries:countableBorelequivalencerelations:theorem:involutions}, there are
  countable groups $\Gamma$ and $\Delta$ of \Borel automorphisms of $X$ and $Y$ whose
  induced orbit equivalence relations are $E$ and $F$. Then $\pair{c}{\mu}$ has the desired
  property if and only if there exists $d \in \functions{\N}{\continuousfunctions{Y}{X}}$ such that the
  following conditions hold:
  \begin{enumerate}
    \item $\mu(\verticalsection{D}{c}) = 1$.
    \item $\forconullmany{\mu} x \in X \forall \gamma \in \Gamma \ \codedfunction{c}(x) \mathrel{F}
      \codedfunction{c}(\gamma \cdot x)$.
    \item $\forconullmany{\mu} x \in X \forall \delta \in \Delta \ (\delta \cdot \codedfunction{c}(x) \in
      \verticalsection{D}{d} \implies x \mathrel{E} \codedfunction{d}(\delta \cdot
        \codedfunction{c}(x)))$.
    \item $\pushforward{(\codedfunction{c})}{\mu}(\verticalsection{D}{d}) = 1$.
  \end{enumerate}
  Clearly the sets determined by conditions (1) and (2) are \Borel, Proposition \ref
  {preliminaries:measures:proposition:evaluation:codes} ensures that the set
  determined by condition (3) is \Borel, and Proposition \ref
  {preliminaries:measures:proposition:pushforward} implies that the set determined
  by condition (4) is \Borel.
\end{propositionproof}

\section{Measure hyperfiniteness} \label{preliminaries:measurehyperfiniteness}

Here we consider connections between hyperfiniteness and measures.

\begin{proposition} \label{preliminaries:measurehyperfiniteness:proposition:nonamenable}
  Suppose that $\Gamma$ is a countable discrete non-amenable group, $X$ is a standard \Borel
  space, $\Gamma \action X$ is a free \Borel action, and $\mu$ is an $\orbitequivalencerelation
  {\Gamma}{X}$-invariant \Borel probability measure on $X$. Then the induced orbit equivalence
  relation $\orbitequivalencerelation{\Gamma}{X}$ is not $\mu$-hyperfinite.
\end{proposition}

\begin{propositionproof}
  See, for example, \cite[Proposition 2.5]{JacksonKechrisLouveau}.
\end{propositionproof}

\begin{theorem}[\Dye, \Krieger]
  \label{preliminaries:measurehyperfiniteness:theorem:increasingunion}
  Suppose that $X$ is a standard \Borel space, $\mu$ is a \Borel probability measure on $X$, and
  $\sequence{E_n}[n \in \N]$ is an increasing sequence of $\mu$-hyperfinite \Borel equivalence
  relations on $X$. Then the equivalence relation $E = \union[n \in \N][E_n]$ is also
  $\mu$-hyperfinite.
\end{theorem}

\begin{theoremproof}
  See, for example, \cite[Propositon 6.11]{KechrisMiller}.
\end{theoremproof}

Given equivalence relations $E$ and $F$ on $X$, define
\begin{equation*}
  \bermetric{X}{\mu}(E, F) = \mu(\set{x \in X}[\equivalenceclass{x}{E} \neq \equivalenceclass{x}{F}]).
\end{equation*}

\begin{proposition}
  Suppose that $X$ is a standard \Borel space and $\mu$ is a \Borel probability measure on $X$.
  Then $\bermetric{X}{\mu}$ is a complete pseudo-metric.
\end{proposition}

\begin{propositionproof}
  To see that $\bermetric{X}{\mu}$ is a pseudo-metric, it is sufficient to check the triangle inequality.
  Towards this end, suppose that $E_1$, $E_2$, and $E_3$ are \Borel equivalence relations on
  $X$, and observe that
  \begin{align*}
    \bermetric{X}{\mu}( & E_1, E_3) \\
      & = 1 - \mu(\set{x \in X}[\equivalenceclass{x}{E_1} = \equivalenceclass{x}{E_3}]) \\
      & \le 1 - \mu(\set{x \in X}[\equivalenceclass{x}{E_1} = \equivalenceclass{x}{E_2}] \intersection
        \set{x \in X}[\equivalenceclass{x}{E_2} = \equivalenceclass{x}{E_3}]) \\
      & = 1 + \mu(\set{x \in X}[\equivalenceclass{x}{E_1} = \equivalenceclass{x}{E_2}] \union \set
        {x \in X}[\equivalenceclass{x}{E_2} = \equivalenceclass{x}{E_3}]) \mathrel{-} \\
      & \mspace{21mu} (\mu(\set{x \in X}[\equivalenceclass{x}{E_1} = \equivalenceclass{x}{E_2}]) +
        \mu(\set{x \in X}[\equivalenceclass{x}{E_2} = \equivalenceclass{x}{E_3}])) \\
      & \le 2 - (\mu(\set{x \in X}[\equivalenceclass{x}{E_1} = \equivalenceclass{x}{E_2}]) +
        \mu(\set{x \in X}[\equivalenceclass{x}{E_2} = \equivalenceclass{x}{E_3}])) \\
      & = \bermetric{X}{\mu}(E_1, E_2) + \bermetric{X}{\mu}(E_2, E_3).
  \end{align*}

  To see that $\bermetric{X}{\mu}$ is complete, suppose that $\sequence{E_n}[n \in \N]$ is an
  $\bermetric{X}{\mu}$-Cauchy sequence, fix a sequence of real numbers $\epsilon_n > 0$ such
  that $\sum_{n \in \N} \epsilon_n < \infty$, and fix a strictly increasing sequence of natural numbers
  $k_n$ such that
  \begin{equation*}
    \forall n \in \N \forall i, j \ge k_n \ \bermetric{X}{\mu}(E_i, E_j) \le \epsilon_n.
  \end{equation*}
  Note that for all $n \in \N$, the set $Y_n = \set{x \in X}[\forall m \ge n \ \equivalenceclass{x}
  {E_{k_m}} = \equivalenceclass{x}{E_{k_n}}]$ has $\mu$-measure at least $1 - \sum_{m \ge n}
  \epsilon_m$. In particular, it follows that the set $Y = \union[n \in \N][Y_n]$ is $\mu$-conull.
  Letting $E$ denote the union of the diagonal on $X$ with the equivalence relations of the form
  $\restriction{E_{k_n}}{Y_n}$ for $n \in \N$, it follows that $E_{k_n} \convergesto_{\bermetric{X}
  {\mu}} E$ as $n \convergesto \infty$, thus $\bermetric{X}{\mu}$ is indeed complete.
\end{propositionproof}

It is not difficult to see that $\bermetric{X}{\mu}$ is not separable, even when restricted to the family
of \Borel equivalence relations on $X$ whose classes are all of cardinality two. In contrast, we
have the following.

\begin{proposition} \label{preliminaries:measurehyperfiniteness:proposition:finiteapproximation}
  Suppose that $X$ is a standard \Borel space and $E$ is a countable \Borel equivalence relation
  on $X$. Then there is a countable family $\calF$ of finite \Borel subequivalence relations of $E$
  such that for all \Borel probability measures $\mu$ on $X$, the family $\calF$ is $\bermetric{X}
  {\mu}$-dense in the set of all finite \Borel subequivalence relations of $E$.
\end{proposition}

\begin{propositionproof}
  Fix an enumeration $\sequence{U_n}[n \in \N]$ of a basis, closed under finite unions, for a
  \Polish topology generating the \Borel structure of $X$. By Theorem \ref
  {preliminaries:countableBorelequivalencerelations:theorem:involutions}, there is a sequence
  $\sequence{f_n}[n \in \N]$ of \Borel automorphisms of $X$ such that $E = \union[n \in \N]
  [\graph{f_n}]$. 
  
  For each $n \in \N$ and $s \in \Bairetree[n]$, let $X_s$ denote the \Borel set of $x \in X$ with the 
  property that whenever $i, j, k < n$, $x \in U_{s(i)} \intersection U_{s(j)}$, $y \in U_{s(k)}$,
  and $f_i(x) = f_k(y)$, there exists $\ell < n$ such that $y \in U_{s(\ell)}$ and $f_j(x) = f_\ell(y)$.
  Let $F_s$ denote the reflexive \Borel relation on $X$ in which distinct points $x$ and $y$ related if 
  there exist $i, j < n$ and $z \in U_{s(i)} \intersection U_{s(j)} \intersection X_s$ such that 
  $x = f_i(z)$ and $y = f_j(z)$.
  
  \begin{lemma}
    Each $F_s$ is an equivalence relation.
  \end{lemma}
  
  \begin{lemmaproof}
    As $F_s$ is clearly reflexive and symmetric, it is sufficient to show that it is transitive. 
    Towards this end, observe that if $x \mathrel{F_s} y \mathrel{F_s} z$ are pairwise distinct, 
    then there exist $i, j < n$ and $v \in U_{s(i)} \intersection U_{s(j)} \intersection X_s$ with 
    the property that $x = f_i(v)$ and $y = f_j(v)$, as well as $k, \ell < n$ and $w \in U_{s(k)} 
    \intersection U_{s(\ell)} \intersection X_s$ with the property that $y = f_k(w)$ and $z = 
    f_\ell(w)$. As $v \in X_s$, there exists $m < n$ with $w \in U_{s(m)}$ and $x = f_i(v) = 
    f_m(w)$, in which case the definition of $F_s$ ensures that $x \mathrel{F_s} z$. 
  \end{lemmaproof}
  
  To see that the family $\calF = \set{F_s}[s \in \Bairetree]$ is as desired, suppose that 
  $\epsilon > 0$, $F$ is a finite \Borel subequivalence relation of $E$, and $\mu$ is a 
  \Borel probability measure on $X$. Fix $n \in \N$ sufficiently large that the $\mu$-measure 
  of the set $Y = \set{x \in X}[\forall y, z \in \equivalenceclass{x}{F} \exists i < n \ f^i(y) = z]$ is 
  strictly greater than $1 - \epsilon$. Set $\delta = \mu(Y) - (1 - \epsilon)$, and define $Y_k = 
  \set{x \in X}[x \mathrel{F} f_k(x)]$ for all $k < n$. As \Borel probability measures on \Polish
  spaces are regular, there exists $s \in \Bairespace[n]$ with the property that the $\mu$-measure
  of the set
  \begin{equation*}
    Z_{i, j, k} = \set{x \in X}[(\inverse{f_i} \composition f_j)(x) \in U_{s(k)} \iff (\inverse{f_i} 
      \composition f_j)(x) \in Y_k)]
  \end{equation*}
  is at least $1 - \delta / n^3$, for all $i, j, k < n$.

  \begin{lemma} \label{lemma:ZcontainedinXs}
    The set $Z = Y \intersection \intersection[i, j, k < n][Z_{i, j, k}]$ is contained in $X_s$.
  \end{lemma}
  
  \begin{lemmaproof}
    We must show that if $i, j, k < n$, $z \in U_{s(i)} \intersection 
    U_{s(j)} \intersection Z$, $y \in U_{s(k)}$, and $f_i(z) = f_k(y)$, then there exists $\ell < n$ such that
    $y \in U_{s(\ell)}$ and $f_j(z) = f_\ell(y)$. Towards this end, note that $y = (\inverse{f_k} 
    \composition f_i)(z)$, so the fact that $z \in Z$ ensures that $z \in Y_i \intersection Y_j$ 
    and $y \in Y_k$. In particular, it follows that $f_j(z) \mathrel{F} z \mathrel{F} f_i(z) = 
    f_k(y) \mathrel{F} y$. The fact that $z \in Y$ then yields $\ell < n$ such that $f_j(z) = 
    f_\ell(y)$. As $y \in Y_\ell$, one more appeal to the fact that $z \in Z$ ensures that $y 
    \in U_{s(\ell)}$.
  \end{lemmaproof}

  \begin{lemma}
    Suppose that $z \in Z$. Then $\equivalenceclass{z}{F} = \equivalenceclass{z}{F_s}$.
  \end{lemma}
  
  \begin{lemmaproof}
    Suppose first that $x \in \equivalenceclass{z}{F}$. As $z \in Y$, there exist $i, j < n$
    such that $x = f_i(z)$ and $z = f_j(z)$. Then $z \in Y_i \intersection Y_j$, so the
    fact that $z \in Z$ ensures that $z \in U_{s(i)} \intersection U_{s(j)}$. As Lemma \ref
    {lemma:ZcontainedinXs} implies that $z \in X_s$, the definition of $F_s$ ensures that 
    $x \in \equivalenceclass{z}{F_s}$.
    
    Suppose now that $x \in \equivalenceclass{z}{F_s}$. The definition of $F_s$ then
    yields $i, j < n$ and $w \in U_{s(i)} \intersection U_{s(j)} \intersection X_s$ such that 
    $x = f_i(w)$ and $z = f_j(w)$. As $z \in Y$, there exists $\ell < n$ such that $z = f_\ell(z)$. As
    $z \in Y_\ell$, the fact that $z \in Z$ ensures that $z \in U_{s(\ell)}$, so the
    fact that $w \in X_s$ yields $k < n$ such that $z \in U_{s(k)}$ and $x = f_k(z)$. One more
    appeal to the fact that $z \in Z$ then ensures that $z \in Y_k$, in which case $x = f_k(z)
    \in \equivalenceclass{z}{F}$.
  \end{lemmaproof}

  As $\mu(Z) \ge 1 - \epsilon$, it follows that $\bermetric{X}{\mu}(F, F_s) \le \epsilon$.
\end{propositionproof}

We use $\hyperfinite{X}{E}$ denote the space of \Borel probability measures $\mu$
on $X$ with respect to which $E$ is $\mu$-hyperfinite. The following fact originally
appeared in \cite{Segal}.

\begin{theorem}[\Segal] \label{preliminaries:measurehyperfiniteness:theorem:uniform}
  Suppose that $X$ is a standard \Borel space and $E$ is a countable \Borel equivalence relation
  on $X$. Then there is a \Borel set $F \subseteq (\N \times (X \times X)) \times \probabilitymeasures
  {X}$ such that for all $\mu \in \probabilitymeasures{X}$, the following conditions hold:
  \begin{enumerate}
    \item The sets $\verticalsection{(\horizontalsection{F}{\mu})}{n}$ form an increasing 
      sequence of finite \Borel subequivalence relations of $E$.
    \item The set $\horizontalsection{B}{\mu} = \set{x \in X}[{\equivalenceclass{x}{E} \neq \union[n 
      \in \N][\equivalenceclass{x}{\verticalsection{(\horizontalsection{F}{\mu})}{n}}]}]$ does not contain 
      a $\mu$-positive \Borel subset on which $E$ is hyperfinite.
  \end{enumerate}
  In particular, it follows that $\hyperfinite{X}{E}$ is \Borel.
\end{theorem}

\begin{theoremproof}
  Fix real numbers $\epsilon_n > 0$ such that $\sum_{n \in \N} \epsilon_n < \infty$.
  By Proposition \ref{preliminaries:measurehyperfiniteness:proposition:finiteapproximation}, 
  there is a family $\calE = \set{E_k}[k \in \N]$ of finite \Borel subequivalence relations 
  of $E$ such that for all \Borel probability measures $\mu$ on $X$, the family $\calE$ 
  is $\bermetric{X}{\mu}$-dense in the set of all finite \Borel subequivalence relations of 
  $E$. Then the functions $m_n \from \probabilitymeasures{X} \to \closedinterval{0}{1}$ given by
  \begin{equation*}
    \textstyle
    m_n(\mu) = \sup_{k \in \N} \mu(\set{x \in X}[\forall i < n \ x \mathrel{E_k} f_i(x)])
  \end{equation*}
  are \Borel, as are the functions $k_n \from \probabilitymeasures{X} \to \N$ given by
  \begin{equation*}
    \textstyle
    k_n(\mu) = \min \set{k \in \N}[{\mu(\set{x \in X}[\forall i < n \ x \mathrel{E_k} f_i(x)]) > 
      m_n(\mu) - \epsilon_n}],
  \end{equation*}
  thus so too is the set $F \subseteq (\N \times (X \times X)) \times \probabilitymeasures
  {X}$ given by
  \begin{equation*}
    x \mathrel{\verticalsection{(\horizontalsection{F}{\mu})}{n}} y \iff \forall m \ge n \ x 
      \mathrel{E_{k_m(\mu)}} y.
  \end{equation*}
  
  To see that $F$ is as desired, suppose that $\mu \in \probabilitymeasures{X}$. As the 
  sets $\verticalsection{(\horizontalsection{F}{\mu})}{n} = \intersection[m \ge n][E_{k_m
  (\mu)}]$ form an increasing sequence of finite \Borel subequivalence relations of $E$, 
  it is enough to show that if $A \subseteq \horizontalsection{B}{\mu}$ is a \Borel set on 
  which $E$ is hyperfinite, then $\mu(A) = 0$. As $\horizontalsection{B}{\mu}$ is $E$-invariant
  and $E$ is countable, the \Lusin-\Novikov uniformization theorem and Proposition 
  \ref{preliminaries:hyperfiniteness:proposition:closure:homomorphisms} allow us to assume
  that $A$ is $E$-invariant.
  
  \begin{lemma}
    Suppose that $n \in \N$. Then
    \begin{equation*}
      \mu(\set{x \in A}[\forall i < n \ x \mathrel{E_{k_n(\mu)}} f_i(x)]) \ge \mu(A) - \epsilon_n.
    \end{equation*}
  \end{lemma}
  
  \begin{lemmaproof}
    As $E_{k_n(\mu)}$ is finite, Remark \ref
    {preliminaries:countableBorelequivalencerelations:remark:smooth:countable} ensures that it 
    has a \Borel transversal $C \subseteq X$ from which the quotient $X / E_{k_n(\mu)}$ 
    inherits a standard \Borel structure, and moreover, that the map associating each 
    $E_{k_n(\mu)}$-class with the unique point of $C$ it contains is a \Borel reduction of $E / 
    E_{k_n(\mu)}$ to $E$. Proposition \ref
    {preliminaries:hyperfiniteness:proposition:closure:homomorphisms} therefore implies that the 
    restriction of $E / E_{k_n(\mu)}$ to $A / E_{k_n(\mu)}$ is hyperfinite.
    
    Given $\epsilon > 0$, observe that all but finitely many relations $E'$ along any 
    sequence witnessing the hyperfiniteness of the restriction of $E / E_{k_n(\mu)}$ to $A / 
    E_{k_n(\mu)}$, when viewed as equivalence relations on $A$, satisfy the condition that
    $\mu(\set{x \in A}[\forall i < n \ x \mathrel{E'} f_i(x)]) > \mu(A) - \epsilon$. The $\bermetric
    {X}{\mu}$-density of $\calE$ therefore yields $k \in \N$ such that
    \begin{equation*}
      \mu(A) - \mu(\set{x \in A}[\forall i < n \ x \mathrel{E_{k_n(\mu)}} f_i(x)]) - \epsilon
    \end{equation*}
    is strictly less than
    \begin{equation*}
        \mu(\set{x \in X}[\forall i < n \ x \mathrel{E_k} f_i(x)]) - \mu(\set{x \in X}[\forall i < n
          \ x \mathrel{E_{k_n(\mu)}} f_i(x)]).
    \end{equation*}
    As the definition of $k_n(\mu)$ ensures that the latter quantity is itself strictly less than 
    $\epsilon_n$, it follows that
    \begin{equation*}
      \mu(\set{x \in A}[\forall i < n \ x \mathrel{E_{k_n(\mu)}} f_i(x)]) > \mu(A) - \epsilon_n - \epsilon,
    \end{equation*} 
    thus $\mu(\set{x \in A}
    [\forall i < n \ x \mathrel{E_{k_n(\mu)}} f_i(x)]) \ge \mu(A) - \epsilon_n$, as the former 
    inequality holds for all $\epsilon > 0$.
  \end{lemmaproof}
  
  Set $A' = \union[n \in \N][{\intersection[m \ge n][{\set{x \in A}[\forall i < m \ x \mathrel
  {E_{k_m(\mu)}} f_i(x)]}]}]$, and note that $\mu(A) = \mu(A')$, since $\sum_{n \in \N} \epsilon_n 
  < \infty$, thus $\mu(A) = 0$, since $A' \intersection \horizontalsection{B}{\mu} = \emptyset$. 
  
  As $\hyperfinite{X}{E} = \set{\mu \in \probabilitymeasures{X}}[\mu
  (\horizontalsection{B}{\mu}) = 0]$ and the \Lusin-\Novikov uniformization theorem 
  ensures that the set $B = \set{\pair{x}{\mu} \in X \times \probabilitymeasures{X}}[x \in 
  \horizontalsection{B}{\mu}]$ is \Borel, it follows that $\hyperfinite{X}{E}$ is \Borel as well.
\end{theoremproof}

\begin{proposition} \label{preliminaries:measurehyperfiniteness:proposition:characterization}
  Suppose that $X$ is a standard \Borel space, $E$ is a countable \Borel equivalence relation
  on $X$, $\rho \from E \to \Rplus$ is a \Borel cocycle, and there is a $\rho$-invariant \Borel
  probability measure $\mu$ on $X$ for which $E$ is not $\mu$-hyperfinite. Then there is such a
  measure which is also $E$-ergodic.
\end{proposition}

\begin{propositionproof}
  This follows from Theorems \ref
  {preliminaries:measuredequivalencerelations:theorem:ergodicdecomposition} and \ref
  {preliminaries:measurehyperfiniteness:theorem:uniform}.
\end{propositionproof}

We use \definedsymbol{\ergodic{X}{E}} to denote the family of all $E$-ergodic \Borel probability
measures on $X$, \definedsymbol{\quasiinvariant{X}{E}} to denote the family of all $E$-quasi-invariant 
\Borel probability measures on $X$, and \definedsymbol{\ergodicquasiinvariant{X}{E}} to denote
$\ergodic{X}{E} \intersection \quasiinvariant{X}{E}$.

\begin{theorem} \label{preliminaries:measurehyperfiniteness:theorem:E0}
  Suppose that $X$ is a standard \Borel space and $E$ is a countable \Borel equivalence
  relation on $X$. Then exactly one of the following holds:
  \begin{enumerate}
    \item The relation $E$ is measure hyperfinite.
    \item The set $\ergodicquasiinvariant{X}{E} \setminus \hyperfinite{X}{E}$ is non-empty.
  \end{enumerate}
\end{theorem}

\begin{theoremproof}
  Suppose that $E$ is not measure hyperfinite. Proposition \ref
  {preliminaries:measuredequivalencerelationsproposition:quasiinvariant} then
  yields an $E$-quasi-invariant \Borel probability measure $\mu$ on $X$ with
  respect to which $E$ is not $\mu$-hyperfinite, and Propositions \ref
  {preliminaries:measuredequivalencerelationsproposition:cocycles} and \ref
  {preliminaries:measurehyperfiniteness:proposition:characterization} give
  rise to an $E$-ergodic such measure.
\end{theoremproof}

We close this section by considering preservation of $\mu$-hyperfiniteness under \Borel 
homomorphisms.

\begin{proposition} \label{preliminaries:measurehyperfiniteness:proposition:smooth}
  Suppose that $X$ and $Y$ are standard \Borel spaces, $E$ is a countable \Borel equivalence
  relation on $X$, $\phi \from X \to Y$ is a \Borel homomorphism from $E$ to equality,
  $\mu$ is a \Borel probability measure on $X$, $\sequence{\mu_y}[y \in Y]$ is a \Borel
  disintegration of $\mu$ through $\phi$, and $\restriction{E}{\preimage{\phi}{y}}$ is
  $\mu_y$-hyperfinite for $(\pushforward{\phi}{\mu})$-almost every $y \in Y$. Then $E$ is
  $\mu$-hyperfinite.
\end{proposition}

\begin{propositionproof}
  By Theorem \ref{preliminaries:measurehyperfiniteness:theorem:uniform}, the set $D = \set{y \in Y}[E
  \text{ is $\mu_y$-hyperfinite}]$ is \Borel, and there is a hyperfinite \Borel equivalence
  relation $F$ on $X$ for which there is a \Borel set $C \subseteq X$ such that $\mu_y(C) = 1$
  and $\restriction{E}{C} = \restriction{F}{C}$ for all $y \in D$. Then $\mu(C) = 1$, so $E$ is
  $\mu$-hyperfinite.
\end{propositionproof}

\begin{proposition} \label{preliminaries:measurehyperfiniteness:proposition:hyperfinite}
  Suppose that $X$ and $Y$ are standard \Borel spaces, $E$ is a countable \Borel
  equivalence relation on $X$, $F$ is a hyperfinite \Borel equivalence relation on $Y$, $\phi \from
  X \to Y$ is a \Borel homomorphism from $E$ to $F$, $\mu$ is a \Borel probability measure on
  $X$, $\sequence{\mu_y}[y \in Y]$ is a \Borel disintegration of $\mu$ through $\phi$, and
  $\restriction{E}{\preimage{\phi}{y}}$ is $\mu_y$-hyperfinite for $(\pushforward{\phi}{\mu})$-almost
  every $y \in Y$. Then $E$ is $\mu$-hyperfinite.
\end{proposition}

\begin{propositionproof}
  Fix an increasing sequence $\sequence{F_n}[n \in \N]$ of finite \Borel equivalence relations on
  $Y$ whose union is $F$. Proposition \ref{preliminaries:measurehyperfiniteness:proposition:smooth}
  then ensures that each of the equivalence relations $E_n = E \intersection \preimage{(\phi \times
  \phi)}{F_n}$ is $\mu$-hyperfinite. As $E = \union[n \in \N][E_n]$, Theorem \ref
  {preliminaries:measurehyperfiniteness:theorem:increasingunion} implies that $E$ is
  $\mu$-hyperfinite.
\end{propositionproof}

\section{Actions of $\SL{2}{\Z}$} \label{preliminaries:SL2}

Let $\sim$ denote the equivalence relation on $\R[2] \setminus \set{\pair{0}{0}}$ given by
\begin{equation*}
  v \sim w \iff \exists r \in \R \ (r > 0 \mathand rv = w),
\end{equation*}
and let $\T$ denote the quotient. Define $\projection[\T] \from \R[2] \setminus \set
{\pair{0}{0}} \to \T$ by $\projection[\T](v) = \equivalenceclass{v}{\sim}$, and let
\definedsymbol{\SL{2}{\Z} \action \T} denote the action induced by $\SL{2}{\Z}
\action \R[2]$.

\begin{proposition}[\Jackson-\Kechris-\Louveau] \label{preliminaries:SL2:proposition:hyperfinite}
  The action $\SL{2}{\Z} \action \T$ is hyperfinite.
\end{proposition}

\begin{propositionproof}
  See the remark following the proof of \cite[Lemma 3.6]{JacksonKechrisLouveau}.
\end{propositionproof}

Let \definedsymbol{\ASL{2}{\Z}} denote the group of all functions $T \from
\R[2] \to \R[2]$ of the form $T(x) = Ax + b$, where $A \in \SL{2}{\Z}$ and $b \in \Z[2]$, under
composition. Define $\projection[\SL{2}{\Z}] \from \ASL{2}{\Z} \to \SL{2}
{\Z}$ by $\projection[\SL{2}{\Z}](Ax + b) = A$.

\begin{proposition} \label{preliminaries:SL2:proposition:nonhyperfinite}
  Suppose that $\mu$ is the \Borel probability measure on $\T[2]$ induced by \Lebesgue measure
  on $\R[2]$. Then the orbit equivalence relation $\orbitequivalencerelation{\SL{2}{\Z}}{\T[2]}$ is
  not $\mu$-hyperfinite.
\end{proposition}

\begin{propositionproof}
  As $\SL{2}{\Z}$ is not amenable and \cite[Lemma 3.6]{JacksonKechrisLouveau} ensures that
  $\SL{2}{\Z} \action \R[2]$ is free off of a $\mu$-null set, this is a consequence of Proposition \ref
  {preliminaries:measurehyperfiniteness:proposition:nonamenable}.
\end{propositionproof}

\begin{proposition}[\Jackson-\Kechris-\Louveau] \label{preliminaries:SL2:proposition:treeable}
  The orbit equivalence relation $\orbitequivalencerelation{\SL{2}{\Z}}{\T[2]}$ is treeable.
\end{proposition}

\begin{propositionproof}
  See \cite[Proposition 3.5]{JacksonKechrisLouveau}.
\end{propositionproof}

\section{Complexity} \label{preliminaries:complexity}

The conclusion of the following summarizes the main results of \cite{AdamsKechris}.

\begin{theorem}[\Adams-\Kechris] \label{preliminaries:complexity:theorem:antichain}
  Suppose that $X$ is a standard \Borel space, $E$ is a countable \Borel equivalence relation
  on $X$, $\sequence{E_r}[r \in \R]$ is a \Borel sequence of subequivalence relations of $E$,
  and $\sequence{\mu_r}[r \in \R]$ is a \Borel sequence of \Borel probability measures on $X$
  such that:
  \begin{enumerate}
    \item Each $\mu_r$ is $E_r$-ergodic and $E_r$-quasi-invariant.
    \item The relation $E_r$ is $\mu_r$-nowhere reducible to the relation $E_s$, for all distinct
      $r, s \in \R$.
  \end{enumerate}
  Then the following hold:
  \begin{enumerate}
    \renewcommand{\theenumi}{\alph{enumi}}
    \item There is an embedding of containment on \Borel subsets of $\R$ into \Borel
      reducibility of countable \Borel equivalence relations with smooth-to-one \Borel
      homomorphisms to $E$ (in the codes).
    \item \Borel bi-reducibility and reducibility of countable \Borel equival\-ence relations with
      smooth-to-one \Borel homomorphisms to $E$ are both $\Sigmaclass[1][2]$-complete
      (in the codes).
    \item Every \Borel quasi-order is \Borel reducible to \Borel reducib\-ility of countable \Borel
      equivalence relations with smooth-to-one \Borel homomorphisms to $E$.
    \item \Borel and $\sigmaclass{\Sigmaclass[1][1]}$-measurable reducibility do not agree on the
      countable \Borel equivalence relations with smooth-to-one \Borel homomorphisms to $E$.
  \end{enumerate}
\end{theorem}

\begin{theoremproof}
  The proof of \cite[Theorem 4.1]{AdamsKechris} yields (a), the proof of \cite[Theorem 5.1]
  {AdamsKechris} yields (b), the final paragraph of \cite[\S7]{AdamsKechris} yields (c), and
  the proof of \cite[Theorem 5.5]{AdamsKechris} yields (d).
\end{theoremproof}

\part{Tools}

Here we introduce the new ideas underlying our arguments. In \S\ref
{productivehyperfiniteness}, we show that $\SL{2}{\Z} \action \T$ satisfies a measureless
strengthening of amenability. In \S\ref{projectiverigidity}, we use this to prove that
$\ASL{2}{\Z} \action \R[2]$ satisfies a measureless local rigidity property. In \S\ref
{projectiveseparability}, we establish a strong separability property for orbit equivalence
relations induced by such actions. In \S\ref{measures}, we show that the latter yields
countability of an appropriate auxiliary equivalence relation on the underlying space of ergodic
quasi-invariant \Borel probability measures witnessing the failure of hyperfiniteness, and we
derive several consequences of this countability. In \S\ref{stratification}, we provide a general
stratification theorem for treeable countable \Borel equivalence relations.

\section{Productive hyperfiniteness} \label{productivehyperfiniteness}

Suppose that $\Gamma$ is a countable discrete group. The \emph{diagonal product} of actions
$\Gamma \action X$ and $\Gamma \action Y$ is the action $\Gamma \action X \times Y$ given by
$\gamma \cdot (x,y) = (\gamma \cdot x, \gamma \cdot y)$. We say that a \Borel action
$\Gamma \action X$ on a standard \Borel space is \definedterm{productively hyperfinite} if
whenever $\Gamma \action Y$ is a \Borel action on a standard \Borel space, the orbit equivalence
relation induced by the diagonal product action $\Gamma \action X \times Y$ is hyperfinite.

\begin{proposition} \label{productivehyperfiniteness:proposition:sufficientcondition}
  Suppose that $\Gamma$ is a countable discrete group, $X$ is a standard \Borel space, and
  $\Gamma \action X$ is a \Borel action such that:
  \begin{enumerate}
    \item The induced orbit equivalence relation is hyperfinite.
    \item The stabilizer of every point is hyperfinite.
    \item Only countably-many points have infinite stabilizers.
  \end{enumerate}
  Then $\Gamma \action X$ is productively hyperfinite.
\end{proposition}

\begin{propositionproof}
  Let $C$ denote the set of points whose stabilizers are infinite, and fix an increasing
  sequence $\sequence{E_n}[n \in \N]$ of finite \Borel equivalence relations whose union is
  $\orbitequivalencerelation{\Gamma}{X}$.

  Suppose now that $Y$ is a standard \Borel space and $\Gamma \action Y$ is a \Borel action. For
  each $n \in \N$, let $F_n$ denote the equivalence relation on $(X \setminus C) \times Y$ for
  which two \heightcorrection{$\orbitequivalencerelation{\Gamma}{(X \setminus C) \times
  Y}$}-equivalent pairs $\pair{x}{y}$ and $\pair{x'}{y'}$ are related exactly when $x \mathrel{E_n}
  x'$. As each $F_n$ is finite and their union is \heightcorrection{$\orbitequivalencerelation
  {\Gamma}{(X \setminus C) \times Y}$}, the latter equivalence relation is hyperfinite.

  It only remains to show that $\orbitequivalencerelation{\Gamma}{C \times Y}$ is hyperfinite. As
  $C$ is countable and Proposition \ref{preliminaries:hyperfiniteness:proposition:closure:unions}
  ensures that the family of \Borel sets on which a \Borel equivalence relation is
  hyperfinite forms a $\sigma$-ideal, we need only show that \heightcorrection
  {$\orbitequivalencerelation{\Gamma}{X \times Y}$} is hyperfinite on $\singleton{x} \times Y$, for
  all $x \in C$. But this follows from the fact that its restriction to such a set is the orbit equivalence 
  relation induced by a \Borel action of the stabilizer of $x$.
\end{propositionproof}

To apply this to $\SL{2}{\Z} \action \T$, we must first consider its stabilizers.

\begin{proposition} \label{productivehyperfiniteness:proposition:stabilizer}
  Suppose that $\theta \in \T$. Then the stabilizer of $\theta$ under $\SL{2}{\Z} \action \T$ is either
  trivial or infinite cyclic.
\end{proposition}

\begin{propositionproof}
  We consider first the case that $\theta \intersection \Z[2] \neq \emptyset$. Let $v$ denote the
  unique element of $\theta \intersection \Z[2]$ of minimal length. Note that the stabilizers of
  $\theta$ and $v$ are one and the same, for if $A$ is in the stabilizer of $\theta$, then $v$ is an
  eigenvector of $A$, so minimality ensures that $Av = v$. Minimality also ensures that the
  coordinates of $v$ are relatively prime, so there exists $a \in \Z[2]$ such that $a \cdot v = 1$, in
  which case $B = \matrix{a_1}{a_2}{-v_2}{v_1}$ is a matrix in $\SL{2}{\Z}$ for which $Bv = \vector
  {1}{0}$, thus conjugation by $B$ yields an isomorphism of the stabilizer of $v$ with that
  of $\vector{1}{0}$, and the latter is the infinite cyclic group $\set{\matrix{1}{n}{0}{1}}[n \in \Z]$.

  It remains to consider the case that $\theta \intersection \Z[2] = \emptyset$. Fix $v \in \theta$. An
  elementary calculation reveals that the stabilizer of $v$ is trivial. Let $\Lambda$ denote the set of
  eigenvalues of matrices in the stabilizer of $\theta$, noting that $\Lambda$ forms a group under
  multiplication.

  \begin{lemma} \label{productivehyperfiniteness:lemma:generator}
    The group $\Lambda$ is cyclic.
  \end{lemma}

  \begin{lemmaproof}
    It is sufficient to show that $1$ is isolated in $\Lambda \intersection \closedopeninterval{1}
    {\infty}$. Towards this end, suppose that $A$ is in the stabilizer of $\theta$ and $v$ is an
    eigenvector of $A$ with eigenvalue $\lambda > 1$. If $\mu$ is the other eigenvalue of $A$,
    then $\lambda \mu = \determinant{A} = 1$, so $\trace{A} = \lambda + \mu = \lambda + 1 /
    \lambda$. As $\trace{A} \in \Z$, another elementary calculation reveals that $\lambda \ge (3 +
    \sqrt{5}) / 2$.
  \end{lemmaproof}

  By Lemma \ref{productivehyperfiniteness:lemma:generator}, there is a matrix $A$ in the
  stabilizer of $\theta$ which has an eigenvalue $\lambda$ generating $\Lambda$. Note that if
  $B$ is any matrix in the stabilizer of $\theta$, then there exists $n \in \Z$ for which $v$ is an
  eigenvector of $B$ with eigenvalue $\lambda^n$, in which case $A^nB^{-1}$ is in the stabilizer
  of $v$, so $B = A^n$, thus $A$ generates the stabilizer of $\theta$, hence the latter is cyclic.

  Observe finally that if $A$ is a non-identity matrix fixing $\theta$, then any two distinct powers of
  $A$ are themselves distinct, since the eigenvalues corresponding to $v$ are distinct. In
  particular, it follows that if the stabilizer of $\theta$ is non-trivial, then it is infinite.
\end{propositionproof}

As a consequence, we can now obtain the main result of this section.

\begin{proposition} \label{productivehyperfiniteness:proposition:SL2Z}
  The action $\SL{2}{\Z} \action \T$ is productively hyperfinite.
\end{proposition}

\begin{propositionproof}
  As Proposition \ref{preliminaries:SL2:proposition:hyperfinite} ensures that the orbit equivalence
  relation induced by $\SL{2}{\Z} \action \T$ is hyperfinite, Proposition \ref
  {productivehyperfiniteness:proposition:stabilizer} ensures that the non-trivial stabilizers of $\SL
  {2}{\Z} \action \T$ are infinite cyclic, and Proposition \ref{preliminaries:hyperfiniteness:proposition:Z}
  ensures that infinite cyclic groups are hyperfinite, it is sufficient to show that only countably
  many $\theta \in \T$ have non-trivial stabilizers, by Proposition \ref
  {productivehyperfiniteness:proposition:sufficientcondition}. As every such $\theta$ is the
  equivalence class of an eigenvector of some non-trivial matrix in $\SL{2}{\Z}$, and every such
  matrix admits at most two such classes of eigenvectors, this follows from the countability of $\SL
  {2}{\Z}$.
\end{propositionproof}

\section{Projective rigidity} \label{projectiverigidity}

Given $R \subseteq X \times X$, $\Delta \action Y$, and $\rho \from R \to \Delta$, we say that
a function $\phi \from X \to Y$ is \definedterm{$\rho$-invariant} if $x_1 \mathrel{R} x_2 \implies \phi
(x_1) = \rho(x_1, x_2) \cdot \phi(x_2)$ for all $x_1, x_2 \in X$. The \definedterm{difference set}
associated with two functions $\phi \from A \subseteq X \to Y$ and $\psi \from B \subseteq X \to Y$
is given by
\begin{equation*}
  \differenceset{\phi}{\psi} = \set{x \in A \intersection B}[\phi(x) \neq \psi(x)] \union (A
    \symmetricdifference B).
\end{equation*}
We say that $\Delta \action Y$ is \definedterm{projectively rigid} if
whenever $X$ is a standard \Borel space, $E$ is a countable \Borel equivalence relation on $X$,
and $\rho \from E \to \Delta$ is a \Borel function, there is essentially at most one countable-to-one
$\rho$-invariant \Borel function, in the sense that for any two such functions $\phi$ and $\psi$, the
relation $\restriction{E}{\differenceset{\phi}{\psi}}$ is hyperfinite.

\begin{theorem} \label{projectiverigidity:theorem:ASL2Z}
  The action $\ASL{2}{\Z} \action \R[2]$ is projectively rigid.
\end{theorem}

\begin{theoremproof}
  Suppose that $X$ is a standard \Borel space, $E$ is a countable \Borel equivalence relation on
  $X$, $\rho \from E \to \ASL{2}{\Z}$ is a \Borel function, $\phi \from X \to \R[2]$ is a
  countable-to-one $\rho$-invariant \Borel function, and $\psi \from X \to \R[2]$ is a
  $\rho$-invariant \Borel function.

  Define $\pi \from \differenceset{\phi}{\psi} \to \T$ by $\pi(x) = \projection[\T](\phi(x) - \psi(x))$,
  and define $\sigma \from \restriction{E}{\differenceset{\phi}{\psi}} \to \SL{2}{\Z}$ by $\sigma(x_1,
  x_2) = \projection[\SL{2}{\Z}](\rho(x_1, x_2))$.

  \begin{lemma}
    The function $\pi$ is $\sigma$-invariant.
  \end{lemma}

  \begin{lemmaproof}
    Simply observe that if $x_1 \mathrel{(\restriction{E}{\differenceset{\phi}{\psi}})} x_2$, then
    \begin{align*}
      \pi(x_1)
        & = \projection[\T](\phi(x_1) - \psi(x_1)) \\
        & = \projection[\T](\rho(x_1, x_2) \cdot \phi(x_2) - \rho(x_1, x_2) \cdot \psi(x_2)) \\
        & = \projection[\T](\sigma(x_1, x_2) \cdot \phi(x_2) - \sigma(x_1, x_2) \cdot \psi(x_2)) \\
        & = \projection[\T](\sigma(x_1, x_2) \cdot (\phi(x_2) - \psi(x_2))) \\
        & = \sigma(x_1, x_2) \cdot \projection[\T](\phi(x_2) - \psi(x_2)) \\
        & = \sigma(x_1, x_2) \cdot \pi(x_2),
    \end{align*}
    thus $\pi$ is $\sigma$-invariant.
  \end{lemmaproof}

  As $\restriction{(\projection[{\T[2]}] \composition \phi)}{\differenceset{\phi}{\psi}}$ is also
  $\sigma$-invariant, it follows that the product $\pi \times (\projection[{\T[2]}] \composition \phi)$ is
  a countable-to-one homomorphism from $\restriction{E}{\differenceset{\phi}{\psi}}$ to the
  orbit equivalence relation induced by the diagonal product action $\SL{2}{\Z} \action \T \times
  \T[2]$. As Proposition \ref{productivehyperfiniteness:proposition:SL2Z} ensures that $\SL{2}{\Z}
  \action \T$ is productively hyperfinite, it follows that the latter relation is hyperfinite. As
  Proposition \ref{preliminaries:hyperfiniteness:proposition:closure:homomorphisms} ensures that the
  family of hyperfinite \Borel equivalence relations is closed downward under countable-to-one
  \Borel homomorphism, it follows that the former relation is also hyperfinite.
\end{theoremproof}

\begin{remark}
  As noted by both Manuel Inselmann and one of the anonymous referees, the productive
  hyperfiniteness of $\SL{2}{\Z} \action \T$ can also be used to show that the orbit equivalence
  relation induced by $\SL{2}{\Z} \action \R[2]$ is hyperfinite. To see this, observe that the
  function $\pi \from \R[2] \setminus \set{0} \to \T \times \R[2]$ given by $\pi(x) = \pair{\projection[\T](x)}{x}$ is
  a reduction of the orbit equivalence relation induced by $\SL{2}{\Z} \action (\R[2] \setminus \set{0})$ to the
  orbit equivalence relation induced by $\SL{2}{\Z} \action \T \times \R[2]$.
\end{remark}

\section{Projective separability} \label{projectiveseparability}

Let $\Borelfunctions{X}{Y}$ denote the set of \Borel functions $\phi \from B \to Y$, where $B$
varies over \Borel subsets of $X$. Let $\functions{X}[\mu]{Y}$ denote $\Borelfunctions{X}{Y}$
equipped with the pseudo-metric $\uniformmetric{\mu}(\phi, \psi) = \mu
(\differenceset{\phi}{\psi})$.

\begin{proposition} \label{projectiveseparability:proposition:characterization}
  Suppose that $X$ and $Y$ are standard \Borel spaces, $\mu$ is a finite \Borel measure on $X$,
  and $\calL \subseteq \functions{X}[\mu]{Y}$. Then the following are equivalent:
  \begin{enumerate}
    \item The space $\calL$ is separable.
    \item There is a \Borel set $R \subseteq X \times Y$, whose vertical sections are countable,
      with the property that
      \begin{equation*}
        \forall \phi \in \calL \ \mu(\set{x \in \domain{\phi}}[\neg x \mathrel{R} \phi(x)]) = 0.
      \end{equation*}
  \end{enumerate}
\end{proposition}

\begin{propositionproof}
  To see $(1) \implies (2)$, note that if $\calD$ is a countable dense subset of $\calL$, then
  the set $R = \bigcup_{\phi \in \calD} \graph{\phi}$ is as desired, since graphs of \Borel functions are
  \Borel. To see $(2) \implies (1)$, it is sufficient to show that if condition (2) holds, then there is a
  countable subset of $\functions{X}[\mu]{Y}$ whose closure contains $\calL$. As the vertical
  sections of $R$ are countable, the \Lusin-\Novikov uniformization theorem yields
  a countable family $\calF$ of \Borel partial functions, the union of whose graphs is $R$. Fix
  a countable algebra $\calB$ of \Borel subsets of $X$, containing the domain of every $\phi
  \in \calF$, such that for all \Borel sets $A \subseteq X$ and all $\epsilon > 0$, there exists
  $B \in \calB$ with $\mu(A \symmetricdifference B) \le \epsilon$.  We then obtain the desired
  countable dense family by considering those $\psi \from B \to Y$, where $B$ ranges over
  $\calB$, for which there is a finite partition $\calA \subseteq \calB$ of $B$ such
  that $\forall A \in \calA \exists \phi \in \calF \ \restriction{\phi}{A} = \restriction{\psi}{A}$.
\end{propositionproof}

We say that a function $\phi \from Y \to Y'$ is a \definedterm{homomorphism} from a set $\calL
\subseteq \functions{X}[\mu]{Y}$ to a set $\calL'  \subseteq \functions{X}[\mu]{Y'}$ if $\forall \psi \in
\calL \ \phi \composition \psi \in \calL'$.

\begin{proposition} \label{projectiveseparability:proposition:closure:general}
  Suppose that $X$, $Y$, and $Y'$ are standard \Borel spaces, $\mu$ is a \Borel probability
  measure on $X$, $\calL \subseteq \functions{X}[\mu]{Y}$ and $\calL' \subseteq \functions{X}
  [\mu]{Y'}$, there is a countable-to-one \Borel homomorphism $\phi \from Y \to Y'$ from $\calL$ to
  $\calL'$, and $\calL'$ is separable. Then $\calL$ is separable.
\end{proposition}

\begin{propositionproof}
  Fix a \Borel set $R' \subseteq X \times Y'$ satisfying the analog of condition (2) of Proposition \ref
  {projectiveseparability:proposition:characterization} for $\calL'$, and observe that the set $R =
  \preimage{(\id \times \phi)}{R'}$ satisfies condition (2) of Proposition \ref
  {projectiveseparability:proposition:characterization} for $\calL$.
\end{propositionproof}

Let \definedsymbol{\homomorphisms{E}{\mu}{F}} denote the subspace of $\functions{X}
[\mu]{Y}$ consisting of all countable-to-one partial homomorphisms $\phi \in \functions{X}
[\mu]{Y}$ from $E$ to $F$.

\begin{proposition} \label{projectiveseparability:proposition:closure:special}
  Suppose that $X$, $Y$, and $Y'$ are standard \Borel spaces, $E$, $F$, and $F'$ are countable
  \Borel equivalence relations on $X$, $Y$, and $Y'$, $\mu$ is a \Borel probability measure on
  $X$, there is a countable-to-one \Borel homomorphism $\phi \from Y \to Y'$ from $F$ to $F'$,
  and $\homomorphisms{E}{\mu}{F'}$ is separable. Then $\homomorphisms{E}{\mu}{F}$ is
  separable.
\end{proposition}

\begin{propositionproof}
  As the function $\phi$ is also a homomorphism from $\homomorphisms{E}{\mu}{F}$ to
  $\homomorphisms{E}{\mu}{F'}$, the desired result follows from Proposition \ref
  {projectiveseparability:proposition:closure:general}.
\end{propositionproof}

We say that $F$ is \definedterm{projectively separable} if whenever $X$ is a standard \Borel
space, $E$ is a countable \Borel equivalence relation on $X$, and $\mu$ is a \Borel probability
measure on $X$ with respect to which $E$ is $\mu$-nowhere hyperfinite, the space
$\homomorphisms{E}{\mu}{F}$ is separable.

\begin{proposition} \label{projectiveseparability:proposition:closure}
  Suppose that $X$ and $Y$ are standard \Borel spaces, $E$ and $F$ are countable \Borel
  equivalence relations on $X$ and $Y$, there is a countable-to-one \Borel homomorphism
  from $E$ to $F$, and $F$ is projectively separable. Then $E$ is projectively separable.
\end{proposition}

\begin{propositionproof}
  This is a direct consequence of Proposition \ref
  {projectiveseparability:proposition:closure:special}.
\end{propositionproof}

We next establish the connection between projective rigidity and projective separability.

\begin{theorem} \label{projectiveseparability:theorem:sufficientcondition}
  Suppose that $\Delta$ is a countable discrete group, $Y$ is a standard \Borel space, and
  $\Delta \action Y$ is a projectively rigid \Borel action. Then the orbit equivalence
  relation $F = \orbitequivalencerelation{\Delta}{Y}$ is projectively separable.
\end{theorem}

\begin{theoremproof}
  Suppose that $X$ is a standard \Borel space, $E$ is a countable \Borel equivalence relation
  on $X$, and $\mu$ is a \Borel probability measure on $X$ with respect to which $E$ is
  $\mu$-nowhere hyperfinite. Let \definedsymbol{\countingmeasure} denote the counting
  measure on $X$. The \Lusin-\Novikov uniformization theorem yields an increasing sequence
  $\sequence{R_n}[n \in \N]$ of \Borel subsets of $X \times X$ such that $E = \union[n \in \N]
  [R_n]$ and every vertical section of every $R_n$ has cardinality at most $n$. Set $\nu_n =
  \restriction{(\mu \times \countingmeasure)}{R_n}$ for all $n \in \N$.
  
  \begin{lemma} \label{projectiveseparability:lemma:limits}
    Suppose that $\phi \in \homomorphisms{E}{\mu}{F}$, $\rho \from \restriction{E}{\domain{\phi}}
    \to \Delta$ is a \Borel function with respect to which $\phi$ is invariant, $\sequence{D_n}[n \in
    \N]$ is a sequence of \Borel subsets of $X$ with $\sum_{n \in \N}
    \mu(\domain{\phi} \symmetricdifference D_n) < \infty$, $\sequence{\rho_n \from \restriction{R_n}{D_n}
    \to \Delta}[n \in \N]$ is a sequence of \Borel functions such that $\sum_{n \in \N} \uniformmetric
    {\nu_n}(\restriction{\rho}{(\restriction{R_n}{\domain{\phi}})}, \rho_n) < \infty$, and $\phi_n \from D_n \to
    Y$ is a $\rho_n$-invariant \Borel function for all $n \in \N$. Then $\uniformmetric{\mu}(\phi,
    \phi_n) \goesto 0$.
  \end{lemma}
  
  \begin{lemmaproof}
    For all $n \in \N$, let $E_n$ be the equivalence relation on
    $\domain{\phi} \intersection D_n$ generated by the relation $S_n =
    (\restriction{R_n}{(\domain{\phi} \intersection D_n)}) \setminus \differenceset{\rho}{\rho_n}$.

    \begin{sublemma}
      For all $n \in \N$, there is a \Borel function $\sigma_n \from E_n \to
      \Delta$ for which every $(\restriction{\rho}{S_n})$-invariant function is
      $\sigma_n$-invariant.
    \end{sublemma}

    \begin{sublemmaproof}
      Note that if $x \mathrel{E_n} y$, then there are only countably many $\ell
      \in \N$ and $\sequence{z_i}[i \le \ell] \in \functions{\ell+1}{X}$ such
      that $x = z_0$, $\forall i < \ell \ z_i \mathrel{S_n} z_{i+1}$, and $y = z_\ell$,
      so the \Lusin-\Novikov uniformization theorem yields \Borel functions $\ell
      \from E_n \to \N$ and $f \from E_n \to \functions{<\N}{X}$ with the property 
      that
      \begin{equation*}
        \forall \pair{x}{y} \in E_n \ x = f_0(x, y) \mathrel{S_n} f_1(x, y) \mathrel{S_n}
          \cdots \mathrel{S_n} f_{\ell(x, y)}(x, y) = y.
      \end{equation*}
      Define $\sigma_n(x, y) = \prod_{i < \ell(x, y)} \rho(f_i(x, y), f_{i+1}(x, y))$.
    \end{sublemmaproof}

    As the restrictions of $\phi$ and $\phi_n$ to $\domain{\phi} \intersection D_n$
    are $(\restriction{\rho}{S_n})$-invariant, they are $\sigma_n$-invariant. Note
    that the set $D = \domain{\phi} \intersection \union[n \in \N][{\intersection[m
    \ge n][D_m]}]$ is $(\restriction{\mu}{\domain{\phi}})$-conull.

    \begin{sublemma} \label{projectiveseparability:sublemma:conullunion}
      There is a $(\restriction{\mu}{D})$-conull \Borel set $C \subseteq D$
      such that $E \intersection (C \times D) \subseteq \union[n \in \N]
      [{\intersection[m \ge n][S_m]}]$.
    \end{sublemma}

    \begin{sublemmaproof}
      The \Lusin-\Novikov uniformization theorem ensures that the sets $C_n
      = \set{x \in \domain{\phi} \intersection D_n}[\exists y \in \domain{\phi}
      \intersection D_n \ x \mathrel{(R_n \setminus S_n)} y]$ are \Borel, and
      \Fubini's theorem (see, for example, \cite[\S17.A]{Kechris}) ensures that
      $\mu(C_n) \le \uniformmetric{\nu_n}(\restriction{\rho}{(\restriction{R_n}
      {\domain{\phi}})}, \rho_n)$ for all $n \in \N$. In particular, it follows that the set
      $C = D \setminus \intersection[n \in \N][{\union[m \ge n][C_m]}]$ is $(\restriction
      {\mu}{D})$-conull. And if $\pair{x}{y} \in E \intersection (C \times D)$,
      then there exists $n \in \N$ for which $x \mathrel{R_n} y$, so the fact that $x \in
      C$ ensures that $x \mathrel{S_m} y$ for sufficiently large $m \ge n$.
    \end{sublemmaproof}
    
    Suppose now that $\epsilon > 0$. Set $F_n = \intersection[m \ge n][E_m]$ for
    all $n \in \N$, and observe that $\restriction{E}{C} = \union[n \in \N][\restriction{F_n}
    {C}]$. As Theorem \ref{preliminaries:measurehyperfiniteness:theorem:increasingunion}
    ensures that the $\mu$-hyperfinite \Borel equivalence relations are closed under
    increasing unions, there are \Borel sets $B_n \subseteq C \intersection \intersection
    [m \ge n][D_m]$ with the property that $\mu(C \setminus B_n) < \epsilon$ and 
    $\restriction{F_n}{B_n}$ is $(\restriction{\mu}{B_n})$-nowhere hyperfinite, thus
    $\restriction{\phi}{B_n} = \restriction{\phi_n}{B_n}$, for sufficiently large $n \in \N$.
  \end{lemmaproof}
  
  Fix a countable family $\calB$ of \Borel subsets of $X$ such that for all \Borel
  sets $A \subseteq X$ and all real numbers $\epsilon > 0$, there exists $B \in
  \calB$ with $\mu(A \symmetricdifference B) \le \epsilon$. Proposition \ref
  {projectiveseparability:proposition:characterization} yields countable dense sets
  $\calD_n \subseteq \functions{R_n}[\nu_n]{\Delta}$. For each $n \in \N$, $B \in
  \calB$, $\epsilon \in \openinterval{0}{\infty} \intersection \Q$, and $\sigma \in
  \calD_n$ for which it is possible, fix a \Borel set $D_{n, B, \epsilon, \sigma}
  \subseteq X$ with $\mu(B \symmetricdifference D_{n, B, \epsilon, \sigma}) \le
  \epsilon$, a \Borel function $\rho_{n, B, \epsilon, \sigma} \from \restriction{R_n}
  {D_{n, B, \epsilon, \sigma}} \to \Delta$ such that $\uniformmetric{\nu_n}(\sigma,
  \rho_{n, B, \epsilon, \sigma}) \le \epsilon$, and a $\rho_{n, B, \epsilon,
  \sigma}$-invariant \Borel function $\phi_{n, B, \epsilon, \sigma} \from
  D_{n, B, \epsilon, \sigma} \to Y$. It only remains to check that the functions of
  the form $\phi_{n, B, \epsilon, \sigma}$ are dense in $\homomorphisms{E}
  {\mu}{F}$.

  Towards this end, suppose that $\phi \in \homomorphisms{E}{\mu}{F}$, and fix
  a \Borel function $\rho \from \restriction{E}{\domain{\phi}} \to \Delta$ for which
  $\phi$ is $\rho$-invariant. Fix a sequence $\sequence{\epsilon_n}[n \in \N]$ of
  positive rational numbers for which $\sum_{n \in \N} \epsilon_n < \infty$, and for
  each $n \in \N$, fix $B_n \in \calB$ with $\mu(B_n \symmetricdifference \domain
  {\phi}) \le \epsilon_n$ and $\sigma_n \in \calD_n$ such that $\uniformmetric
  {\nu_n}(\sigma_n, \restriction{\rho}{(\restriction{R_n}{\domain{\phi}})}) \le
  \epsilon_n$. Then the sets $D_n = D_{n, B_n, \epsilon_n, \sigma_n}$ and the
  functions $\rho_n = \sigma_{n, B_n, \epsilon_n, \sigma_n}$ and $\phi_n =
  \phi_{n, B_n, \epsilon_n, \sigma_n}$ are well-defined. As $\mu(\domain{\phi}
  \symmetricdifference D_n) \le 2\epsilon_n$ and $\uniformmetric{\nu_n}
  (\restriction{\rho}{(\restriction{R_n}{\domain{\phi}})}, \rho_n) \le 2\epsilon_n$
  for all $n \in \N$, Lemma \ref{projectiveseparability:lemma:limits} ensures that
  $\uniformmetric{\mu}(\phi, \phi_n) \goesto 0$.
\end{theoremproof}

In particular, we can now establish the existence of non-trivial projec\-tively-separable countable
\Borel equivalence relations.

\begin{theorem} \label{projectiveseparability:theorem:sl2action}
  The orbit equivalence relation induced by $\SL{2}{\Z} \action \T[2]$ is
  projectively separable.
\end{theorem}

\begin{theoremproof}
  Note that the orbit equivalence relation in question is \Borel reducible to that induced by $\ASL{2}
  {\Z} \action \R[2]$. As Theorem \ref{projectiverigidity:theorem:ASL2Z} ensures that the latter
  action is projectively rigid, its induced orbit equivalence relation is projectively separable by
  Theorem \ref{projectiveseparability:theorem:sufficientcondition}. But Proposition \ref
  {projectiveseparability:proposition:closure} ensures that the projectively-separable
  countable \Borel equivalence relations are closed under \Borel reducibility.
\end{theoremproof}

\section{The space of measures} \label{measures}

Here we consider connections between $E$ and $\ergodicquasiinvariant{X}{E}
\setminus \hyperfinite{X}{E}$. Theorems \ref
{preliminaries:measuredequivalencerelations:theorem:Borel} and \ref
{preliminaries:measurehyperfiniteness:theorem:uniform} ensure that the latter is
a \Borel subset of $\probabilitymeasures{X}$.

\begin{proposition} \label{measures:proposition:successor:sufficientcondition}
  Suppose that $X$ is a standard \Borel space, $E$ is a countable \Borel equivalence
  relation on $X$, and the set $\ergodicquasiinvariant{X}{E} \setminus \hyperfinite{X}{E}$ is a
  single measure-equivalence class. Then $E$ is a successor of $\Ezero$ under measure
  reducibility.
\end{proposition}

\begin{propositionproof}
  Suppose that $Y$ is a standard \Borel space and $F$ is a countable \Borel equivalence relation
  on $Y$ which is measure reducible to $E$, but not to $\Ezero$. We must show that $E$ is
  measure reducible to $F$.

  By Theorem \ref{preliminaries:measurehyperfiniteness:theorem:E0}, there exists $\nu \in 
  \ergodicquasiinvariant{Y}{F} \setminus \hyperfinite{Y}{F}$. By Proposition \ref
  {preliminaries:measuredequivalencerelationsproposition:null}, there is a $\nu$-null \Borel set $N
  \subseteq Y$ on which $F$ is non-smooth. As $F$ is countable, the \Lusin-\Novikov uniformization
  theorem ensures that $\saturation{N}{F}$ is \Borel, so by
  replacing $N$ with $\saturation{N}{F}$, we can assume that $N$ is $F$-invariant. Fix a
  $\nu$-conull \Borel set $C \subseteq \setcomplement{N}$ for which there is a \Borel reduction
  $\phi \from C \to X$ of $\restriction{F}{C}$ to $E$. As $E$ and $F$ are countable, the
  \Lusin-\Novikov uniformization theorem ensures that the set $B = \saturation{\image{\phi}
  {C}}{E}$ is \Borel, and that there is a \Borel function $\psi \from B \to C$ such that $\graph{\phi
  \composition \psi} \subseteq E$. In particular, it follows that $\psi$ is a \Borel reduction of
  $\restriction{E}{B}$ to $\restriction{F}{C}$.

  Suppose now that $\mu$ is a \Borel probability measure on $X$. As Propositon \ref
  {preliminaries:hyperfiniteness:proposition:closure:homomorphisms} ensures that the class of
  hyperfinite \Borel equivalence relations is closed downward under \Borel reducibility, it follows
  that the push-forward $\nu'$ of $\restriction{\nu}{C}$ through $\phi$ is not in $\hyperfinite{X}
  {E}$. By Proposition \ref{preliminaries:measuredequivalencerelationsproposition:quasiinvariant},
  there is an $E$-quasi-invariant \Borel probability measure $\nu''$ on $X$ such that $\nu'
  \absolutelycontinuous \nu''$ and the two measures have the same $E$-invariant null \Borel
  sets. Then $\nu'' \in \ergodicquasiinvariant{X}{E} \setminus \hyperfinite{X}{E}$, so $\restriction
  {E}{\setcomplement{B}}$ is measure hyperfinite, thus there is a \Borel set $A \subseteq
  \setcomplement{B}$ such that $\restriction{E}{A}$ is hyperfinite and $\mu(A \union B) = 1$.
  As Theorem \ref{preliminaries:hyperfiniteness:theorem:embedding} ensures that every hyperfinite
  \Borel equivalence relation is \Borel reducible to every non-smooth \Borel equivalence relation,
  it follows that there is a \Borel reduction $\psi' \from A \to N$ of $\restriction{E}{A}$ to
  $\restriction{F}{N}$. As $\psi \union \psi'$ is a reduction of $\restriction{E}{(A \union B)}$
  to $F$, it follows that $E$ is $\mu$-reducible to $F$, thus $E$ is measure reducible to $F$.
\end{propositionproof}

\begin{proposition} \label{measures:proposition:successor:sufficientcondition:countableunion}
  Suppose that $X$ is a standard \Borel space, $E$ is a countable \Borel equivalence relation on
  $X$, and $\ergodicquasiinvariant{X}{E} \setminus \hyperfinite{X}{E}$ is a non-empty countable
  union of measure-equivalence classes. Then $E$ is a countable disjoint union of successors
  of $\Ezero$ under measure reducibility.
\end{proposition}

\begin{propositionproof}
  Suppose that $N$ is a non-empty countable set and $\ergodicquasiinvariant{X}{E} \setminus
  \hyperfinite{X}{E}$ is the disjoint union of the measure-equivalence classes of \Borel
  probability measures $\mu_n$ on $X$, for $n \in N$. Fix a partition $\sequence{B_n}[n \in N]$ of
  $X$ into $E$-invariant \Borel sets with the property that $\mu_n(B_n) = 1$ for all $n \in N$, and
  observe that Proposition \ref{measures:proposition:successor:sufficientcondition} ensures that
  each $\restriction{E}{B_n}$ is a successor of $\Ezero$ under measure reducibility.
\end{propositionproof}

On the other hand, we have the following.

\begin{proposition} \label{measures:proposition:perfectsequence:sufficientcondition}
  Suppose that $X$ is a standard \Borel space, $E$ is a countable \Borel equivalence relation on
  $X$, and $\ergodicquasiinvariant{X}{E} \setminus \hyperfinite{X}{E}$ is not a countable
  union of measure-equivalence classes. Then there are \Borel sequences $\sequence{B_c}[c \in
  \Cantorspace]$ of pairwise disjoint $E$-invariant subsets of $X$ and $\sequence{\mu_c}[c \in
  \Cantorspace]$ of \Borel probability measures on $X$ in $\ergodicquasiinvariant{X}{E} \setminus
  \hyperfinite{X}{E}$ such that $\mu_c(B_c) = 1$ for all $c \in \Cantorspace$.
\end{proposition}

\begin{propositionproof}
  As measure equivalence is \Borel, Theorem \ref
  {preliminaries:Borelequivalencerelations:theorem:perfect} yields a \Borel sequence $\sequence
  {\mu_c}[c \in \Cantorspace]$ of pairwise orthogonal \Borel probability measures on $X$ in
  $\ergodicquasiinvariant{X}{E} \setminus \hyperfinite{X}{E}$. Theorem \ref
  {preliminaries:measures:theorem:perfect} then implies that by thinning down $\sequence{\mu_c}[c \in
  \Cantorspace]$, we can ensure the existence of a \Borel sequence $\sequence{A_c}[c \in
  \Cantorspace]$ of pairwise disjoint \Borel subsets of $X$ such that $\mu_c(A_c) = 1$ for all
  $c \in \Cantorspace$. Define $B_c = \set{x \in X}[\equivalenceclass{x}{E} \subseteq A_c]$.
\end{propositionproof}

Combining the previous two results yields the following.

\begin{proposition}
  Suppose that $X$ is a standard \Borel space and $E$ is a non-measure-hyperfinite
  countable \Borel equivalence relation on $X$. Then at least one of the following holds:
  \begin{enumerate}
    \item The relation $E$ is a countable disjoint union of successors of $\Ezero$ under
      measure reducibility.
    \item There are \Borel sequences $\sequence{B_c}[c \in \Cantorspace]$ of pairwise disjoint
      $E$-invar\-iant subsets of $X$ and $\sequence{\mu_c}[c \in \Cantorspace]$ of \Borel
      probability measures on $X$ in $\ergodicquasiinvariant{X}{E} \setminus \hyperfinite{X}{E}$
      such that $\mu_c(B_c) = 1$ for all $c \in \Cantorspace$.
  \end{enumerate}
\end{proposition}

\begin{propositionproof}
  This follows from Propositions \ref
  {measures:proposition:successor:sufficientcondition:countableunion} and \ref
  {measures:proposition:perfectsequence:sufficientcondition}.
\end{propositionproof}

Let \definedsymbol{\absolutelycontinuous[E,F]} denote the set of all $\pair{\mu}{\nu} \in
(\ergodicquasiinvariant{X}{E} \setminus \hyperfinite{X}{E}) \times (\ergodicquasiinvariant{Y}{F}
\setminus \hyperfinite{Y}{F})$ for which there is a $\mu$-conull \Borel set $C \subseteq X$ and a
\Borel reduction $\phi \from C \to Y$ of $\restriction{E}{C}$ to $F$ sending $(\restriction{\mu}
{C})$-positive sets to $\nu$-positive sets. Clearly $\absolutelycontinuous[E,F]$ is transitive,
and if $C \subseteq X$ is a $\mu$-conull \Borel set and $\phi \from C \to X$ is a \Borel
reduction of $\restriction{E}{C}$ to $F$, then $\mu \absolutelycontinuous[E,F] 
\pushforward{\phi}{(\restriction{\mu}{C})}$. When $E = F$, we simply write 
\definedsymbol{\absolutelycontinuous[E]}. It is easy to see this is an equivalence relation,
in spite of our adherence to the usual measure-theoretic abuse of notation.

The following fact provides a partial converse to Proposition \ref
{measures:proposition:successor:sufficientcondition}.

\begin{proposition} \label{measures:proposition:successor:characterization}
  Suppose that $X$ is a standard \Borel space, $E$ is a countable \Borel equivalence relation
  on $X$, and some vertical section of $\absolutelycontinuous[E]$ is a countable union of
  measure-equivalence classes. Then the following are equivalent:
  \begin{enumerate}
    \item The set $\ergodicquasiinvariant{X}{E} \setminus \hyperfinite{X}{E}$ is a single
      measure-equivalence class.
    \item The relation $E$ is a successor of $\Ezero$ under measure reducibility.
  \end{enumerate}
\end{proposition}

\begin{propositionproof}
  By Proposition \ref{measures:proposition:successor:sufficientcondition}, it is sufficient to show
  that if $\ergodicquasiinvariant{X}{E} \setminus \hyperfinite{X}{E}$ contains multiple
  measure-equivalence classes, then $E$ is not a successor of $\Ezero$ under measure
  reducibility. Towards this end, fix $\mu \in \ergodicquasiinvariant{X}{E} \setminus 
  \hyperfinite{X}{E}$ for which the corresponding vertical
  section of $\absolutelycontinuous[E]$ is a countable union of measure-equivalence classes, as
  well as a \Borel probability measure $\nu$ on $X$ in $\ergodicquasiinvariant{X}{E} \setminus
  \hyperfinite{X}{E}$ for which $\mu \not \measureequivalence \nu$. Fix an $E$-invariant
  $\nu$-conull \Borel set $D \subseteq X$ which is null with respect to every measure in the
  \textexponent{\mu}{th} vertical section of $\absolutelycontinuous[E]$ which is not measure
  equivalent to $\nu$.

  \begin{lemma}
    Suppose that $A \subseteq X \setminus D$ is a $\mu$-conull \Borel set and $B \subseteq
    D$ is a $\nu$-conull \Borel set. Then there is no \Borel reduction $\phi \from A \union B \to
    D$ of $\restriction{E}{(A \union B)}$ to $\restriction{E}{D}$.
  \end{lemma}

  \begin{lemmaproof}
    Suppose that $\phi$ is such a reduction. Then our choice of $D$ ensures that $\pushforward
    {(\restriction{\phi}{A})}{(\restriction{\mu}{A})} \absolutelycontinuous \nu$, so $\mu
    \absolutelycontinuous[E] \nu$. As $\absolutelycontinuous[E]$ is transitive, it follows that $\mu
    \absolutelycontinuous[E] \pushforward{(\restriction{\phi}{B})}{(\restriction{\nu}{B})}$, so
    our choice of $D$ also ensures that $\pushforward{(\restriction{\phi}{B})}{(\restriction{\nu}{B})}
    \absolutelycontinuous \nu$. Then there exist $x \in A$ and $y \in B$ such that $\phi(x) \mathrel{E}
    \phi(y)$, contradicting the fact that $\phi$ is a reduction.
  \end{lemmaproof}

  In particular, it follows that $E$ is not measure reducible to $\restriction{E}{D}$, and therefore
  cannot be a successor of $\Ezero$ under measure reducibility.
\end{propositionproof}

The following provides a partial converse to Proposition \ref
{measures:proposition:successor:sufficientcondition:countableunion}.

\begin{proposition} \label{measures:proposition:successor:characterization:countableunion}
  Suppose that $X$ is a standard \Borel space, $E$ is a countable \Borel equivalence relation
  on $X$, and every vertical section of $\absolutelycontinuous[E]$ is a countable union of
  measure-equivalence classes. Then the following are equivalent:
  \begin{enumerate}
    \item The set $\ergodicquasiinvariant{X}{E} \setminus \hyperfinite{X}{E}$ is a non-empty
      countable union of measure-equivalence classes.
    \item The relation $E$ is a non-empty countable disjoint union of successors of $\Ezero$
      under measure reducibility.
  \end{enumerate}
\end{proposition}

\begin{propositionproof}
  By Proposition \ref{measures:proposition:successor:sufficientcondition:countableunion}, it is
  sufficient to show that if $N$ is a non-empty countable set and $\sequence{B_n}[n \in N]$ is a
  partition of $X$ into $E$-invariant \Borel sets on which $E$ is a successor of $\Ezero$ under
  measure reducibility, then $\ergodicquasiinvariant{X}{E} \setminus \hyperfinite{X}{E}$ is a
  countable union of measure-equivalence classes. Towards this end, note that for all $n \in N$,
  every vertical section of $\absolutelycontinuous[\restriction{E}{B_n}]$ is a countable union of
  measure-equivalence classes, so Proposition \ref
  {measures:proposition:successor:characterization} ensures that $\ergodicquasiinvariant
  {B_n}{\restriction{E}{B_n}} \setminus \hyperfinite{B_n}{\restriction{E}{B_n}}$ is the
  measure-equivalence class of some \Borel probability measure $\mu_n$ on $B_n$.
  Identifying $\mu_n$ with the corresponding \Borel probability measure on $X$, it follows that
  $\ergodicquasiinvariant{X}{E} \setminus \hyperfinite{X}{E}$ is the union of the
  measure-equivalence classes of $\mu_n$, for $n \in N$.
\end{propositionproof}

Summarizing these results, we obtain the following.

\begin{theorem} \label{measures:theorem:trichotomy}
  Suppose that $X$ is a standard \Borel space, $E$ is a non-measure-hyperfinite countable
  \Borel equivalence relation on $X$, and every vertical section of $\absolutelycontinuous[E]$
  is a countable union of measure-equivalence classes. Then exactly one of the following holds:
  \begin{enumerate}
    \item The relation $E$ is a countable disjoint union of successors of $\Ezero$
      under measure reducibility.
    \item There are \Borel sequences $\sequence{B_c}[c \in \Cantorspace]$ of pairwise disjoint
      $E$-invar\-iant subsets of $X$ and $\sequence{\mu_c}[c \in \Cantorspace]$ of \Borel
      probability measures on $X$ in $\ergodicquasiinvariant{X}{E} \setminus \hyperfinite{X}{E}$
      such that $\mu_c(B_c) = 1$ for all $c \in \Cantorspace$.
  \end{enumerate}
\end{theorem}

\begin{theoremproof}
  If $\ergodicquasiinvariant{X}{E} \setminus \hyperfinite{X}{E}$ is a non-empty countable union
  of measure-equival\-ence classes, then Proposition \ref
  {measures:proposition:successor:characterization:countableunion} ensures that condition
  (2) holds, and its proof implies that condition (3) fails. If $\ergodicquasiinvariant{X}{E}
  \setminus \hyperfinite{X}{E}$ is not a countable union of measure-equivalence classes, then
  Proposition \ref{measures:proposition:successor:characterization:countableunion} ensures that
  condition (2) fails, and Proposition \ref{measures:proposition:perfectsequence:sufficientcondition}
  implies that condition (3) holds.
\end{theoremproof}

In light of our earlier results, the following yields a criterion for ensuring that \Borel
subequivalence relations of successors of $\Ezero$ under measure reducibility are again
successors of $\Ezero$ under measure reducibility.

\begin{proposition} \label{measures:proposition:containment}
  Suppose that $X$ is a standard \Borel space, $E \subseteq F$ are countable \Borel equivalence
  relations on $X$, $\mu$ is an $E$-ergodic $F$-quasi-invariant \Borel probability measure on
  $X$, and $\ergodicquasiinvariant{X}{F} \setminus \hyperfinite{X}{F}$ is contained in the
  measure-equivalence class of $\mu$. Then $\ergodicquasiinvariant{X}{E} \setminus
  \hyperfinite{X}{E}$ is also contained in the measure-equivalence class of $\mu$.
\end{proposition}

\begin{propositionproof}
  Suppose that $\nu \in \ergodicquasiinvariant{X}{E}$ but $\mu \not\measureequivalence \nu$.
  Then there is an $E$-invariant $\mu$-null $\nu$-conull \Borel set $C \subseteq X$, in which
  case Proposition \ref{preliminaries:measuredequivalencerelationsproposition:quasiinvariant}
  yields a \Borel probability measure $\nu' \reverseabsolutelycontinuous \nu$ with the same 
  $F$-invariant null sets. As the $F$-quasi-invariance of $\mu$ ensures that $\saturation{C}{F}$ is
  $\mu$-null, it follows that $\nu' \in \hyperfinite{X}{F}$. As Proposition \ref
  {preliminaries:hyperfiniteness:proposition:closure:homomorphisms} ensures that the class of
  hyperfinite \Borel equivalence relations is closed downward under \Borel subequivalence
  relations, it follows that $\nu' \in \hyperfinite{X}{E}$, thus $\nu \in \hyperfinite{X}{E}$.
\end{propositionproof}

We also have the following criterion for ensuring strong ergodicity.

\begin{proposition} \label{measures:proposition:E0ergodic}
  Suppose that $X$ and $Y$ are standard \Borel spaces, $E$ and $F$ are countable \Borel
  equivalence relations on $X$ and $Y$, $F$ is hyperfinite, and $\mu \in \ergodicquasiinvariant
  {X}{E} \setminus \hyperfinite{X}{E}$ is not $\pair{E}{F}$-ergodic. Then $\ergodicquasiinvariant
  {X}{E} \setminus \hyperfinite{X}{E}$ is not a countable union of measure-equivalence classes.
\end{proposition}

\begin{propositionproof}
  Fix a $\mu$-null-to-one \Borel homomorphism $\phi \from X \to Y$ from $E$ to $F$,
  as well as a \Borel disintegration $\sequence{\mu_y}[y \in Y]$ of $\mu$ through $\phi$.

  Then the set $C = \set{y \in Y}[E \text{ is not $\mu_y$-hyperfinite}]$ is \Borel by Theorem \ref
  {preliminaries:measurehyperfiniteness:theorem:uniform}. As $E$ is $\mu$-nowhere hyperfinite,
  Proposition \ref{preliminaries:measurehyperfiniteness:proposition:hyperfinite} ensures that $C$ is
  $(\pushforward{\phi}{\mu})$-conull.

  In particular, as $\phi$ is $\mu$-null-to-one, it follows that $C$ is uncountable, in which case
  there is an uncountable partial transversal $P \subseteq C$ of $F$. Theorem 
  \ref{preliminaries:measurehyperfiniteness:theorem:E0} then yields \Borel probability
  measures $\nu_y$ on $X$ in $\ergodicquasiinvariant{X}{E} \setminus \hyperfinite{X}{E}$ such
  that $\saturation{\preimage{\phi}{y}}{E}$ is $\nu_y$-conull, for all $y \in P$. As the latter sets are
  pairwise disjoint, it follows that $\ergodicquasiinvariant{X}{E} \setminus \hyperfinite{X}{E}$ is not
  a countable union of measure-equivalence classes.
\end{propositionproof}

We next compute a bound on the complexity of $\absolutelycontinuous[E, F]$.

\begin{proposition} \label{measures:proposition:analytic}
  Suppose that $X$ and $Y$ are standard \Borel spaces and $E$ and $F$ are countable \Borel
  equivalence relations on $X$ and $Y$. Then $\absolutelycontinuous[E, F]$ is analytic.
\end{proposition}

\begin{propositionproof}
  Note that $\mu \absolutelycontinuous[E, F] \nu$ if and only if there is a code $c$ for a
  measurable function $\phi_c \from X \to Y$ such that $\pushforward{(\codedfunction{c})}{(\restriction
  {\mu}{\domain{\codedfunction{c}}})} \absolutelycontinuous \nu$ and $\codedfunction{c}$ is a
  reduction of $E$ to $F$ on a $\mu$-conull set. Proposition \ref
  {preliminaries:measures:proposition:pushforward} ensures that the former relation is \Borel, and
  Proposition \ref{preliminaries:measuredequivalencerelations:proposition:analytic} implies that the
  latter relation is analytic.
\end{propositionproof}

We close this section by noting that our hypothesis on $\absolutelycontinuous[E]$ holds of all
projectively-separable countable \Borel equivalence relations.

\begin{proposition} \label{measures:proposition:countablesections}
  Suppose that $X$ and $Y$ are standard \Borel spaces, $E$ and $F$ are countable
  \Borel equivalence relations on $X$ and $Y$, and $F$ is projectively separable. Then the
  vertical sections of $\absolutelycontinuous[E, F]$ are countable unions of measure-equivalence
  classes.
\end{proposition}

\begin{propositionproof}
  Suppose that $\mu \in \ergodicquasiinvariant{X}{E} \setminus \hyperfinite{X}{E}$, and let $A$
  denote the vertical section of $\absolutelycontinuous[E, F]$ corresponding to $\mu$. As
  Proposition \ref{measures:proposition:analytic} ensures that $\absolutelycontinuous[E, F]$ is
  analytic, so too is $A$. As measure
  equivalence is \Borel, Theorem \ref{preliminaries:Borelequivalencerelations:theorem:perfect}
  implies that if $A$ is not a union of countably-many measure-equivalence classes, then there
  is a \Borel sequence $\sequence{\nu_c}[c \in \Cantorspace]$ of pairwise orthogonal \Borel
  probability measures on $Y$ in $A$. Theorem \ref{preliminaries:measures:theorem:perfect} then
  ensures that by passing to an appropriate subsequence, we can ensure that there is a \Borel
  sequence $\sequence{D_c}[c \in \Cantorspace]$ of pairwise disjoint subsets of $Y$ such that
  $\nu_c(D_c) = 1$ for all $c \in \Cantorspace$. But for each $c \in \Cantorspace$, there is a
  $\mu$-conull \Borel set $C_c \subseteq X$ for which there is a \Borel reduction $\phi_c \from
  C_c \to D_c$ from $\restriction{E}{C_c}$ to $\restriction{F}{D_c}$, contradicting the projective
  separability of $F$.
\end{propositionproof}

\section{Stratification} \label{stratification}

Proposition \ref{preliminaries:hyperfiniteness:proposition:aperiodic} ensures that every aperiodic
countable \Borel equivalence relation has an aperiodic hyperfinite \Borel subequivalence relation.
This is the special case of the following fact, in which $G$ is the difference of $E$ and equality,
and $\rho$ is the constant cocycle.

\begin{proposition} \label{stratification:proposition:hyperfinitesubequivalencerelation:aperiodic}
  Suppose that $X$ is a standard \Borel space, $E$ is a countable \Borel equivalence relation on
  $X$, $G$ is a \Borel graphing of $E$, and $\rho \from E \to \Rplus$ is an aperiodic \Borel cocycle.
  Then there is a \Borel subgraph $H$ of $G$ generating a hyperfinite \Borel equivalence relation
  on which $\rho$ is aperiodic.
\end{proposition}

\begin{propositionproof}
  As graphs of \Borel functions are themselves \Borel, the following observation implies that it is
  sufficient to establish the proposition on an $E$-complete \Borel set.

  \begin{lemma} \label{stratification:lemma:completeextension}
    Suppose that $B \subseteq X$ is an $E$-complete \Borel set and $H$ is a \Borel subgraph of
    $\restriction{G}{B}$ generating a hyperfinite \Borel equivalence relation on which $\rho$ is
    aperiodic. Then there is a \Borel function $f \from \setcomplement{B} \to X$ such that $\graph
    {f^{\pm 1}} \union H$ is a subgraph of $G$ generating a hyperfinite \Borel equivalence relation
    on which $\rho$ is aperiodic.
  \end{lemma}

  \begin{lemmaproof}
    As the vertical sections of $G$ are countable, the \Lusin-\Novikov uniformization theorem
    yields \Borel sets $B_n \subseteq X$ and \Borel
    functions $f_n \from B_n \to X$ with the property that $G = \union[n \in \N][\graph{f_n}]$. Let
    $d_G(x, B)$ denote the length of the shortest $G$-path from $x$ to an element of $B$. Observe
    that this function is \Borel, as it can also be expressed, for $x \notin B$, as the least $n \in \N$ for
    which there exist $k_1, \ldots, k_n \in \N$ such that $f_{k_1} \composition \cdots \composition
    f_{k_n}(x) \in B$. Noting that for each $x \in X$, the set of $y \in \verticalsection{G}{x}$ with the
    property that $d_G(y, B) = d_G(x, B) - 1$ is countable, one more application of the
    \Lusin-\Novikov uniformization theorem yields a \Borel function $f \from \setcomplement{B}
    \to X$, whose graph is contained in $G$, such that $d_G(f(x), B) = d_G(x, B) - 1$ for all $x \in
    \setcomplement{B}$. As every connected component of $\graph{f^{\pm 1}} \union H$ contains
    a connected component of $H$, it follows that $\rho$ is aperiodic on the equivalence
    relation generated by $\graph{f^{\pm 1}} \union H$. As the function sending $x$ to $f^{d_G(x, B)}
    (x)$ is a \Borel reduction of the latter equivalence relation to that generated by $H$,
    and Proposition \ref{preliminaries:hyperfiniteness:proposition:closure:homomorphisms} ensures
    that the class of hyperfinite \Borel equivalence relations is closed under \Borel reducibility, it
    follows that the equivalence relation generated by $\graph{f^{\pm 1}} \union H$ is hyperfinite.
  \end{lemmaproof}

  We will now recursively construct an increasing sequence $\sequence{H_n}[n \in \N]$ of
  approximations to the desired graph, beginning with $H_0 = \emptyset$. Given $H_n$, let $E_n$
  denote the equivalence relation induced by $H_n$, and let $B_n$ denote the set of all $x \in X$
  for which $\rho$ is finite on $\restriction{E_n}{\equivalenceclass{x}{E}}$. As $H_n$ and $E$ are
  countable, the \Lusin-\Novikov uniformization theorem ensures that these sets are
  \Borel. As Proposition \ref{preliminaries:measuredequivalencerelationsproposition:periodic} implies
  that countable \Borel equivalence relations admitting finite \Borel cocycles to $\R$ are smooth, it
  follows that $\restriction{E_n}{B_n}$ is smooth. Remark \ref
  {preliminaries:countableBorelequivalencerelations:remark:smooth:countable} therefore yields a
  \Borel transversal $A_n \subseteq B_n$ of $\restriction{E_n}{B_n}$. Let $R_n$ be the
  relation consisting of all $\pair{x}{\pair{y}{\pair{x'}{y'}}} \in A_n \times (A_n \times (B_n \times
  B_n))$ for which $x \mathrel{E_n} x' \mathrel{(G \setminus E_n)} y' \mathrel{E_n} y$ and $\rho
  (\equivalenceclass{x}{E_n}, \equivalenceclass{y}{E_n}) \le 1$. As the vertical sections of $R_n$
  are countable, the \Lusin-\Novikov uniformization theorem ensures that the set
  $A_n' = \image{\projection[A_n]}{R_n}$ is \Borel, there is a \Borel uniformization $f_n' \from A_n'
  \to A_n \times (B_n \times B_n)$ of $R_n$, and both of the sets $S_n = \image{f_n'}{A_n'}$ and
  $H_n' = \image{\projection[B_n \times B_n]}{S_n}^{\pm 1}$ are \Borel.

  Set $H_{n+1} = H_n \union H_n'$. To see that the equivalence relation $E_{n+1}$ generated
  by $H_{n+1}$ is hyperfinite, we consider the function $f_n \from A_n \to A_n$ given by $f_n =
  (\projection[A_n] \composition f_n') \union (\restriction{\id}{(A_n \setminus A_n')})$. As Theorem
  \ref{preliminaries:hyperfiniteness:theorem:tailrelation} ensures that $\tailequivalencerelation
  {f_n}$ is hypersmooth, Theorem \ref{preliminaries:hyperfiniteness:theorem:hypersmooth} implies
  that it is hyperfinite. As Proposition \ref
  {preliminaries:hyperfiniteness:proposition:closure:homomorphisms} ensures that the class of
  hyperfinite \Borel equivalence relations is closed downward under \Borel reducibility, and the
  unique function $\phi_n \from B_n \to A_n$ such that $\forall x \in B_n \ x \mathrel{E_n}
  \phi_n(x)$ is a \Borel reduction of $\restriction{E_{n+1}}{B_n}$ to $\tailequivalencerelation{f_n}$,
  it follows that $E_{n+1}$ is hyperfinite. This completes the recursive construction.

  As every equivalence class of $\restriction{E}{\setcomplement{B_n}}$ contains a $\rho$-infinite
  equivalence class of $E_n$, it follows from Lemma \ref{stratification:lemma:completeextension}
  that we can construct the desired graph off of the set $B_\infty = \intersection[n \in \N][B_n]$. In
  order to construct the desired graph on $B_\infty$, set $H_\infty = \union[n \in \N][H_n]$ and let
  $E_\infty$ denote the equivalence relation generated by $H_\infty$. As $E_\infty = \union[n \in
  \N][E_n]$, it follows that $\restriction{E_\infty}{B_\infty}$ is hypersmooth, so Theorem \ref
  {preliminaries:hyperfiniteness:theorem:hypersmooth} ensures that it is hyperfinite. By one more
  application of Lemma \ref{stratification:lemma:completeextension}, it is therefore sufficient to
  observe that there do not exist $(G \setminus E_\infty)$-related points $x, y \in B_\infty$ for
  which the corresponding equivalence classes $\equivalenceclass{x}{E_\infty},
  \equivalenceclass{y}{E_\infty}$ are $\rho$-finite.

  Suppose, towards a contradiction, that there are such points. Then there exists $n \in \N$
  such that $\rho(\equivalenceclass{x}{E_\infty}, \equivalenceclass{x}{E_n}), \rho
  (\equivalenceclass{y}{E_\infty}, \equivalenceclass{y}{E_n}) < 2$. As $\rho(\equivalenceclass{x}
  {E_n}, \equivalenceclass{y}{E_n}) \le 1$ or $\rho(\equivalenceclass{y}{E_n}, \equivalenceclass
  {x}{E_n}) \le 1$, it follows that $\phi_n(x) \in A_n'$ or $\phi_n(y) \in A_n'$, so $\rho
  (\equivalenceclass{x}{E_{n+1}}, \equivalenceclass{x}{E_n}) \ge 2$ or $\rho(\equivalenceclass
  {y}{E_{n+1}}, \equivalenceclass{y}{E_n}) \ge 2$, thus $\rho(\equivalenceclass{x}{E_\infty},
  \equivalenceclass{x}{E_{n+1}}) < 1$ or $\rho(\equivalenceclass{y}{E_\infty}, \equivalenceclass
  {y}{E_{n+1}}) < 1$, which is impossible.
\end{propositionproof}

In particular, we obtain the following measure-theoretic corollary.

\begin{proposition} \label{stratification:proposition:hyperfinitesubequivalencerelation:nonsmooth}
  Suppose that $X$ is a standard \Borel space, $E$ is a countable \Borel equivalence relation on
  $X$, $\mu$ is an $E$-quasi-invar\-iant \Borel probability measure on $X$ for which $E$ is
  $\mu$-nowhere smooth, and $G$ is a \Borel graphing of $E$. Then there is a \Borel subgraph
  $H$ of $G$ whose induced equivalence relation is $\mu$-nowhere smooth but hyperfinite.
\end{proposition}

\begin{propositionproof}
  By Proposition \ref{preliminaries:measuredequivalencerelationsproposition:cocycles}, there is a
  \Borel cocycle $\rho \from E \to \Rplus$ with respect to which $\mu$ is invariant. As Proposition
  \ref{preliminaries:measuredequivalencerelationsproposition:periodic} ensures that countable \Borel
  equivalence relations admitting finite \Borel cocycles to $\R$ are smooth, by throwing away an
  $E$-invariant $\mu$-null \Borel set on which $E$ is smooth, we can assume that $\rho$ is
  aperiodic. Proposition \ref{stratification:proposition:hyperfinitesubequivalencerelation:aperiodic}
  then yields a \Borel subgraph $H$ of $G$ generating a hyperfinite equivalence relation on
  which $\rho$ is aperiodic. As Proposition \ref
  {preliminaries:measuredequivalencerelationsproposition:aperiodic} ensures that every such relation
  is $\mu$-nowhere smooth, the result follows.
\end{propositionproof}

The following yields disjoint \Borel sets which, in the measure-theoret\-ic
setting, are complete with respect to different equivalence relations.

\begin{proposition} \label{stratification:proposition:completesection}
  Suppose that $X$ is a standard \Borel space, $E$ and $F$ are aperiodic countable \Borel
  equivalence relations on $X$, and $\mu$ and $\nu$ are \Borel probability measures on $X$.
  Then there are disjoint \Borel sets $A, B \subseteq X$ such that $\mu(\saturation{A}{E}) = \nu
  (\saturation{B}{F}) = 1$.
\end{proposition}

\begin{propositionproof}
  By two applications of Proposition \ref
  {preliminaries:countableBorelequivalencerelations:proposition:markers}, there are decreasing
  sequences $\sequence{A_n}[n \in \N]$ and $\sequence{B_n}[n \in \N]$ of \Borel subsets of $X$
  such that each $A_n$ is $E$-complete, each $B_n$ is $F$-complete, and $\intersection[n \in \N]
  [A_n] = \intersection[n \in \N][B_n] = \emptyset$. Fix real numbers $\epsilon_n > 0$ such that
  $\epsilon_n \goesto 0$ as $n \goesto \infty$, and recursively construct strictly increasing
  sequences $\sequence{i_n}[n \in \N]$ and $\sequence{j_n}[n \in \N]$ of natural numbers by
  setting $i_0 = 0$, and given $n \in \N$ and $i_n \in \N$, choosing $j_n > \max_{m < n} j_m$
  sufficiently large that $\mu(\saturation{A_{i_n} \setminus B_{j_n}}{E}) \ge 1 - \epsilon_n$, as well
  as $i_{n+1} > i_n$ sufficiently large that $\nu(\saturation{B_{j_n} \setminus A_{i_{n+1}}}{F}) \ge 1
  - \epsilon_n$. Define $A = \union[n \in \N][(A_{i_n} \setminus B_{j_n})]$ and $B = \union
  [n \in \N][(B_{j_n} \setminus A_{i_{n+1}})]$.
\end{propositionproof}

A \definedterm{directed graph} on $X$ is an irreflexive subset $G$ of $X \times X$. The
\definedterm{domain} of such a relation is the set of $x$ for which $\verticalsection{G}{x}$ is
non-empty. An \definedterm{oriented graph} on $X$ is an irreflexive antisymmetric subset $H$ of
$X \times X$. An \definedterm{orientation} of a graph $G$ is an oriented graph $H$ with
$G = H^{\pm 1}$. Although the domain of an orientation $H$ of a graph $G$ can be strictly
smaller than the domain of $G$ itself, we do have the following.

\begin{proposition} \label{stratification:proposition:orientation}
  Suppose that $X$ is a standard \Borel space, $E$ is an aperiodic countable \Borel equivalence
  relation on $X$, $G$ is a locally countable \Borel graph on $X$, and $\mu$ is an
  $E$-quasi-invariant \Borel probability measure on $X$ for which $E$ is $\mu$-nowhere smooth
  and the domain of $G$ has $\mu$-conull $E$-saturation. Then there is a \Borel orientation $H$
  of $G$ whose domain has $\mu$-conull $E$-saturation.
\end{proposition}

\begin{propositionproof}
  For each \Borel set $B \subseteq X$, put $X_B = \set{x \in B}[\equivalenceclass{x}{\restriction
  {E}{B}} \text{ is finite}]$. As $E$ is countable, the \Lusin-\Novikov uniformization theorem ensures
  that such sets are \Borel, as are $E$-saturations of \Borel sets.

  \begin{lemma} \label{lemma:periodicpart}
    Suppose that $B \subseteq X$ is \Borel. Then $\saturation{X_B}{E}$ is $\mu$-null.
  \end{lemma}

  \begin{lemmaproof}
    As $E$ is countable, the \Lusin-\Novikov uniformization theorem ensures that
    there is a \Borel reduction of $\restriction{E}{\saturation{X_B}{E}}$ to $\restriction{E}{X_B}$. As
    Proposition \ref{preliminaries:countableBorelequivalencerelations:proposition:smooth:finite}
    ensures that $\restriction{E}{X_B}$ is smooth, so too is $\restriction{E}{\saturation{X_B}{E}}$.
    As $E$ is $\mu$-nowhere smooth, it follows that $\saturation{X_B}{E}$ is $\mu$-null.
  \end{lemmaproof}

  We consider now the special case that $G$ is of the form $\graph{I}$, where $A \subseteq X$
  is a \Borel set and $I \from A \to A$ is a \Borel involution. Proposition \ref
  {preliminaries:countableBorelequivalencerelations:proposition:smooth:finite} and Remark \ref
  {preliminaries:countableBorelequivalencerelations:remark:smooth:countable} yield a \Borel
  transversal $B \subseteq A$ of the equivalence relation generated by $G$. Lemma \ref
  {lemma:periodicpart} ensures that the set $C = \saturation{A}{E} \setminus \saturation{X_B
  \union X_{A \setminus B}}{E}$ is $\mu$-conull.

  We use $E_B$, $E_{A \setminus B}$, $\mu_B$, and $\mu_{A \setminus B}$ to denote the
  restrictions of $E$, $\preimage{(I \times I)}{E}$, $\mu$, and $\pushforward{I}{\mu}$ to $B
  \intersection C$. As $E_B$ and $E_{A \setminus B}$ are aperiodic, Proposition \ref
  {stratification:proposition:completesection} yields a \Borel set $B' \subseteq B$, an
  $E_B$-invariant $\mu_B$-null \Borel set $N_B \subseteq C$, and an $E_{A \setminus
  B}$-invariant $\mu_{A \setminus B}$-null \Borel set $N_{A \setminus B} \subseteq C$ such that
  $B' \union N_B$ is $E_B$-complete and $(B \setminus B') \union N_{A \setminus B}$ is $E_{A
  \setminus B}$-complete. As $\mu$ is $E$-quasi-invariant, the set $D = C \setminus \saturation
  {N_B \union N_{A \setminus B}}{E}$ is $\mu$-conull. Let $H$ denote the graph of the restriction
  of $I$ to $B' \union \image{I}{B \setminus B'}$.

  The fact that $B$ is a transversal of the equivalence relation generated by $I$ ensures that $H$
  is an oriented graph. To see that $H$ is an orientation of $G$, note that if $x \mathrel{G} y$, then
  $x \in B$ or $y \in B$, from which it follows that $(x \in B' \mathor y \in \image{I}{B \setminus B'})$
  or $(y \in B' \mathor x \in \image{I}{B \setminus B'})$, so $(x \mathrel{H} y \mathor y \mathrel{H}
  x)$ or $(y \mathrel{H} x \mathor x \mathrel{H} y)$, thus $x \mathrel{H} y$ or $y \mathrel{H} x$. To
  see that the $E$-saturation of the domain of $H$ is $\mu$-conull, it is enough to show that the
  domain of $H$ intersects the $E$-class of every $x \in D$. Towards this end, note that $A
  \intersection \equivalenceclass{x}{E}$ is non-empty, thus so too is $B \intersection
  \equivalenceclass{x}{E}$ or $(A \setminus B) \intersection \equivalenceclass{x}{E}$, in which
  case $B' \intersection \equivalenceclass{x}{E}$ or $\image{I}{B \setminus B'} \intersection
  \equivalenceclass{x}{E}$ is non-empty as well, hence the domain of $H$ intersects
  $\equivalenceclass{x}{E}$.

  We now consider the general case. As $G$ is locally countable, Theorem \ref
  {preliminaries:countableBorelequivalencerelations:theorem:involutions} yields \Borel sets $A_n \subseteq X$ and \Borel
  involutions $I_n \from A_n \to A_n$, with pairwise disjoint graphs, such that $G = \union[n \in \N]
  [\graph{I_n}]$. Setting $G_n = \graph{I_n}$, $X_n = \saturation{A_n}{E}$, and $\mu_n =
  \restriction{\mu}{X_n}$, the above special case yields \Borel orientations $H_n$ of $G_n$ whose
  domains have $\mu_n$-conull $E$-saturations. Then $H = \union[n \in \N][H_n]$ is a \Borel
  orientation of $G$ whose domain has $\mu$-conull $E$-saturation.
\end{propositionproof}

A \definedterm{$\mu$-stratification} of $E$ is an increasing sequence $\sequence{E_r}[r \in \R]$
of subequivalence relations of $E$ whose union is $E$ and which is strictly increasing on every
$\mu$-positive \Borel set.

\begin{theorem} \label{stratification:theorem:sufficientcondition}
  Suppose that $X$ is a standard \Borel space, $E$ is a treeable countable \Borel equivalence
  relation on $X$, and $\mu$ is an $E$-quasi-invariant \Borel probability measure on $X$ for
  which $E$ is $\mu$-nowhere hyperfinite. Then there is a \Borel $\mu$-stratification of $E$.
\end{theorem}

\begin{theoremproof}
  Fix a \Borel treeing $G$ of $E$. By Proposition \ref
  {stratification:proposition:hyperfinitesubequivalencerelation:nonsmooth}, we can assume
  that there is \Borel subgraph $H$ of $G$ whose induced equivalence relation $F$ is
  $\mu$-nowhere smooth but hyperfinite. As $E$ is $\mu$-nowhere hyperfinite, the $F$-saturation
  of the domain of $G \setminus H$ is $\mu$-conull. As $F$ is $\mu$-nowhere smooth, Proposition
  \ref{stratification:proposition:orientation} ensures that there is a \Borel orientation $K$ of $G
  \setminus H$ whose domain has $\mu$-conull $F$-saturation. As $\mu$ is $E$-quasi-invariant,
  by throwing out an $E$-invariant $\mu$-null \Borel set, we can assume that the domain of $K$
  intersects every $F$-class. By Proposition \ref
  {preliminaries:measuredequivalencerelationsproposition:cocycles}, there is a \Borel cocycle $\rho
  \from E \to \Rplus$ with respect to which $\mu$ is invariant. As $F$ is $\mu$-nowhere smooth
  and Proposition \ref{preliminaries:measuredequivalencerelationsproposition:periodic} ensures that
  $F$ is smooth on the finite part of $\restriction{\rho}{(\restriction{F}{\domain{K}})}$, by throwing
  out another $\mu$-null \Borel set, we can assume that $\restriction{\rho}{(\restriction{F}{\domain
  {K}})}$ is aperiodic, and therefore that $\restriction{F}{\domain{K}}$ is aperiodic. The
  $E$-quasi-invariance of $\mu$ again allows us to ensure that the set we throw out is
  $E$-invariant. Proposition \ref
  {preliminaries:countableBorelequivalencerelations:proposition:markers:disjoint} then yields a
  partition of the domain of $K$ into a sequence $\sequence{B_q}[q \in \Q]$ of pairwise disjoint
  $F$-complete \Borel sets. Set $K_r = \restriction{K}{(\union[q < r][B_q] \times X)}$ and
  $G_r = H \union K_r^{\pm 1}$ for all $r \in \R$. As $G_r$ is locally countable, the \Lusin-\Novikov
  uniformization theorem ensures that the equivalence relations $E_r$ induced by the graphs $G_r$
  are \Borel.

  Suppose now that $B \subseteq X$ is a \Borel set for which there are real numbers $r < s$ with
  $\restriction{E_r}{B} = \restriction{E_s}{B}$. Then $B \intersection \equivalenceclass{x}{E_s}
  \subseteq \equivalenceclass{x}{E_r}$ for all $x \in B$. As $G_r \subseteq G_s$ and the latter
  graph is acyclic, it follows that if $x \in B$ and $y \in \equivalenceclass{x}{E_s} \setminus
  \equivalenceclass{x}{E_r}$, then there is a unique point of $\equivalenceclass{y}{E_r}$ of minimal
  distance to $\equivalenceclass{x}{E_r}$ with respect to the graph metric associated with $G_s$.
  Let $\phi \from \saturation{B}{E_s} \setminus \saturation{B}{E_r} \to \saturation{B}{E_s} \setminus
  \saturation{B}{E_r}$ be the function sending each point of its domain to the corresponding point
  of its $E_r$-class. As $E$ is countable, the \Lusin-\Novikov uniformization theorem ensures that
  $\saturation{B}{E_r}$, $\saturation
  {B}{E_s}$, and $\phi$ are \Borel. As $\phi$ is a selector for the restriction of $E_r$ to
  $\saturation{B}{E_s} \setminus \saturation{B}{E_r}$, it follows that this restriction is smooth. As
  $F$ is $\mu$-nowhere smooth and Proposition \ref
  {preliminaries:countableBorelequivalencerelations:proposition:smooth:closure} ensures that the
  class of smooth countable \Borel equivalence relations is closed downward under \Borel
  subequivalence relations, it follows that $E_r$ is also $\mu$-nowhere smooth. In particular, this
  means that the set $\saturation{B}{E_s} \setminus \saturation{B}{E_r}$ is $\mu$-null, and since
  $\mu$ is $E_s$-quasi-invariant, so too is the $E_s$-saturation of $\saturation{B}{E_s} \setminus
  \saturation{B}{E_r}$. As every $E_r$-class is properly contained in the corresponding
  $E_s$-class, it follows that $B$ is contained in this saturation, and is therefore $\mu$-null as well,
  hence $\sequence{E_r}[r \in \R]$ is indeed a $\mu$-stratification of $E$.
\end{theoremproof}

\part{Applications}

Here we obtain our main results. While our theorems were listed in
the introduction in order of importance, we now proceed according to the amount of new
machinery behind the arguments, with those requiring the least appearing first. In \S\ref
{products}, we use the countability of the vertical sections of $\absolutelycontinuous[E]$
to establish our results on products. In \S\ref{embeddability}, we combine the countability
of the vertical sections of $\absolutelycontinuous[E]$ with facts about compressibility and costs of
equivalence relations to obtain our results on the distinction between embeddability and
reducibility. In \S\ref{antichains}, we combine projective separability, facts about
$\absolutelycontinuous[E]$, and the existence of stratifications to obtain our results on
antichains and the distinction between containment and reducibility. In \S\ref {bases}, we use
these tools to obtain our anti-basis theorems. And in \S\ref{complexity}, we combine these tools
with Theorem \ref{preliminaries:complexity:theorem:antichain} to obtain our complexity results.

\section{Products} \label{products}

We begin this section with an observation concerning measurable reducibility of products.

\begin{proposition} \label{products:proposition:continuum:measurable}
  Suppose that $X$ and $Y$ are standard \Borel spaces, $E$ and $F$ are countable \Borel
  equivalence relations on $X$ and $Y$, $m$ is a continuous \Borel probability measure on
  $\R$, $\mu \in \ergodicquasiinvariant{X}{E} \setminus \hyperfinite{X}{E}$, and the
  \textexponent{\mu}{th} vertical section of $\absolutelycontinuous[E, F]$ is a countable union of
  measure-equivalence classes. Then $E \times \diagonal{\R}$ is $(\mu \times m)$-nowhere
  reducible to $F$.
\end{proposition}

\begin{propositionproof}
  Suppose, towards a contradiction, that there is a $(\mu \times m)$-positive \Borel set $B
  \subseteq X \times \R$ on which there is a \Borel reduction $\phi \from B \to Y$ of $E \times
  \diagonal{\R}$ to $F$. As $E$ is countable, the \Lusin-\Novikov uniformization theorem
  ensures that $\saturation{B}{E \times \diagonal
  {\R}}$ is \Borel, in addition to yielding a \Borel reduction of $\restriction{(E \times \diagonal{\R})}
  {\saturation{B}{E \times \diagonal{\R}}}$ to $\restriction{(E \times \diagonal{\R})}{B}$. By
  replacing $B$ with its $(E \times \diagonal{\R})$-saturation, we can therefore assume that
  $B$ is $(E \times \diagonal{\R})$-invariant. Note that the set $R = \set{r \in \R}[\mu
  (\horizontalsection{B}{r}) > 0]$ is $m$-positive, by \Fubini's theorem. As $m$ is continuous,
  it follows that $R$ is uncountable. For each $r \in R$, Proposition \ref
  {preliminaries:measuredequivalencerelationsproposition:quasiinvariant} yields an
  $F$-quasi-invariant \Borel probability measure $\nu_r$ on $Y$ such that $\pushforward
  {(\horizontalsection{\phi}{r})}{(\restriction{\mu}{\horizontalsection{B}{r}})} \absolutelycontinuous \nu_r$, but the two
  measures have the same $F$-invariant null sets. But then the $\nu_r$ are
  pairwise orthogonal elements of the \textexponent{\mu}{th} vertical section of $\absolutelycontinuous[E, F]$, 
  the desired contradiction.
\end{propositionproof}

This has the following consequences for measure reducibility.

\begin{proposition} \label{products:proposition:continuum:measure}
  Suppose that $X$ and $Y$ are standard \Borel spaces, $E$ is a non-measure-hyperfinite
  countable \Borel equivalence relation on $X$, and $F$ is a projectively separable
  countable \Borel equivalence relation on $Y$. Then $E \times \diagonal{\R}$ is 
  not measure reducible to $F$.
\end{proposition}

\begin{propositionproof}
  By Theorem \ref{preliminaries:measurehyperfiniteness:theorem:E0}, there exists
  $\mu \in \ergodicquasiinvariant{X}{E} \setminus \hyperfinite{X}{E}$.
  Proposition \ref{measures:proposition:countablesections} then implies that the \textexponent
  {\mu}{th} vertical section of $\absolutelycontinuous[E,F]$ is a countable union of
  measure-equivalence classes. Fix a continuous \Borel probability measure $m$ on $\R$. As
  Proposition \ref{products:proposition:continuum:measurable} ensures that $E \times
  \diagonal{\R}$ is $(\mu \times m)$-nowhere reducible to $F$, the former is not measure
  reducible to the latter.
\end{propositionproof}

\begin{theorem}[\Hjorth]
  There is a non-measure-hyperfinite treeable countable \Borel equivalence relation to which
  some treeable countable \Borel equivalence relation is not measure reducible.
\end{theorem}

\begin{theoremproof}
  Proposition \ref{products:proposition:continuum:measure} ensures that every
  non-measure-hyperfinite pro\-jectively-separable treeable countable \Borel
  equivalence relation has the desired property.
\end{theoremproof}

We now consider products with smaller equivalence relations.

\begin{proposition} \label{products:proposition:two:measurable}
  Suppose that $X$ and $Y$ are standard \Borel spaces, $E$ and $F$ are countable \Borel
  equivalence relations on $X$ and $Y$, $m$ is a strictly positive probability measure on $2$,
  $\mu \in \ergodicquasiinvariant{X}{E} \setminus \hyperfinite{X}{E}$, $\nu \in
  \ergodicquasiinvariant{Y}{F} \setminus \hyperfinite{Y}{F}$, and $\mu \absolutelycontinuous[E, F]
  \nu$. If the \textexponent{\mu}{th} vertical section of $\absolutelycontinuous[E, F]$ is a countable
  union of measure-equivalence classes, then there is an $F$-invariant $\nu$-conull \Borel set
  $C \subseteq Y$ for which $E \times \diagonal{2}$ is not $(\mu \times m)$-reducible to
  $\restriction{F}{C}$.
\end{proposition}

\begin{propositionproof}
  Fix an $F$-invariant $\nu$-conull \Borel set $C \subseteq Y$ which is $\nu'$-null for every
  measure $\nu'$ in the vertical section of $\absolutelycontinuous[E, F]$ corresponding to
  $\mu$, other than those which are measure equivalent to $\nu$. Suppose, towards a
  contradiction, that there is a $(\mu \times m)$-positive \Borel set $B \subseteq X \times 2$ on
  which there is a \Borel reduction $\phi \from B \to Y$ of $E \times \diagonal{2}$ to $\restriction
  {F}{C}$. As $E$ is countable, the \Lusin-\Novikov uniformization theorem ensures
  that $\saturation{B}{E \times \diagonal{2}}$ is \Borel, in addition to yielding a \Borel reduction of
  $\restriction{(E \times \diagonal{2})}{\saturation{B}{E \times \diagonal{2}}}$ to $\restriction{(E
  \times \diagonal{2})}{B}$. By replacing $B$ with its $(E \times \diagonal{2})$-saturation, we can
  therefore assume that $B$ is $(E \times \diagonal{2})$-invariant. Proposition \ref
  {preliminaries:measuredequivalencerelationsproposition:quasiinvariant} then yields $(\restriction{F}
  {C})$-quasi-invariant \Borel probability measures $\nu_i$ on $C$ with the property that
  $\pushforward{(\horizontalsection{\phi}{i})}{(\restriction{\mu}{\horizontalsection{B}{i}})} \absolutelycontinuous \nu_i$
  but the two measures have the same $E$-invariant null \Borel sets, for all $i < 2$. As $\nu_0$
  and $\nu_1$ are orthogonal elements of the vertical section of $\absolutelycontinuous[E,
  \restriction{F}{C}]$, this contradicts our choice of $C$.
\end{propositionproof}

This has the following consequence for measure reducibility.

\begin{proposition}
  Suppose that $X$ is a standard \Borel space and $E$ is a non-measure-hyperfinite 
  projectively-separable countable \Borel equivalence relation on $X$. Then there is a \Borel set
  $B \subseteq X$ on which $E$ is not measure hyperfinite such that $(\restriction{E}{B}) \times
  \diagonal{2}$ is not measure reducible to $\restriction{E}{B}$.
\end{proposition}

\begin{propositionproof}
  By Theorem \ref{preliminaries:measurehyperfiniteness:theorem:E0}, there exists $\mu \in
  \ergodicquasiinvariant{X}{E} \setminus \hyperfinite{X}{E}$. Proposition \ref
  {measures:proposition:countablesections} then implies that the \textexponent{\mu}{th} vertical
  section of $\absolutelycontinuous[E]$ is a countable union of measure-equivalence classes. Fix
  a strictly positive probability measure $m$ on $2$. As Proposition \ref
  {products:proposition:two:measurable} yields a $\mu$-conull \Borel set $C \subseteq
  X$ for which $E \times \diagonal{2}$ is not $(\mu \times m)$-reducible to $\restriction{E}{C}$, it
  follows that $(\restriction{E}{C}) \times \diagonal{2}$ is not measure reducible to $\restriction{E}
  {C}$.
\end{propositionproof}

\begin{remark}
  A similar argument can be used to show that if $X$ and $Y$ are standard \Borel spaces, $E$ is
  a non-measure-hyperfinite countable \Borel equivalence relation on $X$, and $F$ is a
  non-measure-hyperfinite projectively-separable countable \Borel equivalence relation on $Y$,
  then there is a \Borel set $B \subseteq Y$ on which $F$ is not measure-hyperfinite such that
  $E \times \diagonal{2}$ is not measure reducible to $\restriction{F}{B}$.
\end{remark}

\section{Reducibility without embeddability} \label{embeddability}

We begin this section with an observation concerning the relationship between measurable
reducibility and measurable embeddability.

\begin{proposition} \label{embeddability:proposition:sufficientcondition:measurable}
  Suppose that $X$ and $Y$ are standard \Borel spaces, $E$ is an invariant-measure-hyperfinite
  countable \Borel equivalence relation on $X$, $F$ is an aperiodic countable \Borel equivalence
  relation $Y$, and $\mu$ is a \Borel probability measure on $X$. Then $E$ is $\mu$-reducible to
  $F$ if and only if $E$ is $\mu$-embeddable into $F$.
\end{proposition}

\begin{propositionproof}
  Suppose that $E$ is $\mu$-reducible to $F$, and fix a $\mu$-conull \Borel set $C \subseteq X$
  on which there is a \Borel reduction $\phi \from C \to Y$ of $E$ to $F$.  As $E$ is countable,
  the \Lusin-\Novikov uniformization theorem ensures that $\saturation{C}{E}$ is
  \Borel, and there is a \Borel reduction of $\restriction{E}{\saturation{C}{E}}$ to $\restriction{E}
  {C}$. By replacing $\phi$ with its composition with such a function, we can therefore assume
  that $C$ is itself $E$-invariant. Proposition \ref
  {preliminaries:measuredequivalencerelationsproposition:quasiinvariant} ensures that there is an
  $E$-quasi-invariant \Borel probability measure on $X$, with respect to which $\mu$ is
  absolutely continuous, which agrees with $\mu$ on all $E$-invariant \Borel sets. By replacing
  $\mu$ with such a measure, we can assume that $\mu$ is $E$-quasi-invariant.

  We handle first the case that $F$ is smooth. Then $\restriction{E}{C}$ is also smooth. As $E$ is
  countable, Remark \ref{preliminaries:countableBorelequivalencerelations:remark:smooth:countable}
  yields partitions $\sequence{C_n}[n \in \N]$ of $C$ into \Borel partial transversals of $E$, and
  $\sequence{D_n}[n \in \N]$ of $Y$ into \Borel transversals of $F$. One then obtains an
  embedding $\pi \from C \to Y$ of $\restriction{E}{C}$ into $F$ by setting
  \begin{equation*}
    \pi(x) = y \iff \exists n \in \N \ ( x \in C_n \mathcomma y \in D_n \mathcommaand \phi(x) \mathrel
      {F} y ).
  \end{equation*}
  As $C$ inherits a standard \Borel structure from $X$ and functions between standard
  \Borel spaces are \Borel if and only if their graphs are \Borel, it follows that $\pi$ is \Borel.

  We next turn to the case that $F$ is non-smooth. As Proposition \ref
  {preliminaries:hyperfiniteness:proposition:smoothavoidance} ensures that there is a
  \Borel reduction of $F$ to the restriction of $F$ to an $F$-invariant \Borel set off of which $F$
  is smooth, by composing such a reduction with $\phi$, we can
  assume that the restriction of $F$ to the set $Z = \setcomplement{\saturation{\image{\phi}{X}}
  {E}}$ is non-smooth. As $\phi$ is countable-to-one, the \Lusin-\Novikov uniformization theorem
  yields an $(\restriction{E}{C})$-complete \Borel set $B \subseteq C$ on which $\phi$ is injective.

  Fix a \Borel set $A \subseteq B$ of maximal $\mu$-measure on which $E$ is compressible.
  As $E$ is countable, the \Lusin-\Novikov uniformization theorem ensures that
  $\saturation{A}{E}$ is \Borel. As Proposition \ref
  {preliminaries:measuredequivalencerelations:proposition:compression} ensures that countable
  \Borel equivalence relations can be \Borel embedded into their restrictions to complete
  compressible \Borel sets, there is a \Borel injection $\psi \from \saturation{A}{E} \to A$
  whose graph is contained in $E$. Then the function $\pi = \phi \composition \psi$ is a \Borel 
  embedding of $\restriction{E}{\saturation{A}{E}}$ into $\restriction{F}{\image{\phi}{C}}$.

  If $\mu(\saturation{A}{E}) = 1$, then it follows that $E$ is $\mu$-embeddable into $F$. Otherwise,
  Theorem \ref{preliminaries:measuredequivalencerelations:theorem:existence} ensures that
  $\restriction{\mu}{(B \setminus \saturation{A}{E})}$ is equivalent to an $\restriction{E}{(B
  \setminus \saturation{A}{E})}$-invariant \Borel probability measure $\nu$ on $B \setminus
  \saturation{A}{E}$. As $E$ is invariant-measure hyperfinite, there is a $\nu$-conull \Borel set 
  $A' \subseteq B \setminus \saturation{A}{E}$ on which $E$ is hyperfinite. As $E$ is countable, the
  \Lusin-\Novikov uniformization theorem ensures that $\saturation{A'}{E}$ is \Borel and
  there is a \Borel reduction of $\restriction{E}{\saturation{A'}{E}}$ to $\restriction{E}{A'}$. As
  Proposition \ref{preliminaries:hyperfiniteness:proposition:closure:homomorphisms} ensures that
  the class of hyperfinite \Borel equivalence relations is closed downward under \Borel
  reducibility, it follows that $\restriction{E}{\saturation{A'}{E}}$ is also hyperfinite. As Theorem
  \ref{preliminaries:hyperfiniteness:theorem:embedding} ensures that every hyperfinite \Borel
  equivalence relation is \Borel embeddable into every non-smooth \Borel equivalence relation,
  there is a \Borel embedding $\pi' \from \saturation{A'}{E} \to Z$ of $\restriction{E}{\saturation{A'}
  {E}}$ into $\restriction{F}{Z}$. As $\mu(\saturation{A \union A'}{E}) = 1$ and $\pi \union \pi'$ is
  an embedding of $\restriction{E}{\saturation{A \union A'}{E}}$ into $F$, the proposition follows.
\end{propositionproof}

This has the following consequence for the relationship between measure embeddability and
measure reducibility.

\begin{proposition} \label{embeddability:proposition:sufficientcondition:measure}
  Suppose that $X$ and $Y$ are standard \Borel spaces, $E$ and $F$ are countable \Borel
  equivalence relations on $X$ and $Y$, $E$ is invariant-measure hyperfinite, and
  $F$ is aperiodic. Then $E$ is measure reducible to $F$ if and only if $E$ is measure
  embeddable into $F$.
\end{proposition}

\begin{propositionproof}
  This is a direct consequence of Proposition \ref
  {embeddability:proposition:sufficientcondition:measurable}.
\end{propositionproof}

In particular, we obtain the following.

\begin{proposition} \label{embeddability:proposition:sufficientcondition:product}
  Suppose that $X$ is a standard \Borel space and $E$ is an aperiodic invariant-measure-hyperfinite
  countable \Borel equivalence relation on $X$. Then $E \times \square{\N}$ is measure embeddable into $E$.
\end{proposition}

\begin{propositionproof}
  This is a direct consequence of Proposition \ref
  {embeddability:proposition:sufficientcondition:measure}.
\end{propositionproof}

We next record a natural obstacle to measurable embeddability. We use $\invariant{X}{E}$ to
denote the family of all $E$-invariant \Borel probability measures on $X$, and
$\ergodicinvariant{X}{E}$ to denote $\ergodic{X}{E} \intersection \invariant{X}{E}$.

\begin{proposition} \label{embeddability:proposition:obstacle:general}
  Suppose that $X$ and $Y$ are standard \Borel spaces, $E$ and $F$ are countable \Borel
  equivalence relations on $X$ and $Y$, $\mu \in \ergodicinvariant{X}{E} \setminus \hyperfinite
  {X}{E}$, $\nu \in \ergodicinvariant{Y}{F} \setminus \hyperfinite{Y}{F}$, $\cost{\mu}{E} < \cost{\nu}
  {F}$, and the \textexponent{\mu}{th} vertical section of $\absolutelycontinuous[E, F]$ is the
  measure-equivalence class of $\nu$. Then $E$ is not $\mu$-embeddable into $F$.
\end{proposition}

\begin{propositionproof}
  Suppose, towards a contradiction, that there is a $\mu$-conull \Borel set $C \subseteq X$ on
  which there is a \Borel embedding $\pi \from C \to Y$ of $E$ into $F$. Then $\pushforward{\pi}
  {(\restriction{\mu}{C})} \absolutelycontinuous \nu$, since otherwise Proposition \ref
  {preliminaries:measuredequivalencerelationsproposition:quasiinvariant} would yield an
  $F$-quasi-invariant \Borel probability measure $\nu'$ on $Y$ with the same $F$-invariant
  \Borel sets as $\pushforward{\pi}{(\restriction{\mu}{C})}$, in which case the $E$-ergodicity of $\mu$ would
  ensure that $\nu'$ is $F$-ergodic, and the downward closure of the family of hyperfinite
  \Borel equivalence relations under \Borel embeddability (see Proposition \ref
  {preliminaries:hyperfiniteness:proposition:closure:homomorphisms}) would imply that $F$ is
  $\nu'$-nowhere hyperfinite, despite the fact that $\nu$ and $\nu'$ are orthogonal. Let $\nu_D$
  be the \Borel probability measure on the set $D = \image{\pi}{C}$ given by $\nu_D(B) = \nu(B)
  / \nu(D)$. As $\pushforward{\pi}{(\restriction{\mu}{C})} \absolutelycontinuous \nu_D$ and both measures are
  $(\restriction{F}{D})$-ergodic and $(\restriction{F}{D})$-invariant, Proposition \ref
  {preliminaries:measuredequivalencerelations:proposition:equality} implies that
  $\pushforward{\pi}{(\restriction{\mu}{C})} = \nu_D$. The formula for the cost of \Borel restrictions given by
  Proposition \ref{preliminaries:measuredequivalencerelations:proposition:costformula} then
  ensures that $\cost{\nu}{F} \le \cost{\nu_D}{\restriction{F}{D}} = \cost{\mu}{E}$, a contradiction.
\end{propositionproof}

As a special case, we obtain the following.

\begin{proposition} \label{embeddability:proposition:obstacle:special}
  Suppose that $X$ is a standard \Borel space, $E$ is a countable \Borel equivalence relation on
  $X$, $\mu \in \ergodicinvariant{X}{E} \setminus \hyperfinite{X}{E}$, $1 < \cost{\mu}{E} < \infty$,
  and the \textexponent{\mu}{th} vertical section of $\absolutelycontinuous[E]$ is the
  measure-equivalence class of $\mu$. Then for no $n \in \N$ is it the case that
  $E \times \square{n+1}$ is $\mu$-embeddable into $E \times \square{n}$.
\end{proposition}

\begin{propositionproof}
  Let $m_n$ denote the uniform probability measure on $n$. Then the formula for the cost of \Borel
  restrictions given by Proposition \ref
  {preliminaries:measuredequivalencerelations:proposition:costformula} ensures that $\cost{\mu
  \times m_{n+1}}{E \times \square{n+1}} < \cost{\mu \times m_n}{E \times \square{n}}$ for all $n \in
  \N$, so Proposition \ref{embeddability:proposition:obstacle:general} implies that $E \times
  \square{n+1}$ is not $(\mu \times m_{n+1})$-embeddable into $E \times \square{n}$.
\end{propositionproof}

Putting these observations together, we obtain the following.

\begin{proposition}
  Suppose that $X$ is a standard \Borel space and $E$ is an aperiodic non-invariant-measure-hyperfinite
  projectively-separable treeable countable \Borel equivalence relation on $X$. Then there is an aperiodic \Borel 
  subequivalence relation $F$ of $E$ such that for no $n \in \N$ is $F \times \square{n+1}$ measure
  embeddable into $F \times \square{n}$.
\end{proposition}

\begin{propositionproof}
  Fix a \Borel set $B \subseteq X$ and an $(\restriction{E}{B})$-invariant \Borel
  probability measure $\mu$ on $B$ such that $\restriction{E}{B}$ is not $\mu$-hyperfinite. Fix a
  \Borel graphing $G$ of $\restriction{E}{B}$. As $G$ is locally countable, the \Lusin-\Novikov
  uniformization theorem yields an increasing sequence $\sequence{G_n}[n
  \in \N]$ of \Borel subgraphs of $G$ of bounded vertex degree whose union is $G$. As Theorem
  \ref{preliminaries:measurehyperfiniteness:theorem:increasingunion} ensures that the increasing
  union of $\mu$-hyperfinite \Borel equivalence relations is $\mu$-hyperfinite, there exists $n \in
  \N$ sufficiently large for which the equivalence relation $F$ generated by $G_n$ is not
  $\mu$-hyperfinite. Note that $\cost{\nu}{F} < \infty$ for every $F$-invariant \Borel probability
  measure $\nu$ on $B$. By Proposition \ref
  {preliminaries:measurehyperfiniteness:proposition:characterization}, there exists $\nu \in
  \ergodicinvariant{B}{F} \setminus \hyperfinite{B}{F}$. As Proposition \ref
  {projectiveseparability:proposition:closure} ensures that the class of projectively-separable
  countable \Borel equivalence relations is closed downward under \Borel restrictions and \Borel
  subequivalence relations, it follows that $F$ is projectively separable. As Proposition \ref
  {measures:proposition:countablesections} ensures that the vertical sections of
  $\absolutelycontinuous[F']$ are countable unions of measure-equivalence classes, there is a
  $\nu$-conull \Borel set $C \subseteq B$ which is null with respect to every measure in the
  \textexponent{\nu}{th} vertical section of
  $\absolutelycontinuous[F']$, with the exception of those in the measure-equivalence class of
  $\nu$. By removing a $\nu$-null \Borel subset of $C$, we can assume that the relation $F' =
  \restriction{F}{C}$ is aperiodic. As Proposition \ref{preliminaries:treeability:proposition:closure}
  ensures that the family of treeable countable \Borel equivalence relations is closed downward
  under \Borel subequivalence relations, it follows that $F'$ is treeable, so $1 < \cost{\nu}{F'} <
  \infty$ by Proposition \ref{preliminaries:measuredequivalencerelations:proposition:costbound}, thus
  Proposition \ref{embeddability:proposition:obstacle:special} implies that for no $n \in \N$ is it the
  case that $F' \times \square{n+1}$ is $\nu$-embeddable into $F' \times \square{n}$.
  Proposition \ref{preliminaries:hyperfiniteness:proposition:aperiodic} then yields an aperiodic 
  hyperfinite \Borel subequivalence relation $F''$ of $\restriction{E}{\setcomplement{C}}$, in which
  case $F' \union F''$ is as desired.
\end{propositionproof}

\section{Antichains} \label{antichains}

In this section, we produce perfect sequences of pairwise non-measure reducible \Borel
subequivalence relations of a given projectively-separable treeable countable \Borel
equivalence relation.

We begin by noting that hyperfiniteness rules out such sequences.

\begin{proposition} \label{antichains:proposition:linear}
  Suppose that $X$ is a standard \Borel space, $E$ is a hyperfinite \Borel equivalence relation on
  $X$, and $E_1$ and $E_2$ are \Borel subequivalence relations of $E$. Then $E_1$ and $E_2$
  are comparable under \Borel reducibility.
\end{proposition}

\begin{propositionproof}
  As Proposition \ref{preliminaries:hyperfiniteness:proposition:closure:homomorphisms} ensures that
  the family of hyperfinite \Borel equivalence relations is closed downward under \Borel
  subequivalence relations, it follows that $E_1$ and $E_2$ are themselves hyperfinite. But
  Theorem \ref{preliminaries:hyperfiniteness:theorem:linearity} implies that any two hyperfinite \Borel
  equivalence relations are comparable under \Borel reducibility.
\end{propositionproof}

We next turn our attention to very special sorts of antichains.

\begin{proposition} \label{antichains:proposition:trichotomy:restriction}
  Suppose that $X$ is a standard \Borel space and $E$ is a non-measure-hyperfinite 
  projectively-separable countable \Borel equivalence relation on $X$. Then exactly one of the
  following holds:
  \begin{enumerate}
    \item The relation $E$ is a non-empty countable disjoint union of successors of $\Ezero$ under
      measure reducibility.
    \item There are \Borel sequences $\sequence{B_c}[c \in \Cantorspace]$ of pairwise disjoint
      $E$-inva\-riant subsets of $X$ and $\sequence{\mu_c}[c \in \Cantorspace]$ of \Borel
      probability measures on $X$ in $\ergodicquasiinvariant{X}{E} \setminus \hyperfinite{X}{E}$
      with the property that $\mu_c(B_c) = 1$ for all $c \in \Cantorspace$, and for no distinct $c, d \in
      \Cantorspace$ is it the case that $\restriction{E}{B_c}$ is $\mu_c$-reducible to $\restriction{E}
      {B_d}$.
  \end{enumerate}
\end{proposition}

\begin{propositionproof}
  By Theorem \ref{measures:theorem:trichotomy}, it is sufficient to show that if $\sequence{B_c}[c
  \in \Cantorspace]$ is a \Borel sequence of pairwise disjoint $E$-invariant sets and $\sequence
  {\mu_c}[c \in \Cantorspace]$ is a \Borel sequence of \Borel probability measures on $X$ in
  $\ergodicquasiinvariant{X}{E} \setminus \hyperfinite{X}{E}$ such that $\mu_c(B_c) = 1$ for all
  $c \in \Cantorspace$, then by passing to a perfect subsequence, one can ensure that for no
  distinct $c, d \in \Cantorspace$ is it the case that $\restriction{E}{B_c}$ is $\mu_c$-reducible to
  $\restriction{E}{B_d}$. Towards this end, let $R$ denote the binary relation on $\Cantorspace$
  in which two sequences $c, d \in \Cantorspace$ are $R$-related if $\restriction{E}{B_c}$ is
  $\mu_c$-reducible to $\restriction{E}{B_d}$. Then Proposition \ref
  {preliminaries:measuredequivalencerelations:proposition:analytic} ensures that $R$ is analytic,
  and therefore has the \Baire property. As the projective separability of $E$ ensures that the vertical
  sections of $R$ are countable, it follows that the vertical sections of $R$ are meager, so the
  \Kuratowski-\Ulam theorem (see, for example, \cite[Theorem 8.41]{Kechris}) ensures that $R$ is
  itself meager, in which case \Mycielski's theorem (see, for example, \cite[Theorem 19.1]{Kechris})
  yields the desired perfect subsequence.
\end{propositionproof}

In particular, this allows us to characterize the circumstances under which there is a perfect
sequence of pairwise non-measure reducible countable \Borel equivalence relations which
are measure reducible to a given projectively-separable countable \Borel equivalence relation.

\begin{proposition} \label{antichains:proposition:trichotomy:reducible}
  Suppose that $X$ is a standard \Borel space and $E$ is a non-measure-hyperfinite
  projectively-separable countable \Borel equivalence relation on $X$. Then exactly one of the
  following holds:
  \begin{enumerate}
    \item There is a finite family $\calF$ of successors of $\Ezero$ under measure reducibility for
      which $E$ is a non-empty countable disjoint union of \Borel equivalence relations
      which are measure bi-reducible with those in $\calF$.
    \item There is a \Borel sequence $\sequence{E_c}[c \in \Cantorspace]$ of pairwise
      non-measure-red\-ucible countable equivalence relations measure reducible to $E$.
  \end{enumerate}
\end{proposition}

\begin{propositionproof}
  In light of Proposition \ref{antichains:proposition:trichotomy:restriction}, we can assume that
  $E$ is a non-empty countable disjoint union of a sequence $\sequence{E_n}[n \in \N]$
  of successors of $\Ezero$ under measure reducibility.
  
  To see that at least one of these conditions holds, note that if condition (1) fails, then by
  passing to an infinite subsequence, we can assume that the relations $E_n$ are
  pairwise non-meas\-ure-reducible. Proposition \ref{measures:proposition:successor:characterization}
  then ensures that if $n \in \N$ and $\mu \in \ergodicquasiinvariant{X}{E_n} \setminus \hyperfinite{X}
  {E_n}$, then $E_n$ is not $\mu$-reducible to $\disjointunion[m \in \N \setminus \set{n}][E_m]$. In
  particular, if $\sequence{N_c}[c \in \Cantorspace]$ is a \Borel sequence  of subsets of $\N$ such
  that $N_c \nsubseteq N_d$ for all distinct $c, d \in \Cantorspace$, then the relations
  $E_c = \disjointunion[n \in N_c][E_n]$ are pairwise non-measure-reducible.
  
  To see that the conditions are mutually exclusive, we will establish the stronger fact that if
  condition (1) holds, then every sequence $\sequence{F_n}[n \in \N]$ of countable \Borel
  equivalence relations measure reducible to $E$ has an infinite subsequence that is
  (not necessarily strictly) increasing under measure reducibility. Towards this end, note that
  for each $n \in \N$, there is a sequence $\sequence{k_{F,n}}[F \in \calF]$ of countable cardinals
  such that $F_n$ is measure bi-reducible with $\disjointunion[F \in \calF][F \times \diagonal{k_{F,n}}]$. 
  A straightforward induction shows that, by passing to an infinite
  subsequence, we can assume that $k_{F, m} \le k_{F, n}$ for all $F \in \calF$ and $m \le n$
  in $\N$. But this implies that $\sequence{F_n}[n \in \N]$ is increasing under measure reducibility.
\end{propositionproof}

As a corollary, we obtain the following.

\begin{proposition}
  Suppose that $X$ is a standard \Borel space and $E$ is a projectively-separable countable
  \Borel equivalence relation on $X$. Then the following are equivalent:
  \begin{enumerate}
    \item There is a sequence $\sequence{E_n}[n \in \N]$ of countable \Borel equivalence
      relations measure reducible to $E$ for which no infinite subsequence is (not necessarily
      strictly) increasing under measure reducibility.
    \item There is a sequence $\sequence{E_n}[n \in \N]$ of pairwise non-measure-reducible
      countable \Borel equivalence relations measure reducible to $E$.
    \item There is a \Borel sequence $\sequence{E_c}[c \in \Cantorspace]$ of pairwise
      non-measure-red\-ucible countable equivalence relations measure reducible to $E$.
  \end{enumerate}
\end{proposition}

\begin{propositionproof}
  This follows from the proof of Proposition \ref{antichains:proposition:trichotomy:reducible}.
\end{propositionproof}

We next turn our attention to subequivalence relations. The main additional tool we will need is
the following observation concerning the power of $\mu$-strat\-ifications in the presence of
projective separability.

\begin{proposition} \label{antichains:proposition:stratification}
  Suppose that $X$ is a standard \Borel space, $E$ is a projectively-separable countable \Borel
  equivalence relation on $X$, $\mu$ is a \Borel probability measure on $X$, $\sequence{B_n}
  [n \in \N]$ is a sequence of $\mu$-positive \Borel subsets of $X$, and $\sequence{E_{n,r}}[r \in
  \R]$ is a \Borel $(\restriction{\mu}{B_n})$-stratification of $\restriction{E}{B_n}$ such that
  $\intersection[r \in \R][E_{n, r}]$ is $(\restriction{\mu}{B_n})$-nowhere hyperfinite, for all $n \in
  \N$. Then there is a \Borel embedding $\pi \from \R \to \R$ of the usual ordering of $\R$
  into itself such that $E_{m, \pi(r)}$ is $(\restriction{\mu}{B_m})$-nowhere reducible to $E_{n, \pi
  (s)}$ for all distinct $\pair{m}{r}, \pair{n}{s} \in \N \times \R$.
\end{proposition}

\begin{propositionproof}
  Let $R_{m,n}$ denote the relation on $\R$ in which two real numbers $r$ and $s$ are related if
  $E_{m, r}$ is $(\restriction{\mu}{B_m})$-somewhere reducible to $E_{n, s}$.

  \begin{lemma}
    Every horizontal section of every $R_{m,n}$ is countable.
  \end{lemma}

  \begin{lemmaproof}
    Suppose, towards a contradiction, that there exist $m, n \in \N$ and $t \in \R$ for which
    $\horizontalsection{R_{m,n}}{t}$ is uncountable. For each $r \in \horizontalsection
    {R_{m,n}}{t}$, fix a $\mu$-positive \Borel set $B_{m,r} \subseteq B_m$ on which there is a
    \Borel reduction $\phi_r \from B_{m,r} \to B_n$ of $E_{m, r}$ to $E_{n, t}$. Then there exists
    $\epsilon > 0$ such that $\mu(B_{m,r}) \ge \epsilon$ for uncountably many $r \in
    \horizontalsection{R_{m,n}}{t}$. As each $\phi_r$ is a homomorphism from $\restriction
    {(\intersection[r \in \R][E_{m, r}])}{B_{m,r}}$ to $E$, the $(\restriction{\mu}{B_m})$-nowhere
    hyperfiniteness of $\restriction{(\intersection[r \in \R][E_{m, r}])}{B_m}$ coupled with the
    projective separability of $E$ ensures the existence of distinct $r, s \in \horizontalsection
    {R_{m,n}}{t}$ for which $\mu(B_{m,r}), \mu(B_{m,s}) \ge \epsilon$ and $\uniformmetric{\mu}
    (\phi_r, \phi_s) < \epsilon$. Then $\set{x \in B_{m,r} \intersection B_{m,s}}[\phi_r(x)
    = \phi_s(x)]$ is a $\mu$-positive \Borel set on which $E_{m,r}$ and $E_{m,s}$ coincide, a
    contradiction.
  \end{lemmaproof}

  Proposition \ref{preliminaries:measuredequivalencerelations:proposition:analytic} ensures that each
  $R_{m,n}$ is analytic, and therefore has the \Baire property. As the horizontal sections of each $R_{m,n}$ are
  countable and therefore meager, the \Kuratowski-\Ulam theorem ensures that each $R_{m,n}$ is
  meager, thus so too is their union $R$, in which case Mycielski's theorem yields a continuous injection $\phi \from
  \Cantorspace \to \R$ with respect to which pairs of distinct sequences in $\Cantorspace$ are mapped to
  $R$-unrelated pairs of real numbers. \Galvin's theorem (see, for example, \cite[Theorem 19.7]{Kechris})
  ensures that by replacing $\phi$ with its composition with an appropriate continuous function from
  $\Cantorspace$ to $\Cantorspace$, we can assume that it is an embedding of the lexicographical
  ordering of $\Cantorspace$ into the usual ordering of $\R$. Fix a \Borel embedding $\psi \from \R \to
  \Cantorspace$ of the usual ordering of $\R$ into the lexicographical ordering of $\Cantorspace$, and
  observe that the function $\pi = \phi \composition \psi$ is as desired.
\end{propositionproof}

In particular, this yields the following measure-theoretic result.

\begin{theorem} \label{antichains:theorem:dichotomy:measurable}
  Suppose that $X$ is a standard \Borel space, $E$ is a projectively-separable treeable
  countable \Borel equivalence relation on $X$, and $\mu$ is a \Borel probability measure on $X$
  for which $E$ is $\mu$-nowhere hyperfinite. Then there is an increasing \Borel sequence
  $\sequence{E_r}[r \in \R]$ of pairwise $\mu$-nowhere reducible subequivalence relations of $E$.
\end{theorem}

\begin{theoremproof}
  As Proposition \ref{preliminaries:measuredequivalencerelationsproposition:quasiinvariant} yields an
  $E$-quasi-invariant \Borel probability measure $\nu$ for which $\mu \absolutelycontinuous \nu$,
  Theorem \ref{stratification:theorem:sufficientcondition} yields a \Borel $\mu$-stratification
  $\sequence{F_r}[r \in \R]$ of $E$. As Theorem \ref
  {preliminaries:measurehyperfiniteness:theorem:increasingunion} ensures that the family of
  $\mu$-hyperfinite countable \Borel equivalence relations is closed under increasing unions,
  there is a partition $\sequence{B_n}[n \in \N]$ of $X$ into $\mu$-positive \Borel sets, as well as
  a sequence $\sequence{r_n}[n \in \N]$ of real numbers, such that $\restriction{F_{r_n}}{B_n}$ is
  $(\restriction{\mu}{B_n})$-nowhere hyperfinite for all $n \in \N$. Fix order-preserving \Borel
  injections $\phi_n \from \R \to \openinterval{r_n}{\infty}$, and appeal to Proposition \ref
  {antichains:proposition:stratification} to obtain a \Borel embedding $\phi \from \R \to \R$ of the
  usual ordering of $\R$ into itself such that $\restriction{F_{(\phi_m \composition \phi)(r)}}{B_m}$
  is $(\restriction{\mu}{B_m})$-nowhere reducible to $F_{(\phi_n \composition \phi)(s)}$ for all
  distinct $\pair{m}{r}, \pair{n}{s} \in \N \times \R$. Then the relations $E_r = \union[n \in \N]
  [(\restriction{F_{(\phi_n \composition \phi)(r)}}{B_n})]$ are as desired.
\end{theoremproof}

In the special case that the equivalence relation in question is a successor of $\Ezero$ under
measure reducibility, we can ensure that the same holds of the subequivalence relations.

\begin{theorem} \label{antichains:theorem:dichotomy:measurable:successor}
  Suppose that $X$ is a standard \Borel space, $E$ is a projectively-separable treeable countable
  \Borel equivalence relation on $X$ which is a successor of $\Ezero$ under measure reducibility,
  and $\mu$ is a \Borel probability measure on $X$ for which $E$ is $\mu$-nowhere hyperfinite.
  Then there is an increasing \Borel sequence $\sequence{E_r}[r \in \R]$ of pairwise
  $\mu$-nowhere reducible subequivalence relations of $E$ consisting of successors of
  $\Ezero$ under measure reducibility with the property that $\mu$ is $(\intersection[r \in \R]
  [E_r])$-ergodic.
\end{theorem}

\begin{theoremproof}
  Proposition \ref{preliminaries:measuredequivalencerelationsproposition:quasiinvariant} yields an
  $E$-quasi-invariant \Borel probability measure $\nu \reverseabsolutelycontinuous \mu$ agreeing
  with $\mu$ on all $E$-invariant \Borel sets. By Proposition \ref
  {preliminaries:measuredequivalencerelationsproposition:cocycles}, there is a \Borel cocycle $\rho
  \from E \to \Rplus$ with respect to which $\nu$ is invariant. As Theorem \ref
  {preliminaries:measuredequivalencerelations:theorem:ergodicdecomposition} ensures the
  existence of a \Borel ergodic decomposition of $\rho$, Proposition \ref
  {measures:proposition:successor:characterization} implies that
  $\ergodicquasiinvariant{X}{E} \setminus \hyperfinite{X}{E}$ consists of a
  single measure-equivalence class, and Proposition \ref
  {preliminaries:measuredequivalencerelations:proposition:equality} implies that $E$ is not
  almost-everywhere hyperfinite with respect to at most one measure along the ergodic decomposition,
  it follows from Proposition \ref{preliminaries:measurehyperfiniteness:proposition:smooth} that
  $\nu$ is $E$-ergodic. Proposition \ref{measures:proposition:E0ergodic} therefore implies that
  $\nu$ is $\pair{E}{\Ezero}$-ergodic.

  Theorem \ref{stratification:theorem:sufficientcondition} yields a \Borel $\nu$-stratification
  $\sequence{F_r}[r \in \R]$ of $E$. Theorem \ref
  {preliminaries:measurehyperfiniteness:theorem:increasingunion} ensures that not every $F_r$ is
  $\nu$-hyperfinite, so by passing to a \Borel subsequence, we can assume that there is a
  $\nu$-positive \Borel set $B \subseteq X$ on which $\intersection[r \in \R][F_r]$ is
  $\nu$-nowhere hyperfinite. Proposition \ref
  {preliminaries:measuredequivalencerelations:proposition:E0ergodic} implies that by passing to a
  further subsequence, we can also assume that $\restriction{\nu}{B}$ is $(\intersection[r \in \R]
  [\restriction{F_r}{B}])$-ergodic. Proposition \ref{antichains:proposition:stratification} therefore
  yields a \Borel embedding $\phi \from \R \to \R$ of the usual ordering of $\R$ into itself such that
  $\restriction{F_{\phi(r)}}{B}$ is $(\restriction{\nu}{B})$-nowhere reducible to $F_{\phi(s)}$ for all
  distinct $r, s \in \R$. As $E$ is countable, the \Lusin-\Novikov uniformization
  theorem ensures that the set $\saturation{B}{E}$ is \Borel,
  and that there is an extension of the identity function on $B$ to a \Borel function $\psi \from
  \saturation{B}{E} \to B$ whose graph is contained in $E$. Let $E_r$ denote the equivalence
  relation given by $x \mathrel{E_r} y \iff \psi(x) \mathrel{F_{\phi(r)}} \psi(y)$ on $\saturation{B}{E}$,
  and which is trivial off of $\saturation{B}{E}$. As Proposition \ref{measures:proposition:containment}
  ensures that each $\ergodicquasiinvariant{X}{E_r} \setminus \hyperfinite{X}{E_r}$ consists of a single
  measure-equivalence class, Proposition \ref{measures:proposition:successor:sufficientcondition}
  implies that each $E_r$ is a successor of $\Ezero$ under measure reducibility.
\end{theoremproof}

We close this section with the \Borel analogs of these results.

\begin{theorem} \label{antichains:theorem:dichotomy:measure}
  Suppose that $X$ is a standard \Borel space and $E$ is a non-measure-hyperfinite 
  projectively-separable treeable countable \Borel equivalence relation on $X$. Then
  there is an increasing \Borel sequence $\sequence{E_r}[r \in \R]$ of pairwise
  non-measure-reducible subequivalence relations of $E$.
\end{theorem}

\begin{theoremproof}
  Appeal to Theorem \ref{preliminaries:measurehyperfiniteness:theorem:E0} to obtain a \Borel
  probability measure $\mu \in \ergodicquasiinvariant{X}{E} \setminus \hyperfinite{X}{E}$, and
  apply Theorem \ref{antichains:theorem:dichotomy:measurable}.
\end{theoremproof}

\begin{theorem} \label{antichains:theorem:dichotomy:measure:successor}
  Suppose that $X$ is a standard \Borel space and $E$ is a projectively-separable treeable countable
  \Borel equivalence relation on $X$ which is a successor of $\Ezero$ under measure reducibility.
  Then there is an increasing \Borel sequence $\sequence{E_r}[r \in \R]$ of pairwise
  non-measure-reducible subequivalence relations of $E$ which are themselves
  successors of $\Ezero$ under measure reducibility.
\end{theorem}

\begin{theoremproof}
  Appeal to Theorem \ref{preliminaries:measurehyperfiniteness:theorem:E0} to obtain a \Borel
  probability measure $\mu \in \ergodicquasiinvariant{X}{E} \setminus \hyperfinite{X}{E}$, and
  apply Theorem \ref{antichains:theorem:dichotomy:measurable:successor}.
\end{theoremproof}

\section{Bases} \label{bases}
Here we establish the nonexistence of small bases $\calB \subseteq \calE$ for $\calE$
under measure reducibility. We obtain the optimal result in this direction when working below
successors of $\Ezero$ under measure reducibility.

\begin{theorem} \label{bases:theorem:successor}
  Suppose that $X$ is a standard \Borel space and $E$ is a projectively-separable treeable countable
  \Borel equivalence relation on $X$ that is a successor of $\Ezero$ under measure reducibility.
  Then every basis for the non-measure-hyperfinite \Borel subequivalence relations of $E$ has
  cardinality at least $\continuum$.
\end{theorem}

\begin{theoremproof}
  By Theorem \ref{antichains:theorem:dichotomy:measure:successor}, there is an increasing \Borel
  sequence $\sequence{E_r}[r \in \R]$ of pairwise non-measure-reducible subequivalence relations
  of $E$, which are also successors of $\Ezero$ under measure reducibility. Then each element of
  $\calB$ is measure reducible to at most one $E_r$, thus $\cardinality{\calB} \ge \continuum$.
\end{theoremproof}

While we can nearly obtain the analogous result without the assumption that $E$ is a successor
of $\Ezero$ under measure reducibility, there is a slight metamathematical wrinkle. Although we have thus
far freely used the axiom of choice throughout the paper, it is not difficult to push through all of our
arguments under the axiom of dependent choice. While the cardinality restriction appearing
below implies only that bases are necessarily uncountable under the axiom of
dependent choice, it yields the full result that bases have size continuum under the
axiom of choice, as well as in models of the axiom of dependent choice where every subset
of the real numbers has the \Baire property and there is an injection of the real numbers into every
non-well-orderable set, such as $L(\R)$ under the axiom of determinacy (see
\cite{CaicedoKetchersid}).

\begin{theorem}
  Suppose that $X$ is a standard \Borel space, $E$ is a non-measure-hyperfinite
  projectively-separable treeable countable \Borel equivalence relation on $X$, and
  $\calB$ is a basis for the non-measure-hyperfinite \Borel subequivalence relations
  of $E$ under measure reducibility. Then $\R$ is a union of $\cardinality
  {\calB}$-many countable sets.
\end{theorem}

\begin{theoremproof}
  By Theorem \ref{antichains:theorem:dichotomy:measure}, there is an increasing \Borel
  sequence $\sequence{E_r}[r \in \R]$ of pairwise non-measure-reducible subequivalence
  relations of $E$. But then each element of $\calB$ is measure reducible to only
  countably-many relations of the form $E_r$.
\end{theoremproof}

\section{Complexity} \label{complexity}

In this section, we establish a technical strengthening of Theorem \ref
{antichains:theorem:dichotomy:measurable} which gives rise to our complexity results.

\begin{theorem} \label{complexity:theorem:strongantichain}
  Suppose that $X$ is a standard \Borel space and $E$ is a non-measure-hyperfinite
  projectively-separable treeable countable \Borel equivalence relation on $X$. Then
  there are \Borel sequences $\sequence{E_r}[r \in \R]$ of subequivalence relations of
  $E$ and $\sequence{\mu_r}[r \in \R]$ of \Borel probability measures on $X$ such that:
  \begin{enumerate}
    \item Each $\mu_r$ is $E_r$-quasi-invariant and $E_r$-ergodic.
    \item The relation $E_r$ is $\mu_r$-nowhere reducible to the relation $E_s$, for all distinct
      $r, s \in \R$.
  \end{enumerate}
\end{theorem}

\begin{theoremproof}
  Note that if $E$ is not a countable disjoint union of successors
  of $\Ezero$ under measure reducibility, then Proposition \ref{antichains:proposition:trichotomy:restriction}
  yields the desired result. On the other hand, if $E$ is a countable disjoint union of successors of
  $\Ezero$ under measure reducibility, then there is an $E$-invariant \Borel set $B \subseteq X$
  on which $E$ is a successor of $\Ezero$ under measure reducibility. Proposition \ref
  {measures:proposition:successor:characterization} then yields a \Borel probability measure
  $\mu$ on $B$ for which $\restriction{E}{B}$ is $\mu$-nowhere hyperfinite, in which case one
  obtains the desired equivalence relations by trivially extending those given by Theorem \ref
  {antichains:theorem:dichotomy:measurable:successor} from $B$ to $X$.
\end{theoremproof}

As a consequence, we obtain the following.

\begin{theorem} \label{theorem:complexity}
  Suppose that $X$ is a standard \Borel space and $E$ is a non-measure-hyperfinite 
  projectively-separable treeable countable \Borel equivalence relation on $X$. Then
  the following hold:
  \begin{enumerate}
    \renewcommand{\theenumi}{\alph{enumi}}
    \item There is an embedding of containment on \Borel subsets of $\R$ into \Borel
      reducibility of countable \Borel equivalence relations with smooth-to-one \Borel
      homomorphisms to $E$ (in the codes).
    \item \Borel bi-reducibility and reducibility of countable \Borel equival\-ence relations with
      smooth-to-one \Borel homomorphisms to $E$ are both $\Sigmaclass[1][2]$-complete
      (in the codes).
    \item Every \Borel quasi-order is \Borel reducible to \Borel reducib\-ility of countable \Borel
      equivalence relations with smooth-to-one \Borel homomorphisms to $E$.
    \item \Borel and $\sigmaclass{\Sigmaclass[1][1]}$-measurable reducibility do not agree on the
      countable \Borel equivalence relations with smooth-to-one \Borel homomorphisms to $E$.
  \end{enumerate}
\end{theorem}

\begin{theoremproof}
  By Theorem \ref{complexity:theorem:strongantichain}, there are \Borel sequences
  $\sequence{E_r}[r \in \R]$ of subequivalence relations of $E$ and $\sequence{\mu_r}[r \in
  \R]$ of \Borel probability measures on $X$ such that:
    \begin{enumerate}
      \item Each $\mu_r$ is $E_r$-quasi-invariant and $E_r$-ergodic.
      \item The relation $E_r$ is $\mu_r$-nowhere reducible to the relation $E_s$, for all distinct
        $r, s \in \R$.
    \end{enumerate}
  But then Theorem \ref{preliminaries:complexity:theorem:antichain} yields the desired result.
\end{theoremproof}

\begin{acknowledgments}
  We would like to thank Manuel Inselmann, Alex\-ander Kechris, Andrew Marks, and the
  anonymous referees for their valuable comments on earlier versions of this paper.
\end{acknowledgments}

\addblanklinetocontents
\bibliographystyle{amsalpha}
\bibliography{bibliography}

\end{document}